\documentclass[12pt]{article}
\usepackage{epsf, epsfig, color}
\usepackage{amssymb,latexsym}

\newtheorem{lem}{Lemma}[section]
\newtheorem{theo}{Theorem}[section]

\newtheorem{cor}{Corollary}[section]
\newtheorem{con}{Conjecture}[section]

\newcommand {\red} {\textcolor{red}}
\newcommand {\green} {\textcolor{green}}
\newcommand {\blue} {\textcolor{blue}}

\def\relabel {\label}  %final version

\renewcommand{\theenumi}{\rm (\roman{enumi})}

\newcommand{\proof}
{{\noindent {\em Proof}.\quad}\setcounter{countclaim}{0}
\setcounter{countcase}{0}}
\newcommand{\proofend}{{\hfill$\Box$}}
\newcommand{\proofends}{\mbox {\hfill    } \Box}
\newcommand{\gap}{\vspace{0.5 cm}}
\newcommand{\sgap}{\vspace{0.3 cm}}

\newcommand{\rem}{\noindent {\bf Remark.}}
\newcommand{\rems}{\noindent {\bf Remarks.}}
\newcommand{\note}{\noindent {\bf Note.}}
\newcommand{\notes}{\noindent {\bf Notes.}}

\newcommand{\acircle}{\makebox(0,0)[l]{\circle*{8}}}

\newcounter{countfig}
\def\infig{\addtocounter{countfig}{1}{Figure \thecountfig}}
\def\itfig{{\rm {\it Figure} \thecountfig}}
\def\norfig{{\rm Figure \thecountfig}}

\newcounter{countclaim}
\def\inclaim{\addtocounter{countclaim}{1}
{\noindent {\bf Claim \thecountclaim}: }}

\newcounter{countcase}
\def\incase{\addtocounter{countcase}{1}
{\noindent {\bf Case \thecountcase}: }}

\newcommand{\beeq}{\begin{equation}}
\newcommand{\eneq}{\end{equation}}

\newcommand{\beeqn}{\begin{eqnarray*}}
\newcommand{\eneqn}{\end{eqnarray*}}

\def \sets{{\cal S}}
\def\setz{{\cal Z}}
\def \setr{{\cal R}}

\def \remarkrevision {\red{checking******}\\}

\begin{document}
  
\newcommand{\resection}[1]
{\section{#1}\setcounter{equation}{0}}

\renewcommand{\theequation}{\thesection.\arabic{equation}}

\renewcommand{\labelenumi}{\rm(\roman{enumi})}

\baselineskip 0.6 cm

\title {On Zero-free Intervals of Flow Polynomials\thanks{Partially 
supported by NIE AcRf funding (RI 2/12 DFM) of Singapore.}}

\author
{ 
F.M. Dong\thanks{Corresponding author.
Email: fengming.dong@nie.edu.sg.} \\
\small Mathematics and Mathematics Education\\
\small National Institute of Education\\
\small Nanyang Technological University, Singapore 637616
\sgap \\
%K.M. Koh\\
%\small Department of Mathematics\\
%\small National University of Singapore,Singapore 117543
}

\date{}

\maketitle

\begin{abstract}
This article studies real roots of 
the flow polynomial $F(G,\lambda)$ of a bridgeless graph $G$. 
For any integer $k\ge 0$, 
let $\xi_k$ be the supremum in $(1,2]$
such that $F(G,\lambda)$ has no real roots in 
$(1,\xi_k)$ for all graphs $G$ with $|W(G)|\le k$,
where $W(G)$ is the set of vertices in $G$ of degrees larger than $3$. 
We prove that $\xi_k$ can be determined by considering 
a finite set of graphs and show that 
$\xi_k=2$ for $k\le 2$, $\xi_3=1.430\cdots$,
$\xi_4=1.361\cdots$ and $\xi_5=1.317\cdots$.
We also prove that for any bridgeless graph $G=(V,E)$, 
if all roots of $F(G,\lambda)$ are real but
some of these roots are not in the set $\{1,2,3\}$, then 
$|E|\ge  |V|+17$ and  
$F(G,\lambda)$ has at least 9 real roots in $(1,2)$.
%then every root of of $F(G,\lambda)$ is contained in the set $\{1,2,3\}$.
\end{abstract}


\section{Introduction}

The graphs considered in this paper 
are undirected and finite,
and may have loops and parallel edges. 
However, the graphs should have no loops
when their chromatic polynomials are considered,
and the graphs should have no bridges when 
their flow polynomials are considered. 
For any graph $G$, let $V(G), E(G), P(G,\lambda)$ 
and $F(G,\lambda)$ be the set of vertices,
the set of edges, the chromatic polynomial 
and the flow polynomial of $G$.
The roots of $P(G,\lambda)$ and $F(G,\lambda)$ 
are called {\it the chromatic roots}
and {\it the flow roots} of $G$ respectively.

A {\it near-triangulation} is 
a loopless connected plane graph in
which at most one face is not bounded by a cycle of order 3. 
Birkhoff and Lewis~\cite{bir} showed that 
$G$ has no real chromatic roots in $(1,2)$
for every near-triangulation $G$. 
Since $P(G,\lambda)=\lambda F(G^*,\lambda)$
for any plane graph $G$, where $G^*$ is its dual, 
this result is equivalent to 
that any connected plane graph $G$ has no flow roots in $(1,2)$
under the condition $|W(G)|\le 1$, 
where $W(G)$ is the set of vertices $x$ in $G$ 
with its degree\footnote{The degree of $x$ in $G$, 
denoted by $d_G(x)$ (or simply $d(x)$), 
is defined to be the sum of the number of non-loop edges in $G$ incident with $x$ and twice the number of loops in $G$ incident with $x$.} larger than 3.

Jackson~\cite{jac3} generalized Birkhoff and Lewis' result 
by showing that any bridgeless connected graph $G$ with $|W(G)|\le 1$
has no real flow roots in $(1,2)$,
no matter whether $G$ is planar or non-planar. 

One of the purposes of this paper is to find 
maximal zero-free intervals 
in $(1,2)$ for the flow polynomials of some families of graphs
and hence extend Jackson's result mentioned above.
For any integer $k\ge 0$, 
let $\Psi_k$ be the 
set of bridgeless connected graphs with $|W(G)|\le k$
and $\xi_k$ be the supremum  in $(1,2]$ such that 
every graph $G$ in $\Psi_k$ has no flow roots 
%(i.e., the roots of $F(G,\lambda)$) 
in $(1,\xi_k)$. 
So $\xi_0, \xi_1, \xi_2,\cdots$ is a non-increasing sequence.  
In Section~\ref{sec4}, 
we will show that $\xi_k$ can be determined 
by considering a finite set of graphs in $\Psi_k\setminus \Psi_{k-1}$ 
and finds that $\xi_k=2$ for $k\le 2$,  
$\xi_3=1.430\cdots$, $\xi_4=1.361\cdots$ and $\xi_5=1.317\cdots$.

By definition,  the flow polynomial $F(G,\lambda)$ is $0$ if $G$ contains 
a bridge (e.g., see (\ref{eq1-2})).
A graph $G=(V,E)$ is said to be 
{\it non-separable} if $G$ is  connected with no 
cut-vertex\footnote{A vertex $x$ in $G$ is called a 
{\it cut-vertex} if $G-x$ has more components that $G$ has.} 
and either it has no loops or $|E|=|V|=1$. 

By the definition, 
a graph with one vertex and at most one edge is non-separable,
and a non-separable graph has a bridge if and only if 
this graph is $K_2$.
A graph %$G$ with $|V(G)|\ge 2$ 
is called {\it separable} 
if it is not non-separable.
A {\it block} of $G$ is a maximal subgraph of $G$ 
with the property that it is non-separable.
%Note that $K_2$ is non-separable and a non-separable graph has a bridge if and only if this graph is $K_2$.
By Lemma~\ref{block-factor}, 
if a graph $G$ is separable, then $F(G,\lambda)$ is the product of 
$F(B,\lambda)$ over all blocks $B$ of $G$.
By Lemmas~\ref{v-edge} and~\ref{3-edge}, 
for a non-separable graph $G$, 
if either $G-e$ \footnote{$G-e$ is the subgraph of $G$ obtained from $G$ by deleting $e$.} 
is separable for some edge $e$ in $G$ 
or $G$ has a proper 3-edge-cut\footnote{A 3-edge-cut $E'$ of $G$ 
is said to be {\it proper}
if the deletion of all edges in $E'$ 
produces more non-empty components than $G$ has.
Thus, if $G$ is non-separable, then a 3-edge-cut of $G$ is proper 
if and only if this 3-edge-cut is not formed by 
three edges incident with a common vertex of degree $3$.}, 
then $(\lambda-1)F(G,\lambda)$ or $(\lambda-1)(\lambda-2)F(G,\lambda)$
is equal to the product of 
the flow polynomials of two graphs with less edges.
Note that if $G$ has a 2-edge-cut, then 
$G-e$ is separable for each $e$ in this cut.
Thus, when we consider the locations of flow roots, we need only 
to study those non-separable graphs 
which contain no proper 3-edge-cut nor an edge $e$ 
with $G-e$ to be separable. 

Another purpose of this paper is to study the existence of bridgeless graphs 
which have real flow roots only but have some flow roots 
not in the set  $\{1,2,3\}$. 
If such graphs do exist, then 
some of them are 
%there are such graphs which are also 
non-separable graphs which have neither 
$2$-edge-cut nor proper $3$-edge-cut.
In Section~\ref{sec6}, we show that 
if a non-separable graph $G=(V,E)$ is such a graph
and contains neither 
$2$-edge-cut nor proper $3$-edge-cut, 
%if all flow roots of $G$ are real and some flow rootof $G$ is not contained in $\{1,2,3\}$, 
then $G$ will satisfy various conditions 
(see Theorem~\ref{main-lem1}), 
including that $|W(G)|\ge 3$, 
$|E(G)|\ge |V(G)|+8|W(G)|-7$ and 
$G$ has at least $\frac {22}{27}(2|W(G)|-1)$
%$2.45|W(G)|-1.23$ %$9$ 
real roots in $(1,2)$. 
In the end of this paper, we pose a conjecture that 
that for any bridgeless graph $G$, 
if all flow roots of $G$ are real, 
then every flow root of $G$ is in the set $\{1,2,3\}$.

\resection{Some fundamental results on flow polynomials}
\relabel{sec2}

\def \tvcu {{\cal G}_2}
\def \sets{{\cal S}}
\def \sett{{\cal T}}

The {\it flow polynomial} $F(G,\lambda)$ of a graph $G$
can be obtained from the following  
properties of $F(G,\lambda)$ (see Tutte~\cite{tut}):
\begin{equation}\relabel{eq1-2}
F(G,\lambda)=
\left \{
\begin{array}{ll}
1, &\mbox{if }E=\emptyset;\\
0, &\mbox{if }G\mbox{ has a bridge};\\
F(G_1,\lambda)F(G_2,\lambda),
&\mbox{if }G=G_1\cup G_2;\\
(\lambda-1)F(G-e,\lambda), &\mbox{if }e\mbox{ is a loop};\\
F(G / e,\lambda)-F(G-e,\lambda),  \quad
&\mbox{otherwise},
\end{array}
\right.
\end{equation}
where $G/e$ \footnote{
If $u$ and $v$ are two vertices of a graph $H$, 
let $H/uv$ denote the graph obtained from $H$ by identifying $u$ and $v$.
So every edge of $H$ is also an edge in $H/uv$ 
and every edge of $H$ joining $u$ and $v$ becomes a loop in $H/uv$.
Then $G/e$ is the graph $(G-e)/uv$, where $u$ and $v$ are the two ends of $e$.}  is the graphs obtained from 
$G$ by contracting $e$ respectively,
and $G_1\cup G_2$ is the disjoint union of graphs 
$G_1$ and $G_2$.

By definition, a loop in $G$ is considered as a block,
and any block with more than one vertex 
has no loops nor cut-vertices. 
Let $b(G)$ be the number of non-trivial blocks 
(i.e., those blocks which are not $K_1$) of $G$.
Thus $b(G)=0$ if and only if $E(G)=\emptyset$,
and if $G$ is connected with $E(G)\ne \emptyset$, then 
$b(G)=1$ if and only if $G$ is non-separable.

For a connected graph $G=(V,E)$ without loops, 
it is well known (see Woodall~\cite{woo}) that 
$(-1)^{|V|}P(G,\lambda)>0$ for all real $\lambda<0$
and
$(-1)^{|V|-1}P(G,\lambda)>0$ for all real $0<\lambda<1$.
Woodall~\cite{woo} and Whitehead and Zhao~\cite{white}
independently showed that 
 $G$ always has a chromatic root of 
multiplicity $b(G)$ at $\lambda =1$.
Jackson~\cite{jac1} also proved that 
$(-1)^{|V|-b(G)+1}P(G,\lambda)>0$ 
for all real $1<\lambda\le 32/27$, 
where the result does not hold if $32/27$ is replaced by any larger 
number.
For flow polynomials, there is 
an analogous result due to Wakelin~\cite{wak}.

\begin{theo}[\cite{wak}]\quad
\relabel{Wakelin}
Let $G=(V,E)$ be a bridgeless connected graph. Then
\begin{enumerate}
\item[(a)] 
$F(G, \lambda)$ is non-zero with sign 
$(-1)^{|E|-|V|+1}$ for $\lambda\in (-\infty, 1)$;
\item[(b)]  
$F(G, \lambda)$ has a zero of multiplicity 
$b(G)$ at $\lambda = 1$;
\item[(c)]  
$F(G, \lambda)$ is non-zero with sign 
$(-1)^{|E|-|V|+b(G)-1}$ for $\lambda\in (1, 32/27]$.
\proofend
\end{enumerate}
\end{theo}

In this paper,
the properties of factorization of flow polynomials
will be applied repeatedly. 
By the result in (\ref{eq1-2}),
the following result can be easily proved by induction. 

\begin{lem}\quad
\relabel{block-factor}
Let $G$ be a bridgeless graph. 
If %$G_1, G_2, \cdots, G_k$ are the components of $G$ or 
%$G$ is connected and 
$G_1, G_2, \cdots, G_k$ are the blocks of $G$, then 
\beeq
F(G,\lambda)=\prod_{1\le i\le k}F(G_i,\lambda).
\eneq
\end{lem}

The next three results on the factorization of flow polynomials
 can be found in~\cite{jac3}
(see \cite{jac2,jac4} also). 
For any graph $G$ and any two vertices $u$ and $v$ in $G$,
let $G+uv$ denote the graph obtained by adding a new edge 
joining $u$ and $v$.

\begin{lem}[\cite{jac3}]\quad
\relabel{v-edge}
Let $G$ be a bridgeless connected graph, $v$ be a vertex of $G$, 
$e = u_1u_2$ be an edge of $G$, and 
$H_1$ and $H_2$ be edge-disjoint subgraphs of 
$G$ such that $E(H_1) \cup E(H_2) = E(G - e)$, 
$V(H_1)\cap V(H_2) = \{v\}$, 
$V(H_1)\cup V(H_2) =V(G)$,
$u_1\in V (H_1)$ and $u_2\in V(H_2)$, as shown in Figure~\ref{f1}. 
Then 
\beeq
F(G, \lambda) =
\frac{F(G_1, \lambda)F(G_2, \lambda)}{\lambda -1}.
\eneq
where $G_i=H_i+vu_i$  %be obtained from $G$ by contracting $E(H_{3-i})$, 
for $i\in \{1, 2\}$.
\end{lem}

\begin{figure}[h!]
\centering 
%\scalebox{1.0}{\input{f1.pic}}
%\input f1.pic
%TeXCAD (http://texcad.sf.net/) Picture. File: [f1.pic]. Options on following lines.
%\grade{\on}
%\emlines{\off}
%\epic{\off}
%\beziermacro{\on}
%\reduce{\on}
%\snapping{\off}
%\pvinsert{% Your \input, \def, etc. here}
%\quality{8.000}
%\graddiff{0.005}
%\snapasp{1}
%\zoom{4.0000}
\unitlength 1mm % = 2.845pt
\linethickness{0.4pt}
\ifx\plotpoint\undefined\newsavebox{\plotpoint}\fi % GNUPLOT compatibility
\begin{picture}(141.75,48)(0,0)
\qbezier(38.5,43.75)(1.5,40.75)(.5,11.75)
\qbezier(137.25,43.5)(100.25,40.5)(99.25,11.5)
\qbezier(.5,11.75)(31,21)(38.5,43.25)
\qbezier(99.25,11.5)(129.75,20.75)(137.25,43)
\qbezier(38.75,43.5)(70,42.875)(69.25,11.75)
\qbezier(69.25,11.75)(45.875,15.75)(39,42.75)
\put(6.25,16){\line(1,0){57.5}}
\put(39,43){\circle*{3.606}}
\put(137.75,42.75){\circle*{3.606}}
\put(6.25,16){\circle*{3.354}}
\put(105,15.75){\circle*{3.354}}
\put(63.5,16){\circle*{3.041}}
%\dashline{1}(46,40.5)(65.75,27.25)
\multiput(45.93,40.43)(.049375,-.033125){16}{\line(1,0){.049375}}
\multiput(47.51,39.37)(.049375,-.033125){16}{\line(1,0){.049375}}
\multiput(49.09,38.31)(.049375,-.033125){16}{\line(1,0){.049375}}
\multiput(50.67,37.25)(.049375,-.033125){16}{\line(1,0){.049375}}
\multiput(52.25,36.19)(.049375,-.033125){16}{\line(1,0){.049375}}
\multiput(53.83,35.13)(.049375,-.033125){16}{\line(1,0){.049375}}
\multiput(55.41,34.07)(.049375,-.033125){16}{\line(1,0){.049375}}
\multiput(56.99,33.01)(.049375,-.033125){16}{\line(1,0){.049375}}
\multiput(58.57,31.95)(.049375,-.033125){16}{\line(1,0){.049375}}
\multiput(60.15,30.89)(.049375,-.033125){16}{\line(1,0){.049375}}
\multiput(61.73,29.83)(.049375,-.033125){16}{\line(1,0){.049375}}
\multiput(63.31,28.77)(.049375,-.033125){16}{\line(1,0){.049375}}
\multiput(64.89,27.71)(.049375,-.033125){16}{\line(1,0){.049375}}
%\end
%\dashline{1}(43.75,38.5)(66.5,23.25)
\multiput(43.68,38.43)(.0477941,-.0320378){17}{\line(1,0){.0477941}}
\multiput(45.305,37.34)(.0477941,-.0320378){17}{\line(1,0){.0477941}}
\multiput(46.93,36.251)(.0477941,-.0320378){17}{\line(1,0){.0477941}}
\multiput(48.555,35.162)(.0477941,-.0320378){17}{\line(1,0){.0477941}}
\multiput(50.18,34.073)(.0477941,-.0320378){17}{\line(1,0){.0477941}}
\multiput(51.805,32.983)(.0477941,-.0320378){17}{\line(1,0){.0477941}}
\multiput(53.43,31.894)(.0477941,-.0320378){17}{\line(1,0){.0477941}}
\multiput(55.055,30.805)(.0477941,-.0320378){17}{\line(1,0){.0477941}}
\multiput(56.68,29.715)(.0477941,-.0320378){17}{\line(1,0){.0477941}}
\multiput(58.305,28.626)(.0477941,-.0320378){17}{\line(1,0){.0477941}}
\multiput(59.93,27.537)(.0477941,-.0320378){17}{\line(1,0){.0477941}}
\multiput(61.555,26.448)(.0477941,-.0320378){17}{\line(1,0){.0477941}}
\multiput(63.18,25.358)(.0477941,-.0320378){17}{\line(1,0){.0477941}}
\multiput(64.805,24.269)(.0477941,-.0320378){17}{\line(1,0){.0477941}}
%\end
%\dashline{1}(44,34.25)(66.5,19.25)
\multiput(43.93,34.18)(.0502232,-.0334821){16}{\line(1,0){.0502232}}
\multiput(45.537,33.108)(.0502232,-.0334821){16}{\line(1,0){.0502232}}
\multiput(47.144,32.037)(.0502232,-.0334821){16}{\line(1,0){.0502232}}
\multiput(48.751,30.965)(.0502232,-.0334821){16}{\line(1,0){.0502232}}
\multiput(50.358,29.894)(.0502232,-.0334821){16}{\line(1,0){.0502232}}
\multiput(51.965,28.823)(.0502232,-.0334821){16}{\line(1,0){.0502232}}
\multiput(53.573,27.751)(.0502232,-.0334821){16}{\line(1,0){.0502232}}
\multiput(55.18,26.68)(.0502232,-.0334821){16}{\line(1,0){.0502232}}
\multiput(56.787,25.608)(.0502232,-.0334821){16}{\line(1,0){.0502232}}
\multiput(58.394,24.537)(.0502232,-.0334821){16}{\line(1,0){.0502232}}
\multiput(60.001,23.465)(.0502232,-.0334821){16}{\line(1,0){.0502232}}
\multiput(61.608,22.394)(.0502232,-.0334821){16}{\line(1,0){.0502232}}
\multiput(63.215,21.323)(.0502232,-.0334821){16}{\line(1,0){.0502232}}
\multiput(64.823,20.251)(.0502232,-.0334821){16}{\line(1,0){.0502232}}
%\end
%\dashline{1}(47.5,27.5)(61,18.5)
\multiput(47.43,27.43)(.05,-.0333333){15}{\line(1,0){.05}}
\multiput(48.93,26.43)(.05,-.0333333){15}{\line(1,0){.05}}
\multiput(50.43,25.43)(.05,-.0333333){15}{\line(1,0){.05}}
\multiput(51.93,24.43)(.05,-.0333333){15}{\line(1,0){.05}}
\multiput(53.43,23.43)(.05,-.0333333){15}{\line(1,0){.05}}
\multiput(54.93,22.43)(.05,-.0333333){15}{\line(1,0){.05}}
\multiput(56.43,21.43)(.05,-.0333333){15}{\line(1,0){.05}}
\multiput(57.93,20.43)(.05,-.0333333){15}{\line(1,0){.05}}
\multiput(59.43,19.43)(.05,-.0333333){15}{\line(1,0){.05}}
%\end
%\dashline{1}(29.75,40.75)(6.75,28)
\multiput(29.68,40.68)(-.0608466,-.0337302){14}{\line(-1,0){.0608466}}
\multiput(27.976,39.735)(-.0608466,-.0337302){14}{\line(-1,0){.0608466}}
\multiput(26.272,38.791)(-.0608466,-.0337302){14}{\line(-1,0){.0608466}}
\multiput(24.569,37.846)(-.0608466,-.0337302){14}{\line(-1,0){.0608466}}
\multiput(22.865,36.902)(-.0608466,-.0337302){14}{\line(-1,0){.0608466}}
\multiput(21.161,35.957)(-.0608466,-.0337302){14}{\line(-1,0){.0608466}}
\multiput(19.457,35.013)(-.0608466,-.0337302){14}{\line(-1,0){.0608466}}
\multiput(17.754,34.069)(-.0608466,-.0337302){14}{\line(-1,0){.0608466}}
\multiput(16.05,33.124)(-.0608466,-.0337302){14}{\line(-1,0){.0608466}}
\multiput(14.346,32.18)(-.0608466,-.0337302){14}{\line(-1,0){.0608466}}
\multiput(12.643,31.235)(-.0608466,-.0337302){14}{\line(-1,0){.0608466}}
\multiput(10.939,30.291)(-.0608466,-.0337302){14}{\line(-1,0){.0608466}}
\multiput(9.235,29.346)(-.0608466,-.0337302){14}{\line(-1,0){.0608466}}
\multiput(7.532,28.402)(-.0608466,-.0337302){14}{\line(-1,0){.0608466}}
%\end
%\dashline{1}(128.5,40.5)(105.5,27.75)
\multiput(128.43,40.43)(-.0608466,-.0337302){14}{\line(-1,0){.0608466}}
\multiput(126.726,39.485)(-.0608466,-.0337302){14}{\line(-1,0){.0608466}}
\multiput(125.022,38.541)(-.0608466,-.0337302){14}{\line(-1,0){.0608466}}
\multiput(123.319,37.596)(-.0608466,-.0337302){14}{\line(-1,0){.0608466}}
\multiput(121.615,36.652)(-.0608466,-.0337302){14}{\line(-1,0){.0608466}}
\multiput(119.911,35.707)(-.0608466,-.0337302){14}{\line(-1,0){.0608466}}
\multiput(118.207,34.763)(-.0608466,-.0337302){14}{\line(-1,0){.0608466}}
\multiput(116.504,33.819)(-.0608466,-.0337302){14}{\line(-1,0){.0608466}}
\multiput(114.8,32.874)(-.0608466,-.0337302){14}{\line(-1,0){.0608466}}
\multiput(113.096,31.93)(-.0608466,-.0337302){14}{\line(-1,0){.0608466}}
\multiput(111.393,30.985)(-.0608466,-.0337302){14}{\line(-1,0){.0608466}}
\multiput(109.689,30.041)(-.0608466,-.0337302){14}{\line(-1,0){.0608466}}
\multiput(107.985,29.096)(-.0608466,-.0337302){14}{\line(-1,0){.0608466}}
\multiput(106.282,28.152)(-.0608466,-.0337302){14}{\line(-1,0){.0608466}}
%\end
%\dashline{1}(36,41)(4.75,23)
\multiput(35.93,40.93)(-.0563063,-.0324324){15}{\line(-1,0){.0563063}}
\multiput(34.241,39.957)(-.0563063,-.0324324){15}{\line(-1,0){.0563063}}
\multiput(32.551,38.984)(-.0563063,-.0324324){15}{\line(-1,0){.0563063}}
\multiput(30.862,38.011)(-.0563063,-.0324324){15}{\line(-1,0){.0563063}}
\multiput(29.173,37.038)(-.0563063,-.0324324){15}{\line(-1,0){.0563063}}
\multiput(27.484,36.065)(-.0563063,-.0324324){15}{\line(-1,0){.0563063}}
\multiput(25.795,35.092)(-.0563063,-.0324324){15}{\line(-1,0){.0563063}}
\multiput(24.105,34.119)(-.0563063,-.0324324){15}{\line(-1,0){.0563063}}
\multiput(22.416,33.146)(-.0563063,-.0324324){15}{\line(-1,0){.0563063}}
\multiput(20.727,32.173)(-.0563063,-.0324324){15}{\line(-1,0){.0563063}}
\multiput(19.038,31.2)(-.0563063,-.0324324){15}{\line(-1,0){.0563063}}
\multiput(17.349,30.227)(-.0563063,-.0324324){15}{\line(-1,0){.0563063}}
\multiput(15.659,29.254)(-.0563063,-.0324324){15}{\line(-1,0){.0563063}}
\multiput(13.97,28.281)(-.0563063,-.0324324){15}{\line(-1,0){.0563063}}
\multiput(12.281,27.308)(-.0563063,-.0324324){15}{\line(-1,0){.0563063}}
\multiput(10.592,26.335)(-.0563063,-.0324324){15}{\line(-1,0){.0563063}}
\multiput(8.903,25.362)(-.0563063,-.0324324){15}{\line(-1,0){.0563063}}
\multiput(7.213,24.389)(-.0563063,-.0324324){15}{\line(-1,0){.0563063}}
\multiput(5.524,23.416)(-.0563063,-.0324324){15}{\line(-1,0){.0563063}}
%\end
%\dashline{1}(134.75,40.75)(103.5,22.75)
\multiput(134.68,40.68)(-.0563063,-.0324324){15}{\line(-1,0){.0563063}}
\multiput(132.991,39.707)(-.0563063,-.0324324){15}{\line(-1,0){.0563063}}
\multiput(131.301,38.734)(-.0563063,-.0324324){15}{\line(-1,0){.0563063}}
\multiput(129.612,37.761)(-.0563063,-.0324324){15}{\line(-1,0){.0563063}}
\multiput(127.923,36.788)(-.0563063,-.0324324){15}{\line(-1,0){.0563063}}
\multiput(126.234,35.815)(-.0563063,-.0324324){15}{\line(-1,0){.0563063}}
\multiput(124.545,34.842)(-.0563063,-.0324324){15}{\line(-1,0){.0563063}}
\multiput(122.855,33.869)(-.0563063,-.0324324){15}{\line(-1,0){.0563063}}
\multiput(121.166,32.896)(-.0563063,-.0324324){15}{\line(-1,0){.0563063}}
\multiput(119.477,31.923)(-.0563063,-.0324324){15}{\line(-1,0){.0563063}}
\multiput(117.788,30.95)(-.0563063,-.0324324){15}{\line(-1,0){.0563063}}
\multiput(116.099,29.977)(-.0563063,-.0324324){15}{\line(-1,0){.0563063}}
\multiput(114.409,29.004)(-.0563063,-.0324324){15}{\line(-1,0){.0563063}}
\multiput(112.72,28.031)(-.0563063,-.0324324){15}{\line(-1,0){.0563063}}
\multiput(111.031,27.058)(-.0563063,-.0324324){15}{\line(-1,0){.0563063}}
\multiput(109.342,26.085)(-.0563063,-.0324324){15}{\line(-1,0){.0563063}}
\multiput(107.653,25.112)(-.0563063,-.0324324){15}{\line(-1,0){.0563063}}
\multiput(105.963,24.139)(-.0563063,-.0324324){15}{\line(-1,0){.0563063}}
\multiput(104.274,23.166)(-.0563063,-.0324324){15}{\line(-1,0){.0563063}}
%\end
%\dashline{1}(32.75,35.75)(6.25,20.25)
\multiput(32.68,35.68)(-.0552083,-.0322917){15}{\line(-1,0){.0552083}}
\multiput(31.023,34.711)(-.0552083,-.0322917){15}{\line(-1,0){.0552083}}
\multiput(29.367,33.742)(-.0552083,-.0322917){15}{\line(-1,0){.0552083}}
\multiput(27.711,32.773)(-.0552083,-.0322917){15}{\line(-1,0){.0552083}}
\multiput(26.055,31.805)(-.0552083,-.0322917){15}{\line(-1,0){.0552083}}
\multiput(24.398,30.836)(-.0552083,-.0322917){15}{\line(-1,0){.0552083}}
\multiput(22.742,29.867)(-.0552083,-.0322917){15}{\line(-1,0){.0552083}}
\multiput(21.086,28.898)(-.0552083,-.0322917){15}{\line(-1,0){.0552083}}
\multiput(19.43,27.93)(-.0552083,-.0322917){15}{\line(-1,0){.0552083}}
\multiput(17.773,26.961)(-.0552083,-.0322917){15}{\line(-1,0){.0552083}}
\multiput(16.117,25.992)(-.0552083,-.0322917){15}{\line(-1,0){.0552083}}
\multiput(14.461,25.023)(-.0552083,-.0322917){15}{\line(-1,0){.0552083}}
\multiput(12.805,24.055)(-.0552083,-.0322917){15}{\line(-1,0){.0552083}}
\multiput(11.148,23.086)(-.0552083,-.0322917){15}{\line(-1,0){.0552083}}
\multiput(9.492,22.117)(-.0552083,-.0322917){15}{\line(-1,0){.0552083}}
\multiput(7.836,21.148)(-.0552083,-.0322917){15}{\line(-1,0){.0552083}}
%\end
%\dashline{1}(131.5,35.5)(105,20)
\multiput(131.43,35.43)(-.0552083,-.0322917){15}{\line(-1,0){.0552083}}
\multiput(129.773,34.461)(-.0552083,-.0322917){15}{\line(-1,0){.0552083}}
\multiput(128.117,33.492)(-.0552083,-.0322917){15}{\line(-1,0){.0552083}}
\multiput(126.461,32.523)(-.0552083,-.0322917){15}{\line(-1,0){.0552083}}
\multiput(124.805,31.555)(-.0552083,-.0322917){15}{\line(-1,0){.0552083}}
\multiput(123.148,30.586)(-.0552083,-.0322917){15}{\line(-1,0){.0552083}}
\multiput(121.492,29.617)(-.0552083,-.0322917){15}{\line(-1,0){.0552083}}
\multiput(119.836,28.648)(-.0552083,-.0322917){15}{\line(-1,0){.0552083}}
\multiput(118.18,27.68)(-.0552083,-.0322917){15}{\line(-1,0){.0552083}}
\multiput(116.523,26.711)(-.0552083,-.0322917){15}{\line(-1,0){.0552083}}
\multiput(114.867,25.742)(-.0552083,-.0322917){15}{\line(-1,0){.0552083}}
\multiput(113.211,24.773)(-.0552083,-.0322917){15}{\line(-1,0){.0552083}}
\multiput(111.555,23.805)(-.0552083,-.0322917){15}{\line(-1,0){.0552083}}
\multiput(109.898,22.836)(-.0552083,-.0322917){15}{\line(-1,0){.0552083}}
\multiput(108.242,21.867)(-.0552083,-.0322917){15}{\line(-1,0){.0552083}}
\multiput(106.586,20.898)(-.0552083,-.0322917){15}{\line(-1,0){.0552083}}
%\end
%\dashline{1}(29.5,30.25)(12.25,19.25)
\multiput(29.43,30.18)(-.0522727,-.0333333){15}{\line(-1,0){.0522727}}
\multiput(27.862,29.18)(-.0522727,-.0333333){15}{\line(-1,0){.0522727}}
\multiput(26.293,28.18)(-.0522727,-.0333333){15}{\line(-1,0){.0522727}}
\multiput(24.725,27.18)(-.0522727,-.0333333){15}{\line(-1,0){.0522727}}
\multiput(23.157,26.18)(-.0522727,-.0333333){15}{\line(-1,0){.0522727}}
\multiput(21.589,25.18)(-.0522727,-.0333333){15}{\line(-1,0){.0522727}}
\multiput(20.021,24.18)(-.0522727,-.0333333){15}{\line(-1,0){.0522727}}
\multiput(18.452,23.18)(-.0522727,-.0333333){15}{\line(-1,0){.0522727}}
\multiput(16.884,22.18)(-.0522727,-.0333333){15}{\line(-1,0){.0522727}}
\multiput(15.316,21.18)(-.0522727,-.0333333){15}{\line(-1,0){.0522727}}
\multiput(13.748,20.18)(-.0522727,-.0333333){15}{\line(-1,0){.0522727}}
%\end
%\dashline{1}(128.25,30)(111,19)
\multiput(128.18,29.93)(-.0522727,-.0333333){15}{\line(-1,0){.0522727}}
\multiput(126.612,28.93)(-.0522727,-.0333333){15}{\line(-1,0){.0522727}}
\multiput(125.043,27.93)(-.0522727,-.0333333){15}{\line(-1,0){.0522727}}
\multiput(123.475,26.93)(-.0522727,-.0333333){15}{\line(-1,0){.0522727}}
\multiput(121.907,25.93)(-.0522727,-.0333333){15}{\line(-1,0){.0522727}}
\multiput(120.339,24.93)(-.0522727,-.0333333){15}{\line(-1,0){.0522727}}
\multiput(118.771,23.93)(-.0522727,-.0333333){15}{\line(-1,0){.0522727}}
\multiput(117.202,22.93)(-.0522727,-.0333333){15}{\line(-1,0){.0522727}}
\multiput(115.634,21.93)(-.0522727,-.0333333){15}{\line(-1,0){.0522727}}
\multiput(114.066,20.93)(-.0522727,-.0333333){15}{\line(-1,0){.0522727}}
\multiput(112.498,19.93)(-.0522727,-.0333333){15}{\line(-1,0){.0522727}}
%\end
\put(5.75,10.5){\makebox(0,0)[cc]{$u_1$}}
\put(104.5,10.25){\makebox(0,0)[cc]{$u_1$}}
\put(62.25,10.25){\makebox(0,0)[cc]{$u_2$}}
\put(35.75,13.5){\makebox(0,0)[cc]{$e$}}
\put(4.25,36.75){\makebox(0,0)[cc]{$H_1$}}
\put(103,36.5){\makebox(0,0)[cc]{$H_1$}}
\put(65.75,35.5){\makebox(0,0)[cc]{$H_2$}}
\put(38,48){\makebox(0,0)[cc]{$v$}}
\put(136.75,47.75){\makebox(0,0)[cc]{$v$}}
\qbezier(105,16)(141.75,13.125)(137.5,42.75)
\put(35,2){$G$}
\put(117.25,4){$G_1$}
\end{picture}

\caption{\relabel{f1} $G-e$ is separable.}
\end{figure}
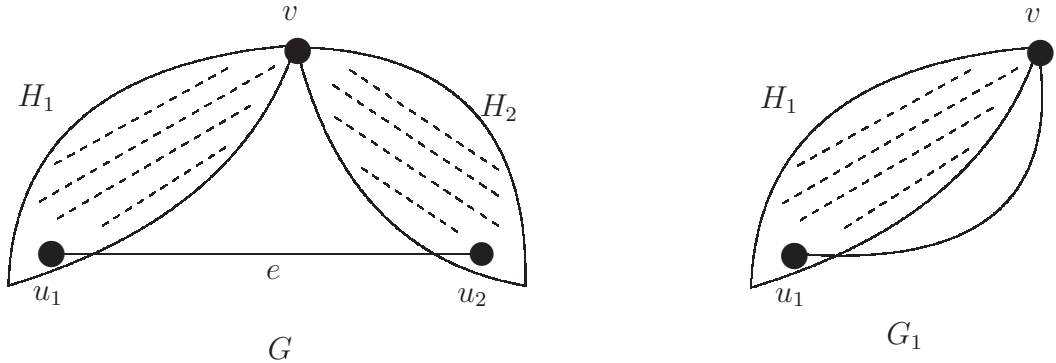

\begin{lem}[\cite{jac3}]\quad
\relabel{2-edge}
Let $G$ be a bridgeless connected graph, 
$S$ be a $2$-edge-cut of $G$, 
and $H_1$ and $H_2$ be the sides of $S$, as shown in 
Figure~\ref{f2}. 
Let $G_i$ be obtained from $G$ by contracting 
$E(H_{3-i})$, for $i\in \{1, 2\}$. Then
\beeq
F(G, \lambda) =
\frac{F(G_1, \lambda)F(G_2, \lambda)}{\lambda -1}.
\eneq
\end{lem}

\begin{figure}[h!]
\centering 
%\scalebox{1.0}{\input{f2.pic}}
%\input f2.pic
%TeXCAD (http://texcad.sf.net/) Picture. File: [f2.pic]. Options on following lines.
%\grade{\on}
%\emlines{\off}
%\epic{\off}
%\beziermacro{\on}
%\reduce{\on}
%\snapping{\off}
%\pvinsert{% Your \input, \def, etc. here}
%\quality{8.000}
%\graddiff{0.005}
%\snapasp{1}
%\zoom{4.0000}
\unitlength 1mm % = 2.845pt
\linethickness{0.4pt}
\ifx\plotpoint\undefined\newsavebox{\plotpoint}\fi % GNUPLOT compatibility
\begin{picture}(146.52,40)(0,0)
\put(11.25,37.25){\line(1,0){49}}
\put(12.25,13){\line(1,0){46.5}}
\put(-130.5,-84.5){\circle*{3.354}}
\put(-82,-84.5){\circle*{3.354}}
\put(12.25,13){\circle*{3.041}}
\put(107.75,13){\circle*{3.041}}
\put(12,36.75){\circle*{3.041}}
\put(107.5,36.75){\circle*{3.041}}
\put(145,25.25){\circle*{3.041}}
\put(-81.75,-85){\circle*{3.5}}
\put(59.75,13){\circle*{3.162}}
\put(59.5,36.75){\circle*{3.162}}
\put(11.5,25.25){\makebox(0,0)[cc]{$H_1$}}
\put(107,25.25){\makebox(0,0)[cc]{$H_1$}}
\put(59.75,25.25){\makebox(0,0)[cc]{$H_2$}}
%\dashline{1}(60,33.5)(69.25,20.25)
\multiput(59.93,33.43)(.0321181,-.0460069){16}{\line(0,-1){.0460069}}
\multiput(60.957,31.957)(.0321181,-.0460069){16}{\line(0,-1){.0460069}}
\multiput(61.985,30.485)(.0321181,-.0460069){16}{\line(0,-1){.0460069}}
\multiput(63.013,29.013)(.0321181,-.0460069){16}{\line(0,-1){.0460069}}
\multiput(64.041,27.541)(.0321181,-.0460069){16}{\line(0,-1){.0460069}}
\multiput(65.069,26.069)(.0321181,-.0460069){16}{\line(0,-1){.0460069}}
\multiput(66.096,24.596)(.0321181,-.0460069){16}{\line(0,-1){.0460069}}
\multiput(67.124,23.124)(.0321181,-.0460069){16}{\line(0,-1){.0460069}}
\multiput(68.152,21.652)(.0321181,-.0460069){16}{\line(0,-1){.0460069}}
%\end
%\dashline{1}(51,27.25)(59.5,17)
\multiput(50.93,27.18)(.0333333,-.0401961){17}{\line(0,-1){.0401961}}
\multiput(52.063,25.813)(.0333333,-.0401961){17}{\line(0,-1){.0401961}}
\multiput(53.196,24.446)(.0333333,-.0401961){17}{\line(0,-1){.0401961}}
\multiput(54.33,23.08)(.0333333,-.0401961){17}{\line(0,-1){.0401961}}
\multiput(55.463,21.713)(.0333333,-.0401961){17}{\line(0,-1){.0401961}}
\multiput(56.596,20.346)(.0333333,-.0401961){17}{\line(0,-1){.0401961}}
\multiput(57.73,18.98)(.0333333,-.0401961){17}{\line(0,-1){.0401961}}
\multiput(58.863,17.613)(.0333333,-.0401961){17}{\line(0,-1){.0401961}}
%\end
%\dashline{1}(6,30.75)(16.5,18)
\multiput(5.93,30.68)(.0324074,-.0393519){18}{\line(0,-1){.0393519}}
\multiput(7.096,29.263)(.0324074,-.0393519){18}{\line(0,-1){.0393519}}
\multiput(8.263,27.846)(.0324074,-.0393519){18}{\line(0,-1){.0393519}}
\multiput(9.43,26.43)(.0324074,-.0393519){18}{\line(0,-1){.0393519}}
\multiput(10.596,25.013)(.0324074,-.0393519){18}{\line(0,-1){.0393519}}
\multiput(11.763,23.596)(.0324074,-.0393519){18}{\line(0,-1){.0393519}}
\multiput(12.93,22.18)(.0324074,-.0393519){18}{\line(0,-1){.0393519}}
\multiput(14.096,20.763)(.0324074,-.0393519){18}{\line(0,-1){.0393519}}
\multiput(15.263,19.346)(.0324074,-.0393519){18}{\line(0,-1){.0393519}}
%\end
%\dashline{1}(101.5,30.75)(112,18)
\multiput(101.43,30.68)(.0324074,-.0393519){18}{\line(0,-1){.0393519}}
\multiput(102.596,29.263)(.0324074,-.0393519){18}{\line(0,-1){.0393519}}
\multiput(103.763,27.846)(.0324074,-.0393519){18}{\line(0,-1){.0393519}}
\multiput(104.93,26.43)(.0324074,-.0393519){18}{\line(0,-1){.0393519}}
\multiput(106.096,25.013)(.0324074,-.0393519){18}{\line(0,-1){.0393519}}
\multiput(107.263,23.596)(.0324074,-.0393519){18}{\line(0,-1){.0393519}}
\multiput(108.43,22.18)(.0324074,-.0393519){18}{\line(0,-1){.0393519}}
\multiput(109.596,20.763)(.0324074,-.0393519){18}{\line(0,-1){.0393519}}
\multiput(110.763,19.346)(.0324074,-.0393519){18}{\line(0,-1){.0393519}}
%\end
%\dashline{1}(4.75,24)(12.5,15)
\multiput(4.68,23.93)(.0331197,-.0384615){18}{\line(0,-1){.0384615}}
\multiput(5.872,22.545)(.0331197,-.0384615){18}{\line(0,-1){.0384615}}
\multiput(7.064,21.16)(.0331197,-.0384615){18}{\line(0,-1){.0384615}}
\multiput(8.257,19.776)(.0331197,-.0384615){18}{\line(0,-1){.0384615}}
\multiput(9.449,18.391)(.0331197,-.0384615){18}{\line(0,-1){.0384615}}
\multiput(10.641,17.007)(.0331197,-.0384615){18}{\line(0,-1){.0384615}}
\multiput(11.834,15.622)(.0331197,-.0384615){18}{\line(0,-1){.0384615}}
%\end
%\dashline{1}(100.25,24)(108,15)
\multiput(100.18,23.93)(.0331197,-.0384615){18}{\line(0,-1){.0384615}}
\multiput(101.372,22.545)(.0331197,-.0384615){18}{\line(0,-1){.0384615}}
\multiput(102.564,21.16)(.0331197,-.0384615){18}{\line(0,-1){.0384615}}
\multiput(103.757,19.776)(.0331197,-.0384615){18}{\line(0,-1){.0384615}}
\multiput(104.949,18.391)(.0331197,-.0384615){18}{\line(0,-1){.0384615}}
\multiput(106.141,17.007)(.0331197,-.0384615){18}{\line(0,-1){.0384615}}
\multiput(107.334,15.622)(.0331197,-.0384615){18}{\line(0,-1){.0384615}}
%\end
\put(12.5,25.125){\oval(20,29.75)[]}
\put(108,25.125){\oval(20,29.75)[]}
\put(61,25.125){\oval(20,29.75)[]}
%\dashline{1}(52,31.25)(65.75,13.75)
\multiput(51.93,31.18)(.0332126,-.0422705){18}{\line(0,-1){.0422705}}
\multiput(53.125,29.658)(.0332126,-.0422705){18}{\line(0,-1){.0422705}}
\multiput(54.321,28.136)(.0332126,-.0422705){18}{\line(0,-1){.0422705}}
\multiput(55.517,26.614)(.0332126,-.0422705){18}{\line(0,-1){.0422705}}
\multiput(56.712,25.093)(.0332126,-.0422705){18}{\line(0,-1){.0422705}}
\multiput(57.908,23.571)(.0332126,-.0422705){18}{\line(0,-1){.0422705}}
\multiput(59.104,22.049)(.0332126,-.0422705){18}{\line(0,-1){.0422705}}
\multiput(60.299,20.528)(.0332126,-.0422705){18}{\line(0,-1){.0422705}}
\multiput(61.495,19.006)(.0332126,-.0422705){18}{\line(0,-1){.0422705}}
\multiput(62.691,17.484)(.0332126,-.0422705){18}{\line(0,-1){.0422705}}
\multiput(63.886,15.962)(.0332126,-.0422705){18}{\line(0,-1){.0422705}}
\multiput(65.082,14.441)(.0332126,-.0422705){18}{\line(0,-1){.0422705}}
%\end
%\dashline{1}(53.5,35)(68.75,15.75)
\multiput(53.43,34.93)(.0325855,-.0411325){18}{\line(0,-1){.0411325}}
\multiput(54.603,33.449)(.0325855,-.0411325){18}{\line(0,-1){.0411325}}
\multiput(55.776,31.968)(.0325855,-.0411325){18}{\line(0,-1){.0411325}}
\multiput(56.949,30.487)(.0325855,-.0411325){18}{\line(0,-1){.0411325}}
\multiput(58.122,29.007)(.0325855,-.0411325){18}{\line(0,-1){.0411325}}
\multiput(59.295,27.526)(.0325855,-.0411325){18}{\line(0,-1){.0411325}}
\multiput(60.468,26.045)(.0325855,-.0411325){18}{\line(0,-1){.0411325}}
\multiput(61.641,24.564)(.0325855,-.0411325){18}{\line(0,-1){.0411325}}
\multiput(62.814,23.084)(.0325855,-.0411325){18}{\line(0,-1){.0411325}}
\multiput(63.987,21.603)(.0325855,-.0411325){18}{\line(0,-1){.0411325}}
\multiput(65.16,20.122)(.0325855,-.0411325){18}{\line(0,-1){.0411325}}
\multiput(66.334,18.641)(.0325855,-.0411325){18}{\line(0,-1){.0411325}}
\multiput(67.507,17.16)(.0325855,-.0411325){18}{\line(0,-1){.0411325}}
%\end
%\dashline{1}(62.5,35.25)(69.75,25)
\multiput(62.43,35.18)(.0323661,-.0457589){16}{\line(0,-1){.0457589}}
\multiput(63.465,33.715)(.0323661,-.0457589){16}{\line(0,-1){.0457589}}
\multiput(64.501,32.251)(.0323661,-.0457589){16}{\line(0,-1){.0457589}}
\multiput(65.537,30.787)(.0323661,-.0457589){16}{\line(0,-1){.0457589}}
\multiput(66.573,29.323)(.0323661,-.0457589){16}{\line(0,-1){.0457589}}
\multiput(67.608,27.858)(.0323661,-.0457589){16}{\line(0,-1){.0457589}}
\multiput(68.644,26.394)(.0323661,-.0457589){16}{\line(0,-1){.0457589}}
%\end
%\dashline{1}(52.25,20.75)(56.75,15.75)
\multiput(52.18,20.68)(.0330882,-.0367647){17}{\line(0,-1){.0367647}}
\multiput(53.305,19.43)(.0330882,-.0367647){17}{\line(0,-1){.0367647}}
\multiput(54.43,18.18)(.0330882,-.0367647){17}{\line(0,-1){.0367647}}
\multiput(55.555,16.93)(.0330882,-.0367647){17}{\line(0,-1){.0367647}}
%\end
%\dashline{1}(65,36.5)(70,29.75)
\multiput(64.93,36.43)(.0333333,-.045){15}{\line(0,-1){.045}}
\multiput(65.93,35.08)(.0333333,-.045){15}{\line(0,-1){.045}}
\multiput(66.93,33.73)(.0333333,-.045){15}{\line(0,-1){.045}}
\multiput(67.93,32.38)(.0333333,-.045){15}{\line(0,-1){.045}}
\multiput(68.93,31.03)(.0333333,-.045){15}{\line(0,-1){.045}}
%\end
%\dashline{1}(4.75,16.75)(8.75,12.75)
\multiput(4.68,16.68)(.0333333,-.0333333){15}{\line(0,-1){.0333333}}
\multiput(5.68,15.68)(.0333333,-.0333333){15}{\line(0,-1){.0333333}}
\multiput(6.68,14.68)(.0333333,-.0333333){15}{\line(0,-1){.0333333}}
\multiput(7.68,13.68)(.0333333,-.0333333){15}{\line(0,-1){.0333333}}
%\end
%\dashline{1}(100.25,16.75)(104.25,12.75)
\multiput(100.18,16.68)(.0333333,-.0333333){15}{\line(0,-1){.0333333}}
\multiput(101.18,15.68)(.0333333,-.0333333){15}{\line(0,-1){.0333333}}
\multiput(102.18,14.68)(.0333333,-.0333333){15}{\line(0,-1){.0333333}}
\multiput(103.18,13.68)(.0333333,-.0333333){15}{\line(0,-1){.0333333}}
%\end
%\dashline{1}(6.25,35.25)(21.5,17.5)
\multiput(6.18,35.18)(.0321053,-.0373684){19}{\line(0,-1){.0373684}}
\multiput(7.4,33.76)(.0321053,-.0373684){19}{\line(0,-1){.0373684}}
\multiput(8.62,32.34)(.0321053,-.0373684){19}{\line(0,-1){.0373684}}
\multiput(9.84,30.92)(.0321053,-.0373684){19}{\line(0,-1){.0373684}}
\multiput(11.06,29.5)(.0321053,-.0373684){19}{\line(0,-1){.0373684}}
\multiput(12.28,28.08)(.0321053,-.0373684){19}{\line(0,-1){.0373684}}
\multiput(13.5,26.66)(.0321053,-.0373684){19}{\line(0,-1){.0373684}}
\multiput(14.72,25.24)(.0321053,-.0373684){19}{\line(0,-1){.0373684}}
\multiput(15.94,23.82)(.0321053,-.0373684){19}{\line(0,-1){.0373684}}
\multiput(17.16,22.4)(.0321053,-.0373684){19}{\line(0,-1){.0373684}}
\multiput(18.38,20.98)(.0321053,-.0373684){19}{\line(0,-1){.0373684}}
\multiput(19.6,19.56)(.0321053,-.0373684){19}{\line(0,-1){.0373684}}
\multiput(20.82,18.14)(.0321053,-.0373684){19}{\line(0,-1){.0373684}}
%\end
%\dashline{1}(101.75,35.25)(117,17.5)
\multiput(101.68,35.18)(.0321053,-.0373684){19}{\line(0,-1){.0373684}}
\multiput(102.9,33.76)(.0321053,-.0373684){19}{\line(0,-1){.0373684}}
\multiput(104.12,32.34)(.0321053,-.0373684){19}{\line(0,-1){.0373684}}
\multiput(105.34,30.92)(.0321053,-.0373684){19}{\line(0,-1){.0373684}}
\multiput(106.56,29.5)(.0321053,-.0373684){19}{\line(0,-1){.0373684}}
\multiput(107.78,28.08)(.0321053,-.0373684){19}{\line(0,-1){.0373684}}
\multiput(109,26.66)(.0321053,-.0373684){19}{\line(0,-1){.0373684}}
\multiput(110.22,25.24)(.0321053,-.0373684){19}{\line(0,-1){.0373684}}
\multiput(111.44,23.82)(.0321053,-.0373684){19}{\line(0,-1){.0373684}}
\multiput(112.66,22.4)(.0321053,-.0373684){19}{\line(0,-1){.0373684}}
\multiput(113.88,20.98)(.0321053,-.0373684){19}{\line(0,-1){.0373684}}
\multiput(115.1,19.56)(.0321053,-.0373684){19}{\line(0,-1){.0373684}}
\multiput(116.32,18.14)(.0321053,-.0373684){19}{\line(0,-1){.0373684}}
%\end
%\dashline{1}(4.75,27.75)(15.75,15.25)
\multiput(4.68,27.68)(.0321637,-.0365497){19}{\line(0,-1){.0365497}}
\multiput(5.902,26.291)(.0321637,-.0365497){19}{\line(0,-1){.0365497}}
\multiput(7.124,24.902)(.0321637,-.0365497){19}{\line(0,-1){.0365497}}
\multiput(8.346,23.513)(.0321637,-.0365497){19}{\line(0,-1){.0365497}}
\multiput(9.569,22.124)(.0321637,-.0365497){19}{\line(0,-1){.0365497}}
\multiput(10.791,20.735)(.0321637,-.0365497){19}{\line(0,-1){.0365497}}
\multiput(12.013,19.346)(.0321637,-.0365497){19}{\line(0,-1){.0365497}}
\multiput(13.235,17.957)(.0321637,-.0365497){19}{\line(0,-1){.0365497}}
\multiput(14.457,16.569)(.0321637,-.0365497){19}{\line(0,-1){.0365497}}
%\end
%\dashline{1}(100.25,27.75)(111.25,15.25)
\multiput(100.18,27.68)(.0321637,-.0365497){19}{\line(0,-1){.0365497}}
\multiput(101.402,26.291)(.0321637,-.0365497){19}{\line(0,-1){.0365497}}
\multiput(102.624,24.902)(.0321637,-.0365497){19}{\line(0,-1){.0365497}}
\multiput(103.846,23.513)(.0321637,-.0365497){19}{\line(0,-1){.0365497}}
\multiput(105.069,22.124)(.0321637,-.0365497){19}{\line(0,-1){.0365497}}
\multiput(106.291,20.735)(.0321637,-.0365497){19}{\line(0,-1){.0365497}}
\multiput(107.513,19.346)(.0321637,-.0365497){19}{\line(0,-1){.0365497}}
\multiput(108.735,17.957)(.0321637,-.0365497){19}{\line(0,-1){.0365497}}
\multiput(109.957,16.569)(.0321637,-.0365497){19}{\line(0,-1){.0365497}}
%\end
%\dashline{1}(4.25,20)(10,14.25)
\multiput(4.18,19.93)(.0336257,-.0336257){19}{\line(0,-1){.0336257}}
\multiput(5.457,18.652)(.0336257,-.0336257){19}{\line(0,-1){.0336257}}
\multiput(6.735,17.374)(.0336257,-.0336257){19}{\line(0,-1){.0336257}}
\multiput(8.013,16.096)(.0336257,-.0336257){19}{\line(0,-1){.0336257}}
\multiput(9.291,14.819)(.0336257,-.0336257){19}{\line(0,-1){.0336257}}
%\end
%\dashline{1}(99.75,20)(105.5,14.25)
\multiput(99.68,19.93)(.0336257,-.0336257){19}{\line(0,-1){.0336257}}
\multiput(100.957,18.652)(.0336257,-.0336257){19}{\line(1,0){.0336257}}
\multiput(102.235,17.374)(.0336257,-.0336257){19}{\line(0,-1){.0336257}}
\multiput(103.513,16.096)(.0336257,-.0336257){19}{\line(1,0){.0336257}}
\multiput(104.791,14.819)(.0336257,-.0336257){19}{\line(0,-1){.0336257}}
%\end
%\dashline{1}(12.5,33.25)(22.25,21.25)
\multiput(12.43,33.18)(.033737,-.0415225){17}{\line(0,-1){.0415225}}
\multiput(13.577,31.768)(.033737,-.0415225){17}{\line(0,-1){.0415225}}
\multiput(14.724,30.356)(.033737,-.0415225){17}{\line(0,-1){.0415225}}
\multiput(15.871,28.944)(.033737,-.0415225){17}{\line(0,-1){.0415225}}
\multiput(17.018,27.533)(.033737,-.0415225){17}{\line(0,-1){.0415225}}
\multiput(18.165,26.121)(.033737,-.0415225){17}{\line(0,-1){.0415225}}
\multiput(19.312,24.709)(.033737,-.0415225){17}{\line(0,-1){.0415225}}
\multiput(20.459,23.297)(.033737,-.0415225){17}{\line(0,-1){.0415225}}
\multiput(21.606,21.886)(.033737,-.0415225){17}{\line(0,-1){.0415225}}
%\end
%\dashline{1}(108,33.25)(117.75,21.25)
\multiput(107.93,33.18)(.033737,-.0415225){17}{\line(0,-1){.0415225}}
\multiput(109.077,31.768)(.033737,-.0415225){17}{\line(0,-1){.0415225}}
\multiput(110.224,30.356)(.033737,-.0415225){17}{\line(0,-1){.0415225}}
\multiput(111.371,28.944)(.033737,-.0415225){17}{\line(0,-1){.0415225}}
\multiput(112.518,27.533)(.033737,-.0415225){17}{\line(0,-1){.0415225}}
\multiput(113.665,26.121)(.033737,-.0415225){17}{\line(0,-1){.0415225}}
\multiput(114.812,24.709)(.033737,-.0415225){17}{\line(0,-1){.0415225}}
\multiput(115.959,23.297)(.033737,-.0415225){17}{\line(0,-1){.0415225}}
\multiput(117.106,21.886)(.033737,-.0415225){17}{\line(0,-1){.0415225}}
%\end
%\dashline{1}(14.75,34.75)(22.5,25.5)
\multiput(14.68,34.68)(.0331197,-.0395299){18}{\line(0,-1){.0395299}}
\multiput(15.872,33.257)(.0331197,-.0395299){18}{\line(0,-1){.0395299}}
\multiput(17.064,31.834)(.0331197,-.0395299){18}{\line(0,-1){.0395299}}
\multiput(18.257,30.41)(.0331197,-.0395299){18}{\line(0,-1){.0395299}}
\multiput(19.449,28.987)(.0331197,-.0395299){18}{\line(0,-1){.0395299}}
\multiput(20.641,27.564)(.0331197,-.0395299){18}{\line(0,-1){.0395299}}
\multiput(21.834,26.141)(.0331197,-.0395299){18}{\line(0,-1){.0395299}}
%\end
%\dashline{1}(110.25,34.75)(118,25.5)
\multiput(110.18,34.68)(.0331197,-.0395299){18}{\line(0,-1){.0395299}}
\multiput(111.372,33.257)(.0331197,-.0395299){18}{\line(0,-1){.0395299}}
\multiput(112.564,31.834)(.0331197,-.0395299){18}{\line(0,-1){.0395299}}
\multiput(113.757,30.41)(.0331197,-.0395299){18}{\line(0,-1){.0395299}}
\multiput(114.949,28.987)(.0331197,-.0395299){18}{\line(0,-1){.0395299}}
\multiput(116.141,27.564)(.0331197,-.0395299){18}{\line(0,-1){.0395299}}
\multiput(117.334,26.141)(.0331197,-.0395299){18}{\line(0,-1){.0395299}}
%\end
%\dashline{1}(16.75,35.75)(21.75,29.75)
\multiput(16.68,35.68)(.0333333,-.04){15}{\line(0,-1){.04}}
\multiput(17.68,34.48)(.0333333,-.04){15}{\line(0,-1){.04}}
\multiput(18.68,33.28)(.0333333,-.04){15}{\line(0,-1){.04}}
\multiput(19.68,32.08)(.0333333,-.04){15}{\line(0,-1){.04}}
\multiput(20.68,30.88)(.0333333,-.04){15}{\line(0,-1){.04}}
%\end
%\dashline{1}(112.25,35.75)(117.25,29.75)
\multiput(112.18,35.68)(.0333333,-.04){15}{\line(0,-1){.04}}
\multiput(113.18,34.48)(.0333333,-.04){15}{\line(0,-1){.04}}
\multiput(114.18,33.28)(.0333333,-.04){15}{\line(0,-1){.04}}
\multiput(115.18,32.08)(.0333333,-.04){15}{\line(0,-1){.04}}
\multiput(116.18,30.88)(.0333333,-.04){15}{\line(0,-1){.04}}
%\end
%\emline(107.25,37)(144.25,25.5)
\multiput(107.25,37)(.1085043988,-.0337243402){341}{\line(1,0){.1085043988}}
%\end
%\emline(144.25,25.5)(107.5,13)
\multiput(144.25,25.5)(-.0990566038,-.0336927224){371}{\line(-1,0){.0990566038}}
%\end
\put(34.5,3.5){$G$}
\put(100,5.75){$G_1$ when $H_2$ is connected}
\end{picture}

\caption{\relabel{f2} $G$ has a $2$-edge-cut.}
\end{figure}
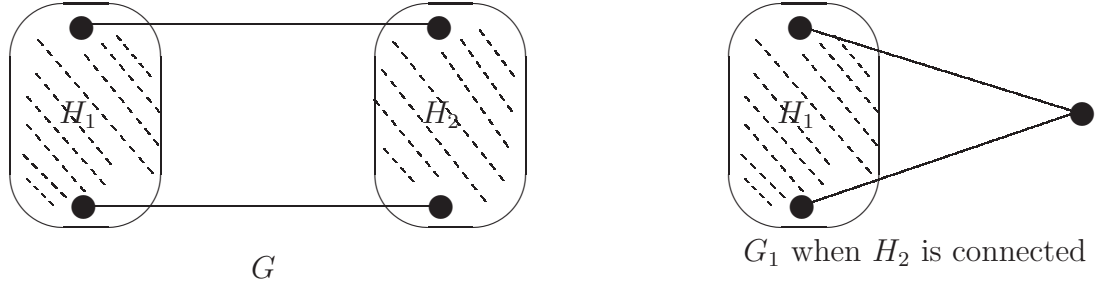

\begin{lem}[\cite{jac3}]\quad
\relabel{3-edge}
Let $G$ be a bridgeless connected  graph, 
$S$ be a $3$-edge-cut of $G$, 
and $H_1$ and $H_2$ be the sides of $S$. 
Let $G_i$ be obtained from $G$ by contracting 
$E(H_{3-i})$, for $i\in \{1, 2\}$.
Then
\beeq
F(G, \lambda) = 
\frac{F(G_1, \lambda)F(G_2, \lambda)}
{(\lambda -1)(\lambda -2)}. 
\eneq
\end{lem}

\noindent{\bf Remark}: 
For a non-separable graph $G$, if 
$G$ contains a 2-edge-cut, then $G-e$ is separable for 
each $e$ in this cut and thus 
Lemma~\ref{2-edge} is a special case of Lemma~\ref{v-edge}.
Also note that the graph in Lemma~\ref{3-edge} 
has a structure similar to the one in Figure~\ref{f2}.

We end this section with the following result 
which will be applied many times in this paper.

\begin{lem}\quad
\relabel{2-blocks}
Let $G$ be a non-separable 
graph with subgraphs $G_1$ and $G_2$ such that 
$V(G_1)\cap V(G_2)=\{u,v\}$,
$V(G_1)\cup V(G_2)=V(G)$,
$E(G_1)\cap E(G_2)=\emptyset$
and $E(G_1)\cup E(G_2)=E(G)$, as shown in Figure~\ref{f3}(a).
Then 
\beeq
F(G, \lambda) = 
\frac{F(G_1+uv, \lambda)F(G_2+uv, \lambda)}{\lambda -1}
+F(G_1, \lambda)F(G_2, \lambda),
\eneq
where $u$ and $v$ be two vertives of $G$.
\end{lem}

\begin{figure}[h!]
\centering 
\scalebox{0.75}
%{\input{f3.pic}}
%\input f5.pic
%TeXCAD (http://texcad.sf.net/) Picture. File: [f7.pic]. Options on following lines.
%\grade{\on}
%\emlines{\off}
%\epic{\off}
%\beziermacro{\on}
%\reduce{\on}
%\snapping{\off}
%\pvinsert{% Your \input, \def, etc. here}
%\quality{8.000}
%\graddiff{0.005}
%\snapasp{1}
%\zoom{4.0000}
{
\unitlength 1mm % = 2.845pt
\linethickness{0.4pt}
\ifx\plotpoint\undefined\newsavebox{\plotpoint}\fi % GNUPLOT compatibility
\begin{picture}(205,66.5)(0,0)
\put(23.75,54.25){\circle*{4.61}}
\put(24,54.5){\circle*{4.61}}
\put(83.5,54.5){\circle*{4.61}}
\put(151.25,54.5){\circle*{4.61}}
\put(23.75,14.5){\circle*{4.61}}
\put(93.5,14.5){\circle*{4.61}}
\put(161.25,15.25){\circle*{4.61}}
\put(93.5,14.5){\circle*{4.61}}
\put(161.25,15.25){\circle*{4.61}}
\put(93.5,14.5){\circle*{4.61}}
\put(161.25,15.25){\circle*{4.61}}
\put(73.5,14.5){\circle*{4.61}}
\put(141.25,15.25){\circle*{4.61}}
\put(24,14.75){\circle*{4.61}}
\qbezier(24,54.75)(9.5,31.375)(24,14.5)
\qbezier(23.75,55)(38.25,31.625)(23.75,14.75)
\qbezier(23.75,54)(-20.875,32.125)(24,14.75)
\qbezier(24,54.25)(68.625,32.375)(23.75,15)
%\dashline{1}(17.75,48.75)(3.25,36.5)
\multiput(17.68,48.68)(-.0381579,-.0322368){19}{\line(-1,0){.0381579}}
\multiput(16.23,47.455)(-.0381579,-.0322368){19}{\line(-1,0){.0381579}}
\multiput(14.78,46.23)(-.0381579,-.0322368){19}{\line(-1,0){.0381579}}
\multiput(13.33,45.005)(-.0381579,-.0322368){19}{\line(-1,0){.0381579}}
\multiput(11.88,43.78)(-.0381579,-.0322368){19}{\line(-1,0){.0381579}}
\multiput(10.43,42.555)(-.0381579,-.0322368){19}{\line(-1,0){.0381579}}
\multiput(8.98,41.33)(-.0381579,-.0322368){19}{\line(-1,0){.0381579}}
\multiput(7.53,40.105)(-.0381579,-.0322368){19}{\line(-1,0){.0381579}}
\multiput(6.08,38.88)(-.0381579,-.0322368){19}{\line(-1,0){.0381579}}
\multiput(4.63,37.655)(-.0381579,-.0322368){19}{\line(-1,0){.0381579}}
%\end
%\dashline{1}(30,49)(44.5,36.75)
\multiput(29.93,48.93)(.0381579,-.0322368){19}{\line(1,0){.0381579}}
\multiput(31.38,47.705)(.0381579,-.0322368){19}{\line(1,0){.0381579}}
\multiput(32.83,46.48)(.0381579,-.0322368){19}{\line(1,0){.0381579}}
\multiput(34.28,45.255)(.0381579,-.0322368){19}{\line(1,0){.0381579}}
\multiput(35.73,44.03)(.0381579,-.0322368){19}{\line(1,0){.0381579}}
\multiput(37.18,42.805)(.0381579,-.0322368){19}{\line(1,0){.0381579}}
\multiput(38.63,41.58)(.0381579,-.0322368){19}{\line(1,0){.0381579}}
\multiput(40.08,40.355)(.0381579,-.0322368){19}{\line(1,0){.0381579}}
\multiput(41.53,39.13)(.0381579,-.0322368){19}{\line(1,0){.0381579}}
\multiput(42.98,37.905)(.0381579,-.0322368){19}{\line(1,0){.0381579}}
%\end
%\dashline{1}(97,49.25)(111.5,37)
\multiput(96.93,49.18)(.0381579,-.0322368){19}{\line(1,0){.0381579}}
\multiput(98.38,47.955)(.0381579,-.0322368){19}{\line(1,0){.0381579}}
\multiput(99.83,46.73)(.0381579,-.0322368){19}{\line(1,0){.0381579}}
\multiput(101.28,45.505)(.0381579,-.0322368){19}{\line(1,0){.0381579}}
\multiput(102.73,44.28)(.0381579,-.0322368){19}{\line(1,0){.0381579}}
\multiput(104.18,43.055)(.0381579,-.0322368){19}{\line(1,0){.0381579}}
\multiput(105.63,41.83)(.0381579,-.0322368){19}{\line(1,0){.0381579}}
\multiput(107.08,40.605)(.0381579,-.0322368){19}{\line(1,0){.0381579}}
\multiput(108.53,39.38)(.0381579,-.0322368){19}{\line(1,0){.0381579}}
\multiput(109.98,38.155)(.0381579,-.0322368){19}{\line(1,0){.0381579}}
%\end
%\dashline{1}(164.75,50)(179.25,37.75)
\multiput(164.68,49.93)(.0381579,-.0322368){19}{\line(1,0){.0381579}}
\multiput(166.13,48.705)(.0381579,-.0322368){19}{\line(1,0){.0381579}}
\multiput(167.58,47.48)(.0381579,-.0322368){19}{\line(1,0){.0381579}}
\multiput(169.03,46.255)(.0381579,-.0322368){19}{\line(1,0){.0381579}}
\multiput(170.48,45.03)(.0381579,-.0322368){19}{\line(1,0){.0381579}}
\multiput(171.93,43.805)(.0381579,-.0322368){19}{\line(1,0){.0381579}}
\multiput(173.38,42.58)(.0381579,-.0322368){19}{\line(1,0){.0381579}}
\multiput(174.83,41.355)(.0381579,-.0322368){19}{\line(1,0){.0381579}}
\multiput(176.28,40.13)(.0381579,-.0322368){19}{\line(1,0){.0381579}}
\multiput(177.73,38.905)(.0381579,-.0322368){19}{\line(1,0){.0381579}}
%\end
%\dashline{1}(97,49.25)(111.5,37)
\multiput(96.93,49.18)(.0381579,-.0322368){19}{\line(1,0){.0381579}}
\multiput(98.38,47.955)(.0381579,-.0322368){19}{\line(1,0){.0381579}}
\multiput(99.83,46.73)(.0381579,-.0322368){19}{\line(1,0){.0381579}}
\multiput(101.28,45.505)(.0381579,-.0322368){19}{\line(1,0){.0381579}}
\multiput(102.73,44.28)(.0381579,-.0322368){19}{\line(1,0){.0381579}}
\multiput(104.18,43.055)(.0381579,-.0322368){19}{\line(1,0){.0381579}}
\multiput(105.63,41.83)(.0381579,-.0322368){19}{\line(1,0){.0381579}}
\multiput(107.08,40.605)(.0381579,-.0322368){19}{\line(1,0){.0381579}}
\multiput(108.53,39.38)(.0381579,-.0322368){19}{\line(1,0){.0381579}}
\multiput(109.98,38.155)(.0381579,-.0322368){19}{\line(1,0){.0381579}}
%\end
%\dashline{1}(164.75,50)(179.25,37.75)
\multiput(164.68,49.93)(.0381579,-.0322368){19}{\line(1,0){.0381579}}
\multiput(166.13,48.705)(.0381579,-.0322368){19}{\line(1,0){.0381579}}
\multiput(167.58,47.48)(.0381579,-.0322368){19}{\line(1,0){.0381579}}
\multiput(169.03,46.255)(.0381579,-.0322368){19}{\line(1,0){.0381579}}
\multiput(170.48,45.03)(.0381579,-.0322368){19}{\line(1,0){.0381579}}
\multiput(171.93,43.805)(.0381579,-.0322368){19}{\line(1,0){.0381579}}
\multiput(173.38,42.58)(.0381579,-.0322368){19}{\line(1,0){.0381579}}
\multiput(174.83,41.355)(.0381579,-.0322368){19}{\line(1,0){.0381579}}
\multiput(176.28,40.13)(.0381579,-.0322368){19}{\line(1,0){.0381579}}
\multiput(177.73,38.905)(.0381579,-.0322368){19}{\line(1,0){.0381579}}
%\end
%\dashline{1}(97,49.25)(111.5,37)
\multiput(96.93,49.18)(.0381579,-.0322368){19}{\line(1,0){.0381579}}
\multiput(98.38,47.955)(.0381579,-.0322368){19}{\line(1,0){.0381579}}
\multiput(99.83,46.73)(.0381579,-.0322368){19}{\line(1,0){.0381579}}
\multiput(101.28,45.505)(.0381579,-.0322368){19}{\line(1,0){.0381579}}
\multiput(102.73,44.28)(.0381579,-.0322368){19}{\line(1,0){.0381579}}
\multiput(104.18,43.055)(.0381579,-.0322368){19}{\line(1,0){.0381579}}
\multiput(105.63,41.83)(.0381579,-.0322368){19}{\line(1,0){.0381579}}
\multiput(107.08,40.605)(.0381579,-.0322368){19}{\line(1,0){.0381579}}
\multiput(108.53,39.38)(.0381579,-.0322368){19}{\line(1,0){.0381579}}
\multiput(109.98,38.155)(.0381579,-.0322368){19}{\line(1,0){.0381579}}
%\end
%\dashline{1}(164.75,50)(179.25,37.75)
\multiput(164.68,49.93)(.0381579,-.0322368){19}{\line(1,0){.0381579}}
\multiput(166.13,48.705)(.0381579,-.0322368){19}{\line(1,0){.0381579}}
\multiput(167.58,47.48)(.0381579,-.0322368){19}{\line(1,0){.0381579}}
\multiput(169.03,46.255)(.0381579,-.0322368){19}{\line(1,0){.0381579}}
\multiput(170.48,45.03)(.0381579,-.0322368){19}{\line(1,0){.0381579}}
\multiput(171.93,43.805)(.0381579,-.0322368){19}{\line(1,0){.0381579}}
\multiput(173.38,42.58)(.0381579,-.0322368){19}{\line(1,0){.0381579}}
\multiput(174.83,41.355)(.0381579,-.0322368){19}{\line(1,0){.0381579}}
\multiput(176.28,40.13)(.0381579,-.0322368){19}{\line(1,0){.0381579}}
\multiput(177.73,38.905)(.0381579,-.0322368){19}{\line(1,0){.0381579}}
%\end
%\dashline{1}(70,49.25)(55.5,37)
\multiput(69.93,49.18)(-.0381579,-.0322368){19}{\line(-1,0){.0381579}}
\multiput(68.48,47.955)(-.0381579,-.0322368){19}{\line(-1,0){.0381579}}
\multiput(67.03,46.73)(-.0381579,-.0322368){19}{\line(-1,0){.0381579}}
\multiput(65.58,45.505)(-.0381579,-.0322368){19}{\line(-1,0){.0381579}}
\multiput(64.13,44.28)(-.0381579,-.0322368){19}{\line(-1,0){.0381579}}
\multiput(62.68,43.055)(-.0381579,-.0322368){19}{\line(-1,0){.0381579}}
\multiput(61.23,41.83)(-.0381579,-.0322368){19}{\line(-1,0){.0381579}}
\multiput(59.78,40.605)(-.0381579,-.0322368){19}{\line(-1,0){.0381579}}
\multiput(58.33,39.38)(-.0381579,-.0322368){19}{\line(-1,0){.0381579}}
\multiput(56.88,38.155)(-.0381579,-.0322368){19}{\line(-1,0){.0381579}}
%\end
%\dashline{1}(137.75,50)(123.25,37.75)
\multiput(137.68,49.93)(-.0381579,-.0322368){19}{\line(-1,0){.0381579}}
\multiput(136.23,48.705)(-.0381579,-.0322368){19}{\line(-1,0){.0381579}}
\multiput(134.78,47.48)(-.0381579,-.0322368){19}{\line(-1,0){.0381579}}
\multiput(133.33,46.255)(-.0381579,-.0322368){19}{\line(-1,0){.0381579}}
\multiput(131.88,45.03)(-.0381579,-.0322368){19}{\line(-1,0){.0381579}}
\multiput(130.43,43.805)(-.0381579,-.0322368){19}{\line(-1,0){.0381579}}
\multiput(128.98,42.58)(-.0381579,-.0322368){19}{\line(-1,0){.0381579}}
\multiput(127.53,41.355)(-.0381579,-.0322368){19}{\line(-1,0){.0381579}}
\multiput(126.08,40.13)(-.0381579,-.0322368){19}{\line(-1,0){.0381579}}
\multiput(124.63,38.905)(-.0381579,-.0322368){19}{\line(-1,0){.0381579}}
%\end
%\dashline{1}(3.25,36.5)(3.25,36.5)
\put(3.18,36.43){\line(0,1){0}}
%\end
%\dashline{1}(44.5,36.75)(44.5,36.75)
\put(44.43,36.68){\line(0,1){0}}
%\end
%\dashline{1}(113,36.75)(113,36.75)
\put(112.93,36.68){\line(0,1){0}}
%\end
%\dashline{1}(180.75,37.5)(180.75,37.5)
\put(180.68,37.43){\line(0,1){0}}
%\end
%\dashline{1}(113,36.75)(113,36.75)
\put(112.93,36.68){\line(0,1){0}}
%\end
%\dashline{1}(180.75,37.5)(180.75,37.5)
\put(180.68,37.43){\line(0,1){0}}
%\end
%\dashline{1}(113,36.75)(113,36.75)
\put(112.93,36.68){\line(0,1){0}}
%\end
%\dashline{1}(180.75,37.5)(180.75,37.5)
\put(180.68,37.43){\line(0,1){0}}
%\end
%\dashline{1}(54,36.75)(54,36.75)
\put(53.93,36.68){\line(0,1){0}}
%\end
%\dashline{1}(121.75,37.5)(121.75,37.5)
\put(121.68,37.43){\line(0,1){0}}
%\end
%\dashline{1}(17.5,44.5)(3,32)
\multiput(17.43,44.43)(-.0381579,-.0328947){19}{\line(-1,0){.0381579}}
\multiput(15.98,43.18)(-.0381579,-.0328947){19}{\line(-1,0){.0381579}}
\multiput(14.53,41.93)(-.0381579,-.0328947){19}{\line(-1,0){.0381579}}
\multiput(13.08,40.68)(-.0381579,-.0328947){19}{\line(-1,0){.0381579}}
\multiput(11.63,39.43)(-.0381579,-.0328947){19}{\line(-1,0){.0381579}}
\multiput(10.18,38.18)(-.0381579,-.0328947){19}{\line(-1,0){.0381579}}
\multiput(8.73,36.93)(-.0381579,-.0328947){19}{\line(-1,0){.0381579}}
\multiput(7.28,35.68)(-.0381579,-.0328947){19}{\line(-1,0){.0381579}}
\multiput(5.83,34.43)(-.0381579,-.0328947){19}{\line(-1,0){.0381579}}
\multiput(4.38,33.18)(-.0381579,-.0328947){19}{\line(-1,0){.0381579}}
%\end
%\dashline{1}(30.25,44.75)(44.75,32.25)
\multiput(30.18,44.68)(.0381579,-.0328947){19}{\line(1,0){.0381579}}
\multiput(31.63,43.43)(.0381579,-.0328947){19}{\line(1,0){.0381579}}
\multiput(33.08,42.18)(.0381579,-.0328947){19}{\line(1,0){.0381579}}
\multiput(34.53,40.93)(.0381579,-.0328947){19}{\line(1,0){.0381579}}
\multiput(35.98,39.68)(.0381579,-.0328947){19}{\line(1,0){.0381579}}
\multiput(37.43,38.43)(.0381579,-.0328947){19}{\line(1,0){.0381579}}
\multiput(38.88,37.18)(.0381579,-.0328947){19}{\line(1,0){.0381579}}
\multiput(40.33,35.93)(.0381579,-.0328947){19}{\line(1,0){.0381579}}
\multiput(41.78,34.68)(.0381579,-.0328947){19}{\line(1,0){.0381579}}
\multiput(43.23,33.43)(.0381579,-.0328947){19}{\line(1,0){.0381579}}
%\end
%\dashline{1}(97.25,45)(111.75,32.5)
\multiput(97.18,44.93)(.0381579,-.0328947){19}{\line(1,0){.0381579}}
\multiput(98.63,43.68)(.0381579,-.0328947){19}{\line(1,0){.0381579}}
\multiput(100.08,42.43)(.0381579,-.0328947){19}{\line(1,0){.0381579}}
\multiput(101.53,41.18)(.0381579,-.0328947){19}{\line(1,0){.0381579}}
\multiput(102.98,39.93)(.0381579,-.0328947){19}{\line(1,0){.0381579}}
\multiput(104.43,38.68)(.0381579,-.0328947){19}{\line(1,0){.0381579}}
\multiput(105.88,37.43)(.0381579,-.0328947){19}{\line(1,0){.0381579}}
\multiput(107.33,36.18)(.0381579,-.0328947){19}{\line(1,0){.0381579}}
\multiput(108.78,34.93)(.0381579,-.0328947){19}{\line(1,0){.0381579}}
\multiput(110.23,33.68)(.0381579,-.0328947){19}{\line(1,0){.0381579}}
%\end
%\dashline{1}(165,45.75)(179.5,33.25)
\multiput(164.93,45.68)(.0381579,-.0328947){19}{\line(1,0){.0381579}}
\multiput(166.38,44.43)(.0381579,-.0328947){19}{\line(1,0){.0381579}}
\multiput(167.83,43.18)(.0381579,-.0328947){19}{\line(1,0){.0381579}}
\multiput(169.28,41.93)(.0381579,-.0328947){19}{\line(1,0){.0381579}}
\multiput(170.73,40.68)(.0381579,-.0328947){19}{\line(1,0){.0381579}}
\multiput(172.18,39.43)(.0381579,-.0328947){19}{\line(1,0){.0381579}}
\multiput(173.63,38.18)(.0381579,-.0328947){19}{\line(1,0){.0381579}}
\multiput(175.08,36.93)(.0381579,-.0328947){19}{\line(1,0){.0381579}}
\multiput(176.53,35.68)(.0381579,-.0328947){19}{\line(1,0){.0381579}}
\multiput(177.98,34.43)(.0381579,-.0328947){19}{\line(1,0){.0381579}}
%\end
%\dashline{1}(97.25,45)(111.75,32.5)
\multiput(97.18,44.93)(.0381579,-.0328947){19}{\line(1,0){.0381579}}
\multiput(98.63,43.68)(.0381579,-.0328947){19}{\line(1,0){.0381579}}
\multiput(100.08,42.43)(.0381579,-.0328947){19}{\line(1,0){.0381579}}
\multiput(101.53,41.18)(.0381579,-.0328947){19}{\line(1,0){.0381579}}
\multiput(102.98,39.93)(.0381579,-.0328947){19}{\line(1,0){.0381579}}
\multiput(104.43,38.68)(.0381579,-.0328947){19}{\line(1,0){.0381579}}
\multiput(105.88,37.43)(.0381579,-.0328947){19}{\line(1,0){.0381579}}
\multiput(107.33,36.18)(.0381579,-.0328947){19}{\line(1,0){.0381579}}
\multiput(108.78,34.93)(.0381579,-.0328947){19}{\line(1,0){.0381579}}
\multiput(110.23,33.68)(.0381579,-.0328947){19}{\line(1,0){.0381579}}
%\end
%\dashline{1}(165,45.75)(179.5,33.25)
\multiput(164.93,45.68)(.0381579,-.0328947){19}{\line(1,0){.0381579}}
\multiput(166.38,44.43)(.0381579,-.0328947){19}{\line(1,0){.0381579}}
\multiput(167.83,43.18)(.0381579,-.0328947){19}{\line(1,0){.0381579}}
\multiput(169.28,41.93)(.0381579,-.0328947){19}{\line(1,0){.0381579}}
\multiput(170.73,40.68)(.0381579,-.0328947){19}{\line(1,0){.0381579}}
\multiput(172.18,39.43)(.0381579,-.0328947){19}{\line(1,0){.0381579}}
\multiput(173.63,38.18)(.0381579,-.0328947){19}{\line(1,0){.0381579}}
\multiput(175.08,36.93)(.0381579,-.0328947){19}{\line(1,0){.0381579}}
\multiput(176.53,35.68)(.0381579,-.0328947){19}{\line(1,0){.0381579}}
\multiput(177.98,34.43)(.0381579,-.0328947){19}{\line(1,0){.0381579}}
%\end
%\dashline{1}(97.25,45)(111.75,32.5)
\multiput(97.18,44.93)(.0381579,-.0328947){19}{\line(1,0){.0381579}}
\multiput(98.63,43.68)(.0381579,-.0328947){19}{\line(1,0){.0381579}}
\multiput(100.08,42.43)(.0381579,-.0328947){19}{\line(1,0){.0381579}}
\multiput(101.53,41.18)(.0381579,-.0328947){19}{\line(1,0){.0381579}}
\multiput(102.98,39.93)(.0381579,-.0328947){19}{\line(1,0){.0381579}}
\multiput(104.43,38.68)(.0381579,-.0328947){19}{\line(1,0){.0381579}}
\multiput(105.88,37.43)(.0381579,-.0328947){19}{\line(1,0){.0381579}}
\multiput(107.33,36.18)(.0381579,-.0328947){19}{\line(1,0){.0381579}}
\multiput(108.78,34.93)(.0381579,-.0328947){19}{\line(1,0){.0381579}}
\multiput(110.23,33.68)(.0381579,-.0328947){19}{\line(1,0){.0381579}}
%\end
%\dashline{1}(165,45.75)(179.5,33.25)
\multiput(164.93,45.68)(.0381579,-.0328947){19}{\line(1,0){.0381579}}
\multiput(166.38,44.43)(.0381579,-.0328947){19}{\line(1,0){.0381579}}
\multiput(167.83,43.18)(.0381579,-.0328947){19}{\line(1,0){.0381579}}
\multiput(169.28,41.93)(.0381579,-.0328947){19}{\line(1,0){.0381579}}
\multiput(170.73,40.68)(.0381579,-.0328947){19}{\line(1,0){.0381579}}
\multiput(172.18,39.43)(.0381579,-.0328947){19}{\line(1,0){.0381579}}
\multiput(173.63,38.18)(.0381579,-.0328947){19}{\line(1,0){.0381579}}
\multiput(175.08,36.93)(.0381579,-.0328947){19}{\line(1,0){.0381579}}
\multiput(176.53,35.68)(.0381579,-.0328947){19}{\line(1,0){.0381579}}
\multiput(177.98,34.43)(.0381579,-.0328947){19}{\line(1,0){.0381579}}
%\end
%\dashline{1}(69.75,45)(55.25,32.5)
\multiput(69.68,44.93)(-.0381579,-.0328947){19}{\line(-1,0){.0381579}}
\multiput(68.23,43.68)(-.0381579,-.0328947){19}{\line(-1,0){.0381579}}
\multiput(66.78,42.43)(-.0381579,-.0328947){19}{\line(-1,0){.0381579}}
\multiput(65.33,41.18)(-.0381579,-.0328947){19}{\line(-1,0){.0381579}}
\multiput(63.88,39.93)(-.0381579,-.0328947){19}{\line(-1,0){.0381579}}
\multiput(62.43,38.68)(-.0381579,-.0328947){19}{\line(-1,0){.0381579}}
\multiput(60.98,37.43)(-.0381579,-.0328947){19}{\line(-1,0){.0381579}}
\multiput(59.53,36.18)(-.0381579,-.0328947){19}{\line(-1,0){.0381579}}
\multiput(58.08,34.93)(-.0381579,-.0328947){19}{\line(-1,0){.0381579}}
\multiput(56.63,33.68)(-.0381579,-.0328947){19}{\line(-1,0){.0381579}}
%\end
%\dashline{1}(137.5,45.75)(123,33.25)
\multiput(137.43,45.68)(-.0381579,-.0328947){19}{\line(-1,0){.0381579}}
\multiput(135.98,44.43)(-.0381579,-.0328947){19}{\line(-1,0){.0381579}}
\multiput(134.53,43.18)(-.0381579,-.0328947){19}{\line(-1,0){.0381579}}
\multiput(133.08,41.93)(-.0381579,-.0328947){19}{\line(-1,0){.0381579}}
\multiput(131.63,40.68)(-.0381579,-.0328947){19}{\line(-1,0){.0381579}}
\multiput(130.18,39.43)(-.0381579,-.0328947){19}{\line(-1,0){.0381579}}
\multiput(128.73,38.18)(-.0381579,-.0328947){19}{\line(-1,0){.0381579}}
\multiput(127.28,36.93)(-.0381579,-.0328947){19}{\line(-1,0){.0381579}}
\multiput(125.83,35.68)(-.0381579,-.0328947){19}{\line(-1,0){.0381579}}
\multiput(124.38,34.43)(-.0381579,-.0328947){19}{\line(-1,0){.0381579}}
%\end
%\dashline{1}(3,32)(3,32)
\put(2.93,31.93){\line(0,1){0}}
%\end
%\dashline{1}(44.75,32.25)(44.75,32.25)
\put(44.68,32.18){\line(0,1){0}}
%\end
%\dashline{1}(113.25,32.25)(113.25,32.25)
\put(113.18,32.18){\line(0,1){0}}
%\end
%\dashline{1}(181,33)(181,33)
\put(180.93,32.93){\line(0,1){0}}
%\end
%\dashline{1}(113.25,32.25)(113.25,32.25)
\put(113.18,32.18){\line(0,1){0}}
%\end
%\dashline{1}(181,33)(181,33)
\put(180.93,32.93){\line(0,1){0}}
%\end
%\dashline{1}(113.25,32.25)(113.25,32.25)
\put(113.18,32.18){\line(0,1){0}}
%\end
%\dashline{1}(181,33)(181,33)
\put(180.93,32.93){\line(0,1){0}}
%\end
%\dashline{1}(53.75,32.25)(53.75,32.25)
\put(53.68,32.18){\line(0,1){0}}
%\end
%\dashline{1}(121.5,33)(121.5,33)
\put(121.43,32.93){\line(0,1){0}}
%\end
%\dashline{1}(16.25,38.5)(4.5,28.25)
\multiput(16.18,38.43)(-.0383987,-.0334967){18}{\line(-1,0){.0383987}}
\multiput(14.797,37.224)(-.0383987,-.0334967){18}{\line(-1,0){.0383987}}
\multiput(13.415,36.018)(-.0383987,-.0334967){18}{\line(-1,0){.0383987}}
\multiput(12.033,34.812)(-.0383987,-.0334967){18}{\line(-1,0){.0383987}}
\multiput(10.65,33.606)(-.0383987,-.0334967){18}{\line(-1,0){.0383987}}
\multiput(9.268,32.4)(-.0383987,-.0334967){18}{\line(-1,0){.0383987}}
\multiput(7.886,31.194)(-.0383987,-.0334967){18}{\line(-1,0){.0383987}}
\multiput(6.503,29.989)(-.0383987,-.0334967){18}{\line(-1,0){.0383987}}
\multiput(5.121,28.783)(-.0383987,-.0334967){18}{\line(-1,0){.0383987}}
%\end
%\dashline{1}(31.5,38.75)(43.25,28.5)
\multiput(31.43,38.68)(.0383987,-.0334967){18}{\line(1,0){.0383987}}
\multiput(32.812,37.474)(.0383987,-.0334967){18}{\line(1,0){.0383987}}
\multiput(34.194,36.268)(.0383987,-.0334967){18}{\line(1,0){.0383987}}
\multiput(35.577,35.062)(.0383987,-.0334967){18}{\line(1,0){.0383987}}
\multiput(36.959,33.856)(.0383987,-.0334967){18}{\line(1,0){.0383987}}
\multiput(38.341,32.65)(.0383987,-.0334967){18}{\line(1,0){.0383987}}
\multiput(39.724,31.444)(.0383987,-.0334967){18}{\line(1,0){.0383987}}
\multiput(41.106,30.239)(.0383987,-.0334967){18}{\line(1,0){.0383987}}
\multiput(42.489,29.033)(.0383987,-.0334967){18}{\line(1,0){.0383987}}
%\end
%\dashline{1}(98.5,39)(110.25,28.75)
\multiput(98.43,38.93)(.0383987,-.0334967){18}{\line(1,0){.0383987}}
\multiput(99.812,37.724)(.0383987,-.0334967){18}{\line(1,0){.0383987}}
\multiput(101.194,36.518)(.0383987,-.0334967){18}{\line(1,0){.0383987}}
\multiput(102.577,35.312)(.0383987,-.0334967){18}{\line(1,0){.0383987}}
\multiput(103.959,34.106)(.0383987,-.0334967){18}{\line(1,0){.0383987}}
\multiput(105.341,32.9)(.0383987,-.0334967){18}{\line(1,0){.0383987}}
\multiput(106.724,31.694)(.0383987,-.0334967){18}{\line(1,0){.0383987}}
\multiput(108.106,30.489)(.0383987,-.0334967){18}{\line(1,0){.0383987}}
\multiput(109.489,29.283)(.0383987,-.0334967){18}{\line(1,0){.0383987}}
%\end
%\dashline{1}(166.25,39.75)(178,29.5)
\multiput(166.18,39.68)(.0383987,-.0334967){18}{\line(1,0){.0383987}}
\multiput(167.562,38.474)(.0383987,-.0334967){18}{\line(1,0){.0383987}}
\multiput(168.944,37.268)(.0383987,-.0334967){18}{\line(1,0){.0383987}}
\multiput(170.327,36.062)(.0383987,-.0334967){18}{\line(1,0){.0383987}}
\multiput(171.709,34.856)(.0383987,-.0334967){18}{\line(1,0){.0383987}}
\multiput(173.091,33.65)(.0383987,-.0334967){18}{\line(1,0){.0383987}}
\multiput(174.474,32.444)(.0383987,-.0334967){18}{\line(1,0){.0383987}}
\multiput(175.856,31.239)(.0383987,-.0334967){18}{\line(1,0){.0383987}}
\multiput(177.239,30.033)(.0383987,-.0334967){18}{\line(1,0){.0383987}}
%\end
%\dashline{1}(98.5,39)(110.25,28.75)
\multiput(98.43,38.93)(.0383987,-.0334967){18}{\line(1,0){.0383987}}
\multiput(99.812,37.724)(.0383987,-.0334967){18}{\line(1,0){.0383987}}
\multiput(101.194,36.518)(.0383987,-.0334967){18}{\line(1,0){.0383987}}
\multiput(102.577,35.312)(.0383987,-.0334967){18}{\line(1,0){.0383987}}
\multiput(103.959,34.106)(.0383987,-.0334967){18}{\line(1,0){.0383987}}
\multiput(105.341,32.9)(.0383987,-.0334967){18}{\line(1,0){.0383987}}
\multiput(106.724,31.694)(.0383987,-.0334967){18}{\line(1,0){.0383987}}
\multiput(108.106,30.489)(.0383987,-.0334967){18}{\line(1,0){.0383987}}
\multiput(109.489,29.283)(.0383987,-.0334967){18}{\line(1,0){.0383987}}
%\end
%\dashline{1}(166.25,39.75)(178,29.5)
\multiput(166.18,39.68)(.0383987,-.0334967){18}{\line(1,0){.0383987}}
\multiput(167.562,38.474)(.0383987,-.0334967){18}{\line(1,0){.0383987}}
\multiput(168.944,37.268)(.0383987,-.0334967){18}{\line(1,0){.0383987}}
\multiput(170.327,36.062)(.0383987,-.0334967){18}{\line(1,0){.0383987}}
\multiput(171.709,34.856)(.0383987,-.0334967){18}{\line(1,0){.0383987}}
\multiput(173.091,33.65)(.0383987,-.0334967){18}{\line(1,0){.0383987}}
\multiput(174.474,32.444)(.0383987,-.0334967){18}{\line(1,0){.0383987}}
\multiput(175.856,31.239)(.0383987,-.0334967){18}{\line(1,0){.0383987}}
\multiput(177.239,30.033)(.0383987,-.0334967){18}{\line(1,0){.0383987}}
%\end
%\dashline{1}(98.5,39)(110.25,28.75)
\multiput(98.43,38.93)(.0383987,-.0334967){18}{\line(1,0){.0383987}}
\multiput(99.812,37.724)(.0383987,-.0334967){18}{\line(1,0){.0383987}}
\multiput(101.194,36.518)(.0383987,-.0334967){18}{\line(1,0){.0383987}}
\multiput(102.577,35.312)(.0383987,-.0334967){18}{\line(1,0){.0383987}}
\multiput(103.959,34.106)(.0383987,-.0334967){18}{\line(1,0){.0383987}}
\multiput(105.341,32.9)(.0383987,-.0334967){18}{\line(1,0){.0383987}}
\multiput(106.724,31.694)(.0383987,-.0334967){18}{\line(1,0){.0383987}}
\multiput(108.106,30.489)(.0383987,-.0334967){18}{\line(1,0){.0383987}}
\multiput(109.489,29.283)(.0383987,-.0334967){18}{\line(1,0){.0383987}}
%\end
%\dashline{1}(166.25,39.75)(178,29.5)
\multiput(166.18,39.68)(.0383987,-.0334967){18}{\line(1,0){.0383987}}
\multiput(167.562,38.474)(.0383987,-.0334967){18}{\line(1,0){.0383987}}
\multiput(168.944,37.268)(.0383987,-.0334967){18}{\line(1,0){.0383987}}
\multiput(170.327,36.062)(.0383987,-.0334967){18}{\line(1,0){.0383987}}
\multiput(171.709,34.856)(.0383987,-.0334967){18}{\line(1,0){.0383987}}
\multiput(173.091,33.65)(.0383987,-.0334967){18}{\line(1,0){.0383987}}
\multiput(174.474,32.444)(.0383987,-.0334967){18}{\line(1,0){.0383987}}
\multiput(175.856,31.239)(.0383987,-.0334967){18}{\line(1,0){.0383987}}
\multiput(177.239,30.033)(.0383987,-.0334967){18}{\line(1,0){.0383987}}
%\end
%\dashline{1}(68.5,39)(56.75,28.75)
\multiput(68.43,38.93)(-.0383987,-.0334967){18}{\line(-1,0){.0383987}}
\multiput(67.047,37.724)(-.0383987,-.0334967){18}{\line(-1,0){.0383987}}
\multiput(65.665,36.518)(-.0383987,-.0334967){18}{\line(-1,0){.0383987}}
\multiput(64.283,35.312)(-.0383987,-.0334967){18}{\line(-1,0){.0383987}}
\multiput(62.9,34.106)(-.0383987,-.0334967){18}{\line(-1,0){.0383987}}
\multiput(61.518,32.9)(-.0383987,-.0334967){18}{\line(-1,0){.0383987}}
\multiput(60.136,31.694)(-.0383987,-.0334967){18}{\line(-1,0){.0383987}}
\multiput(58.753,30.489)(-.0383987,-.0334967){18}{\line(-1,0){.0383987}}
\multiput(57.371,29.283)(-.0383987,-.0334967){18}{\line(-1,0){.0383987}}
%\end
%\dashline{1}(136.25,39.75)(124.5,29.5)
\multiput(136.18,39.68)(-.0383987,-.0334967){18}{\line(-1,0){.0383987}}
\multiput(134.797,38.474)(-.0383987,-.0334967){18}{\line(-1,0){.0383987}}
\multiput(133.415,37.268)(-.0383987,-.0334967){18}{\line(-1,0){.0383987}}
\multiput(132.033,36.062)(-.0383987,-.0334967){18}{\line(-1,0){.0383987}}
\multiput(130.65,34.856)(-.0383987,-.0334967){18}{\line(-1,0){.0383987}}
\multiput(129.268,33.65)(-.0383987,-.0334967){18}{\line(-1,0){.0383987}}
\multiput(127.886,32.444)(-.0383987,-.0334967){18}{\line(-1,0){.0383987}}
\multiput(126.503,31.239)(-.0383987,-.0334967){18}{\line(-1,0){.0383987}}
\multiput(125.121,30.033)(-.0383987,-.0334967){18}{\line(-1,0){.0383987}}
%\end
%\dashline{1}(15.75,32.75)(7,24.75)
\multiput(15.68,32.68)(-.0354251,-.0323887){19}{\line(-1,0){.0354251}}
\multiput(14.334,31.449)(-.0354251,-.0323887){19}{\line(-1,0){.0354251}}
\multiput(12.987,30.218)(-.0354251,-.0323887){19}{\line(-1,0){.0354251}}
\multiput(11.641,28.987)(-.0354251,-.0323887){19}{\line(-1,0){.0354251}}
\multiput(10.295,27.757)(-.0354251,-.0323887){19}{\line(-1,0){.0354251}}
\multiput(8.949,26.526)(-.0354251,-.0323887){19}{\line(-1,0){.0354251}}
\multiput(7.603,25.295)(-.0354251,-.0323887){19}{\line(-1,0){.0354251}}
%\end
%\dashline{1}(32,33)(40.75,25)
\multiput(31.93,32.93)(.0354251,-.0323887){19}{\line(1,0){.0354251}}
\multiput(33.276,31.699)(.0354251,-.0323887){19}{\line(1,0){.0354251}}
\multiput(34.622,30.468)(.0354251,-.0323887){19}{\line(1,0){.0354251}}
\multiput(35.968,29.237)(.0354251,-.0323887){19}{\line(1,0){.0354251}}
\multiput(37.314,28.007)(.0354251,-.0323887){19}{\line(1,0){.0354251}}
\multiput(38.66,26.776)(.0354251,-.0323887){19}{\line(1,0){.0354251}}
\multiput(40.007,25.545)(.0354251,-.0323887){19}{\line(1,0){.0354251}}
%\end
%\dashline{1}(99,33.25)(107.75,25.25)
\multiput(98.93,33.18)(.0354251,-.0323887){19}{\line(1,0){.0354251}}
\multiput(100.276,31.949)(.0354251,-.0323887){19}{\line(1,0){.0354251}}
\multiput(101.622,30.718)(.0354251,-.0323887){19}{\line(1,0){.0354251}}
\multiput(102.968,29.487)(.0354251,-.0323887){19}{\line(1,0){.0354251}}
\multiput(104.314,28.257)(.0354251,-.0323887){19}{\line(1,0){.0354251}}
\multiput(105.66,27.026)(.0354251,-.0323887){19}{\line(1,0){.0354251}}
\multiput(107.007,25.795)(.0354251,-.0323887){19}{\line(1,0){.0354251}}
%\end
%\dashline{1}(166.75,34)(175.5,26)
\multiput(166.68,33.93)(.0354251,-.0323887){19}{\line(1,0){.0354251}}
\multiput(168.026,32.699)(.0354251,-.0323887){19}{\line(1,0){.0354251}}
\multiput(169.372,31.468)(.0354251,-.0323887){19}{\line(1,0){.0354251}}
\multiput(170.718,30.237)(.0354251,-.0323887){19}{\line(1,0){.0354251}}
\multiput(172.064,29.007)(.0354251,-.0323887){19}{\line(1,0){.0354251}}
\multiput(173.41,27.776)(.0354251,-.0323887){19}{\line(1,0){.0354251}}
\multiput(174.757,26.545)(.0354251,-.0323887){19}{\line(1,0){.0354251}}
%\end
%\dashline{1}(99,33.25)(107.75,25.25)
\multiput(98.93,33.18)(.0354251,-.0323887){19}{\line(1,0){.0354251}}
\multiput(100.276,31.949)(.0354251,-.0323887){19}{\line(1,0){.0354251}}
\multiput(101.622,30.718)(.0354251,-.0323887){19}{\line(1,0){.0354251}}
\multiput(102.968,29.487)(.0354251,-.0323887){19}{\line(1,0){.0354251}}
\multiput(104.314,28.257)(.0354251,-.0323887){19}{\line(1,0){.0354251}}
\multiput(105.66,27.026)(.0354251,-.0323887){19}{\line(1,0){.0354251}}
\multiput(107.007,25.795)(.0354251,-.0323887){19}{\line(1,0){.0354251}}
%\end
%\dashline{1}(166.75,34)(175.5,26)
\multiput(166.68,33.93)(.0354251,-.0323887){19}{\line(1,0){.0354251}}
\multiput(168.026,32.699)(.0354251,-.0323887){19}{\line(1,0){.0354251}}
\multiput(169.372,31.468)(.0354251,-.0323887){19}{\line(1,0){.0354251}}
\multiput(170.718,30.237)(.0354251,-.0323887){19}{\line(1,0){.0354251}}
\multiput(172.064,29.007)(.0354251,-.0323887){19}{\line(1,0){.0354251}}
\multiput(173.41,27.776)(.0354251,-.0323887){19}{\line(1,0){.0354251}}
\multiput(174.757,26.545)(.0354251,-.0323887){19}{\line(1,0){.0354251}}
%\end
%\dashline{1}(99,33.25)(107.75,25.25)
\multiput(98.93,33.18)(.0354251,-.0323887){19}{\line(1,0){.0354251}}
\multiput(100.276,31.949)(.0354251,-.0323887){19}{\line(1,0){.0354251}}
\multiput(101.622,30.718)(.0354251,-.0323887){19}{\line(1,0){.0354251}}
\multiput(102.968,29.487)(.0354251,-.0323887){19}{\line(1,0){.0354251}}
\multiput(104.314,28.257)(.0354251,-.0323887){19}{\line(1,0){.0354251}}
\multiput(105.66,27.026)(.0354251,-.0323887){19}{\line(1,0){.0354251}}
\multiput(107.007,25.795)(.0354251,-.0323887){19}{\line(1,0){.0354251}}
%\end
%\dashline{1}(166.75,34)(175.5,26)
\multiput(166.68,33.93)(.0354251,-.0323887){19}{\line(1,0){.0354251}}
\multiput(168.026,32.699)(.0354251,-.0323887){19}{\line(1,0){.0354251}}
\multiput(169.372,31.468)(.0354251,-.0323887){19}{\line(1,0){.0354251}}
\multiput(170.718,30.237)(.0354251,-.0323887){19}{\line(1,0){.0354251}}
\multiput(172.064,29.007)(.0354251,-.0323887){19}{\line(1,0){.0354251}}
\multiput(173.41,27.776)(.0354251,-.0323887){19}{\line(1,0){.0354251}}
\multiput(174.757,26.545)(.0354251,-.0323887){19}{\line(1,0){.0354251}}
%\end
%\dashline{1}(68,33.25)(59.25,25.25)
\multiput(67.93,33.18)(-.0354251,-.0323887){19}{\line(-1,0){.0354251}}
\multiput(66.584,31.949)(-.0354251,-.0323887){19}{\line(-1,0){.0354251}}
\multiput(65.237,30.718)(-.0354251,-.0323887){19}{\line(-1,0){.0354251}}
\multiput(63.891,29.487)(-.0354251,-.0323887){19}{\line(-1,0){.0354251}}
\multiput(62.545,28.257)(-.0354251,-.0323887){19}{\line(-1,0){.0354251}}
\multiput(61.199,27.026)(-.0354251,-.0323887){19}{\line(-1,0){.0354251}}
\multiput(59.853,25.795)(-.0354251,-.0323887){19}{\line(-1,0){.0354251}}
%\end
%\dashline{1}(135.75,34)(127,26)
\multiput(135.68,33.93)(-.0354251,-.0323887){19}{\line(-1,0){.0354251}}
\multiput(134.334,32.699)(-.0354251,-.0323887){19}{\line(-1,0){.0354251}}
\multiput(132.987,31.468)(-.0354251,-.0323887){19}{\line(-1,0){.0354251}}
\multiput(131.641,30.237)(-.0354251,-.0323887){19}{\line(-1,0){.0354251}}
\multiput(130.295,29.007)(-.0354251,-.0323887){19}{\line(-1,0){.0354251}}
\multiput(128.949,27.776)(-.0354251,-.0323887){19}{\line(-1,0){.0354251}}
\multiput(127.603,26.545)(-.0354251,-.0323887){19}{\line(-1,0){.0354251}}
%\end
%\dashline{1}(16.25,27.25)(10.75,21.5)
\multiput(16.18,27.18)(-.0321637,-.0336257){19}{\line(0,-1){.0336257}}
\multiput(14.957,25.902)(-.0321637,-.0336257){19}{\line(0,-1){.0336257}}
\multiput(13.735,24.624)(-.0321637,-.0336257){19}{\line(0,-1){.0336257}}
\multiput(12.513,23.346)(-.0321637,-.0336257){19}{\line(0,-1){.0336257}}
\multiput(11.291,22.069)(-.0321637,-.0336257){19}{\line(0,-1){.0336257}}
%\end
%\dashline{1}(31.5,27.5)(37,21.75)
\multiput(31.43,27.43)(.0321637,-.0336257){19}{\line(0,-1){.0336257}}
\multiput(32.652,26.152)(.0321637,-.0336257){19}{\line(0,-1){.0336257}}
\multiput(33.874,24.874)(.0321637,-.0336257){19}{\line(0,-1){.0336257}}
\multiput(35.096,23.596)(.0321637,-.0336257){19}{\line(0,-1){.0336257}}
\multiput(36.319,22.319)(.0321637,-.0336257){19}{\line(0,-1){.0336257}}
%\end
%\dashline{1}(98.5,27.75)(104,22)
\multiput(98.43,27.68)(.0321637,-.0336257){19}{\line(0,-1){.0336257}}
\multiput(99.652,26.402)(.0321637,-.0336257){19}{\line(0,-1){.0336257}}
\multiput(100.874,25.124)(.0321637,-.0336257){19}{\line(0,-1){.0336257}}
\multiput(102.096,23.846)(.0321637,-.0336257){19}{\line(0,-1){.0336257}}
\multiput(103.319,22.569)(.0321637,-.0336257){19}{\line(0,-1){.0336257}}
%\end
%\dashline{1}(166.25,28.5)(171.75,22.75)
\multiput(166.18,28.43)(.0321637,-.0336257){19}{\line(0,-1){.0336257}}
\multiput(167.402,27.152)(.0321637,-.0336257){19}{\line(0,-1){.0336257}}
\multiput(168.624,25.874)(.0321637,-.0336257){19}{\line(0,-1){.0336257}}
\multiput(169.846,24.596)(.0321637,-.0336257){19}{\line(0,-1){.0336257}}
\multiput(171.069,23.319)(.0321637,-.0336257){19}{\line(0,-1){.0336257}}
%\end
%\dashline{1}(98.5,27.75)(104,22)
\multiput(98.43,27.68)(.0321637,-.0336257){19}{\line(0,-1){.0336257}}
\multiput(99.652,26.402)(.0321637,-.0336257){19}{\line(0,-1){.0336257}}
\multiput(100.874,25.124)(.0321637,-.0336257){19}{\line(0,-1){.0336257}}
\multiput(102.096,23.846)(.0321637,-.0336257){19}{\line(0,-1){.0336257}}
\multiput(103.319,22.569)(.0321637,-.0336257){19}{\line(0,-1){.0336257}}
%\end
%\dashline{1}(166.25,28.5)(171.75,22.75)
\multiput(166.18,28.43)(.0321637,-.0336257){19}{\line(0,-1){.0336257}}
\multiput(167.402,27.152)(.0321637,-.0336257){19}{\line(0,-1){.0336257}}
\multiput(168.624,25.874)(.0321637,-.0336257){19}{\line(0,-1){.0336257}}
\multiput(169.846,24.596)(.0321637,-.0336257){19}{\line(0,-1){.0336257}}
\multiput(171.069,23.319)(.0321637,-.0336257){19}{\line(0,-1){.0336257}}
%\end
%\dashline{1}(98.5,27.75)(104,22)
\multiput(98.43,27.68)(.0321637,-.0336257){19}{\line(0,-1){.0336257}}
\multiput(99.652,26.402)(.0321637,-.0336257){19}{\line(0,-1){.0336257}}
\multiput(100.874,25.124)(.0321637,-.0336257){19}{\line(0,-1){.0336257}}
\multiput(102.096,23.846)(.0321637,-.0336257){19}{\line(0,-1){.0336257}}
\multiput(103.319,22.569)(.0321637,-.0336257){19}{\line(0,-1){.0336257}}
%\end
%\dashline{1}(166.25,28.5)(171.75,22.75)
\multiput(166.18,28.43)(.0321637,-.0336257){19}{\line(0,-1){.0336257}}
\multiput(167.402,27.152)(.0321637,-.0336257){19}{\line(0,-1){.0336257}}
\multiput(168.624,25.874)(.0321637,-.0336257){19}{\line(0,-1){.0336257}}
\multiput(169.846,24.596)(.0321637,-.0336257){19}{\line(0,-1){.0336257}}
\multiput(171.069,23.319)(.0321637,-.0336257){19}{\line(0,-1){.0336257}}
%\end
%\dashline{1}(68.5,27.75)(63,22)
\multiput(68.43,27.68)(-.0321637,-.0336257){19}{\line(0,-1){.0336257}}
\multiput(67.207,26.402)(-.0321637,-.0336257){19}{\line(0,-1){.0336257}}
\multiput(65.985,25.124)(-.0321637,-.0336257){19}{\line(0,-1){.0336257}}
\multiput(64.763,23.846)(-.0321637,-.0336257){19}{\line(0,-1){.0336257}}
\multiput(63.541,22.569)(-.0321637,-.0336257){19}{\line(0,-1){.0336257}}
%\end
%\dashline{1}(136.25,28.5)(130.75,22.75)
\multiput(136.18,28.43)(-.0321637,-.0336257){19}{\line(0,-1){.0336257}}
\multiput(134.957,27.152)(-.0321637,-.0336257){19}{\line(0,-1){.0336257}}
\multiput(133.735,25.874)(-.0321637,-.0336257){19}{\line(0,-1){.0336257}}
\multiput(132.513,24.596)(-.0321637,-.0336257){19}{\line(0,-1){.0336257}}
\multiput(131.291,23.319)(-.0321637,-.0336257){19}{\line(0,-1){.0336257}}
%\end
%\dashline{1}(10.75,21.5)(10.75,21.5)
\put(10.68,21.43){\line(0,1){0}}
%\end
%\dashline{1}(37,21.75)(37,21.75)
\put(36.93,21.68){\line(0,1){0}}
%\end
%\dashline{1}(116.25,21.75)(116.25,21.75)
\put(116.18,21.68){\line(0,1){0}}
%\end
%\dashline{1}(184,22.5)(184,22.5)
\put(183.93,22.43){\line(0,1){0}}
%\end
%\dashline{1}(116.25,21.75)(116.25,21.75)
\put(116.18,21.68){\line(0,1){0}}
%\end
%\dashline{1}(184,22.5)(184,22.5)
\put(183.93,22.43){\line(0,1){0}}
%\end
%\dashline{1}(116.25,21.75)(116.25,21.75)
\put(116.18,21.68){\line(0,1){0}}
%\end
%\dashline{1}(184,22.5)(184,22.5)
\put(183.93,22.43){\line(0,1){0}}
%\end
%\dashline{1}(50.75,21.75)(50.75,21.75)
\put(50.68,21.68){\line(0,1){0}}
%\end
%\dashline{1}(118.5,22.5)(118.5,22.5)
\put(118.43,22.43){\line(0,1){0}}
%\end
%\dashline{1}(16.75,22.75)(14.25,19.5)
\multiput(16.68,22.68)(-.0333333,-.0433333){15}{\line(0,-1){.0433333}}
\multiput(15.68,21.38)(-.0333333,-.0433333){15}{\line(0,-1){.0433333}}
\multiput(14.68,20.08)(-.0333333,-.0433333){15}{\line(0,-1){.0433333}}
%\end
%\dashline{1}(31,23)(33.5,19.75)
\multiput(30.93,22.93)(.0333333,-.0433333){15}{\line(0,-1){.0433333}}
\multiput(31.93,21.63)(.0333333,-.0433333){15}{\line(0,-1){.0433333}}
\multiput(32.93,20.33)(.0333333,-.0433333){15}{\line(0,-1){.0433333}}
%\end
%\dashline{1}(98,23.25)(100.5,20)
\multiput(97.93,23.18)(.0333333,-.0433333){15}{\line(0,-1){.0433333}}
\multiput(98.93,21.88)(.0333333,-.0433333){15}{\line(0,-1){.0433333}}
\multiput(99.93,20.58)(.0333333,-.0433333){15}{\line(0,-1){.0433333}}
%\end
%\dashline{1}(165.75,24)(168.25,20.75)
\multiput(165.68,23.93)(.0333333,-.0433333){15}{\line(0,-1){.0433333}}
\multiput(166.68,22.63)(.0333333,-.0433333){15}{\line(0,-1){.0433333}}
\multiput(167.68,21.33)(.0333333,-.0433333){15}{\line(0,-1){.0433333}}
%\end
%\dashline{1}(98,23.25)(100.5,20)
\multiput(97.93,23.18)(.0333333,-.0433333){15}{\line(0,-1){.0433333}}
\multiput(98.93,21.88)(.0333333,-.0433333){15}{\line(0,-1){.0433333}}
\multiput(99.93,20.58)(.0333333,-.0433333){15}{\line(0,-1){.0433333}}
%\end
%\dashline{1}(165.75,24)(168.25,20.75)
\multiput(165.68,23.93)(.0333333,-.0433333){15}{\line(0,-1){.0433333}}
\multiput(166.68,22.63)(.0333333,-.0433333){15}{\line(0,-1){.0433333}}
\multiput(167.68,21.33)(.0333333,-.0433333){15}{\line(0,-1){.0433333}}
%\end
%\dashline{1}(98,23.25)(100.5,20)
\multiput(97.93,23.18)(.0333333,-.0433333){15}{\line(0,-1){.0433333}}
\multiput(98.93,21.88)(.0333333,-.0433333){15}{\line(0,-1){.0433333}}
\multiput(99.93,20.58)(.0333333,-.0433333){15}{\line(0,-1){.0433333}}
%\end
%\dashline{1}(165.75,24)(168.25,20.75)
\multiput(165.68,23.93)(.0333333,-.0433333){15}{\line(0,-1){.0433333}}
\multiput(166.68,22.63)(.0333333,-.0433333){15}{\line(0,-1){.0433333}}
\multiput(167.68,21.33)(.0333333,-.0433333){15}{\line(0,-1){.0433333}}
%\end
%\dashline{1}(69,23.25)(66.5,20)
\multiput(68.93,23.18)(-.0333333,-.0433333){15}{\line(0,-1){.0433333}}
\multiput(67.93,21.88)(-.0333333,-.0433333){15}{\line(0,-1){.0433333}}
\multiput(66.93,20.58)(-.0333333,-.0433333){15}{\line(0,-1){.0433333}}
%\end
%\dashline{1}(136.75,24)(134.25,20.75)
\multiput(136.68,23.93)(-.0333333,-.0433333){15}{\line(0,-1){.0433333}}
\multiput(135.68,22.63)(-.0333333,-.0433333){15}{\line(0,-1){.0433333}}
\multiput(134.68,21.33)(-.0333333,-.0433333){15}{\line(0,-1){.0433333}}
%\end
%\dashline{1}(14.25,19.5)(14.25,19.5)
\put(14.18,19.43){\line(0,1){0}}
%\end
%\dashline{1}(33.5,19.75)(33.5,19.75)
\put(33.43,19.68){\line(0,1){0}}
%\end
%\dashline{1}(112.75,19.75)(112.75,19.75)
\put(112.68,19.68){\line(0,1){0}}
%\end
%\dashline{1}(180.5,20.5)(180.5,20.5)
\put(180.43,20.43){\line(0,1){0}}
%\end
%\dashline{1}(112.75,19.75)(112.75,19.75)
\put(112.68,19.68){\line(0,1){0}}
%\end
%\dashline{1}(180.5,20.5)(180.5,20.5)
\put(180.43,20.43){\line(0,1){0}}
%\end
%\dashline{1}(112.75,19.75)(112.75,19.75)
\put(112.68,19.68){\line(0,1){0}}
%\end
%\dashline{1}(180.5,20.5)(180.5,20.5)
\put(180.43,20.43){\line(0,1){0}}
%\end
%\dashline{1}(54.25,19.75)(54.25,19.75)
\put(54.18,19.68){\line(0,1){0}}
%\end
%\dashline{1}(122,20.5)(122,20.5)
\put(121.93,20.43){\line(0,1){0}}
%\end
\put(23,59.25){\makebox(0,0)[cc]{$u$}}
\put(83.5,59.25){\makebox(0,0)[cc]{$u$}}
\put(151.25,60){\makebox(0,0)[cc]{$u$}}
\put(22.5,7.75){\makebox(0,0)[cc]{$v$}}
\put(4.25,66.5){\makebox(0,0)[]{}}
\put(64.75,63.5){\makebox(0,0)[]{}}
\put(132.5,64.25){\makebox(0,0)[]{}}
\put(43.5,63.75){\makebox(0,0)[cc]{}}
\put(22.75,.75){\makebox(0,0)[cc]{(a) Graph $G$}}
\put(83.25,.75){\makebox(0,0)[cc]{(b) Graph $H$}}
\put(150.75,1.5){\makebox(0,0)[cc]{(c) Graph $H+v_1v_2$}}
\put(7.5,33){$G_1$}
\put(61.75,32.5){$G_1$}
\put(129.5,33.25){$G_1$}
\put(40.5,34.75){$G_2$}
\put(106.25,35){$G_2$}
\put(174,35.75){$G_2$}
\qbezier(83,54.25)(61.375,36.125)(73.25,14.5)
\qbezier(150.75,55)(129.125,36.875)(141,15.25)
\qbezier(93.75,14.5)(137.25,38)(84.75,54.5)
\qbezier(161.5,15.25)(205,38.75)(152.5,55.25)
\qbezier(93.75,14.5)(137.25,38)(84.75,54.5)
\qbezier(161.5,15.25)(205,38.75)(152.5,55.25)
\qbezier(93.75,14.5)(137.25,38)(84.75,54.5)
\qbezier(161.5,15.25)(205,38.75)(152.5,55.25)
\qbezier(73.25,14.5)(29.75,38)(82.25,54.5)
\qbezier(141,15.25)(97.5,38.75)(150,55.25)
\put(72.75,8.75){$v_1$}
\put(140.5,9.5){$v_1$}
\put(92.75,8.5){$v_2$}
\put(160.5,9.25){$v_2$}
\put(158.25,50.5){\line(0,1){0}}
\qbezier(84,54.25)(105.875,36.625)(93.25,14.5)
\qbezier(152,54.75)(173.875,37.125)(161.25,15)
\put(141.25,15){\line(1,0){20}}
\end{picture}
}

\caption{\relabel{f3} 
$G$ is formed by proper subgraphs $G_1$ and $G_2$,
and $H/v_1v_2=G$}
\end{figure}
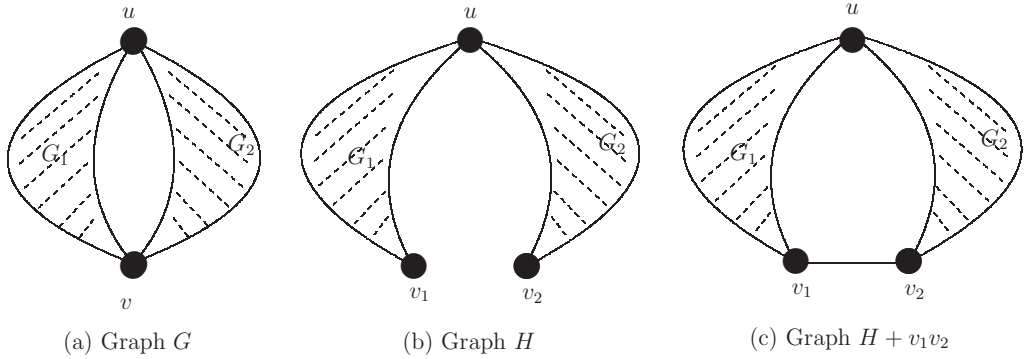

\proof 
Let $H$ be the graph obtained $G$ 
by replacing $v$ by two new vertices $v_1$ and $v_2$ 
and for all edges in $G_i$ incident with $v$, 
changing their common end $v$ to $v_i$,
as shown in Figure~\ref{f3}(b).
Thus $H/v_1v_2$ is the graph $G$. 
By (\ref{eq1-2}), we have 
\beeq
F(G,\lambda)=F(H,\lambda)+F(H+v_1v_2,\lambda).
\eneq
By Lemma~\ref{block-factor}, 
\beeq
F(H,\lambda)=F(G_1,\lambda)F(G_2,\lambda)
\eneq
and by Lemma~\ref{v-edge}, 
\beeq
F(H+v_1v_2,\lambda)=
\frac{F(G_1+uv,\lambda)F(G_2+uv,\lambda)}{\lambda-1}.
\eneq
Thus the result holds.
\proofend

\resection{A theorem on a zero-free interval}
\relabel{sec3}

In this section, we shall 
provide a sufficient condition for determining  
a zero-free interval $(1,\beta)$ of $F(G,\lambda)$, 
where $\beta\in (1,2)$, for 
all graphs $G$ in a family $\sets$.  
We shall first obtain a sufficient condition 
for a real number $\lambda$ in $(1,2)$ 
such that $F(G,\lambda)\ne 0$ for all graphs $G$ in $\sets$.  
In proving this result, we use 
some techniques that have appeared in~\cite{jac1} 
where Jackson proved that 
every chromatic polynomial has no real roots in $(1,32/27]$.
For any connected graph $G$, let
\beeq
Q(G,\lambda)=(-1)^{p(G)}F(G,\lambda)
\eneq
where  
$p(G)=|E(G)|-|V(G)|+b(G)-1$.
So $p(G)=|E(G)|-|V(G)|$ if $G$ is non-separable with $E(G)\ne \emptyset$. 
Theorem~\ref{Wakelin} implies that 
$Q(G,\lambda)>0$ for any bridgeless connected graph $G$
and real number $\lambda\in (1,32/27]$. 
It is also clear that $F(G,\lambda)\ne 0$ if and only if 
$Q(G,\lambda)\ne 0$.

\begin{lem}\quad
\relabel{general-interval-le}
Let $\sets$ be a family of bridgeless connected graphs and
$\lambda$ be any real number in $(1,2)$.
Assume that 
$\sets$ contains a subfamily $\sets'$ of non-separable
graphs  
such that conditions (i)-(iii) below are satisfied:
\begin{enumerate}
\item $Q(G,\lambda)>0$ for all graphs $G\in \sets'$;
%and all real $\lambda\in (1,\beta)$;
\item for every separable graph $G\in \sets$, 
all blocks of $G$ belong to $\sets$;
\item for every non-separable graph $G\in \sets\setminus \sets'$,
one of the following cases occurs:
\begin{enumerate}
\item \relabel{con-a}
for some edge $e$ in $G$,  
$G-e$ has a cut-vertex $u$ and each $G_i$ belongs to 
$\sets$ for $i=1,2$,
where $G_1$ and $G_2$ are graphs stated in Lemma~\ref{v-edge};

\item \relabel{con-b}
for some edge $e$ in $G$,  
both $G-e$ and $G/e$ belong to $\sets$ 
and both $b(G-e)$ and $b(G/e)$ are odd numbers;

\item \relabel{con-c}
there are subgraphs $G_1$ and $G_2$ of $G$ 
with $V(G_1)\cap V(G_2)=\{u_1,u_2\}$,
$V(G_1)\cup V(G_2)=V(G)$,
$E(G_1)\cap E(G_2)=\emptyset$
and $E(G_1)\cup E(G_2)=E(G)$, as shown in Figure~\ref{f3}(a),
such that $b(G_1)+b(G_2)$ is even,
and for $i=1,2$, $|E(G_i)|\ge 2$ and both 
$G_i+u_1u_2$ and $G_i$ belong to $\sets$,
where $G_i+u_1u_2$ is the graph obtained from $G_i$ by
adding a new edge joining $u_1$ and $u_2$; and

\item \relabel{con-d}
there are subgraphs $G_1$ and $G_2$ of $G$ 
with $|E(G_1)|\ge 3$, $|E(G_2)|\ge 2$, 
$V(G_1)\cap V(G_2)=\{u_1,u_2\}$,
$V(G_1)\cup V(G_2)=V(G)$,
$E(G_1)\cap E(G_2)=\emptyset$
and $E(G_1)\cup E(G_2)=E(G)$, as shown in Figure~\ref{f3}(a),
such that %$Q(G_2+u_1u_2,\lambda)\ge (\lambda-1)Q(G_2)$,
$b(G_1/u_1u_2)+b(G_2)$ is an odd number
and $G_1+u_1u_2$, $G_1/u_1u_2$, $G_2$, $G_2+u_1u_2$ and 
$G_2+2u_1u_2$ all belong to $\sets$, where 
$G_2+2u_1u_2$ is the graph obtained from $G_2$ by
adding two parallel edges joining $u_1$ and $u_2$.
\end{enumerate}
\end{enumerate}
Then 
$Q(G,\lambda)>0$ for all graphs $G\in \sets$.
%and all real $\lambda\in (1,\beta)$.
\end{lem}

\proof Suppose the result does not hold. 
Then there exists $G\in \sets$ 
such that $Q(G,\lambda)\le 0$ 
but $Q(H,\lambda)>0$ for all $H\in \sets$ with 
$|E(H)|<m$, where $m=|E(G)|$.
Now let $G$ be fixed.
By Condition (i), 
either $G$ is separable or $G\in \sets\setminus \sets'$.
We shall complete the proof by proving the following claims.

\inclaim $G$ is non-separable.

Suppose that $G$ is separable with blocks 
$G_1,G_2,\cdots, G_k$, 
where $k=b(G)\ge 2$.
For all $i=1,2,\cdots,k$, 
since $|E(G_i)|<m$ and $G_i\in \sets$ by Condition (ii),
we have $Q(G_i,\lambda)>0$. 
Note that 
\beeqn
p(G)
&=&|E(G)|-|V(G)|+k-1\\
&=&\sum_{i=1}^k |E(G_i)|
-\left (-(k-1)+\sum_{i=1}^k |V(G_i)|\right )+k-1\\
&=&2(k-1)+\sum_{i=1}^k p(G_i).
\eneqn
By Lemma~\ref{block-factor}, 
\beeq
F(G,\lambda) 
=\prod_{i=1}^k F(G_i,\lambda).
\eneq

Thus
\beeq
Q(G,\lambda)=(-1)^{p(G)}F(G,\lambda)
=\prod_{i=1}^k (-1)^{p(G_i)}F(G_i,\lambda)
=\prod_{i=1}^k Q(G_i,\lambda)>0,
\eneq
a contradiction. 
Hence Claim~\thecountclaim\ holds.

\inclaim Condition (a) of (iii) is not satisfied.

Suppose that 
$G$ contains an edge $e$ such that 
$G-e$ has a cut-vertex $u$ and $G_i\in \sets$ for $i=1,2$,
where $G_1$ and $G_2$ are graphs stated in Lemma~\ref{v-edge}.
As $|E(G_i)|<m$, we have $Q(G_i,\lambda)>0$ for $i=1,2$. 
By Lemma~\ref{v-edge},
\beeq
F(G,\lambda)=\frac{F(G_1,\lambda)F(G_2,\lambda)}{\lambda-1}.
\eneq
Since $G$ is non-separable by Claim 1, 
both $G_1$ and $G_2$ are non-separable. 
Thus
$$
p(G_1)+p(G_2)=|E(G_1)|-|V(G_1)|+|E(G_2)|-|V(G_2)|=
(|E(G)|+1)-(|V(G)|+1)=p(G),
$$
implying that 
\beeq
Q(G,\lambda)=\frac{Q(G_1,\lambda)Q(G_2,\lambda)}{\lambda-1}>0,
\eneq
a contradiction.
Hence Claim~\thecountclaim\ holds.

\inclaim Condition (b) of (iii) is not satisfied.

Suppose that $G$ contains an edge $e$ such that 
both $b(G/e)$ and $b(G-e)$ is odd and 
%all blocks of $G/e$ and all blocks of $G-e$ belong to $\sets$.
both $G/e$ and $G-e$ belong to $\sets$.

Note that 
\beeqn
p(G/e)&=&|E(G/e)|-|V(G/e)|+b(G/e)-1\\
&=&|E(G)|-1-(|V(G)|-1)+b(G/e)-1
=p(G)+b(G/e)-1
\eneqn
and 
\beeqn
p(G-e)&=&|E(G-e)|-|V(G-e)|+b(G-e)-1\\
&=&|E(G)|-1-|V(G)|+b(G-e)-1
=p(G)+b(G-e)-2.
\eneqn
As $G$ is non-separable, $e$ is not a loop.
By (\ref{eq1-2}), we have 
\beeqn
F(G,\lambda)=F(G/e,\lambda)-F(G-e,\lambda).
\eneqn
Since both $b(G/e)$ and $b(G-e)$ are odd, we have 
\beeqn
Q(G,\lambda)=Q(G/e,\lambda)+Q(G-e,\lambda).
\eneqn
Since both $G/e$ and $G-e$ belong to $\sets$ and 
both have less edges than $G$, 
by the assumption on $G$, we have $Q(G/e,\lambda)>0$ 
and $Q(G-e,\lambda)>0$.
Thus $Q(G,\lambda)>0$, a contradiction.
Hence Claim~\thecountclaim\ holds.

\inclaim Condition (c) of (iii) is not satisfied.

Suppose that condition (c) of (iii) is satisfied.
Let $G_1$ and $G_2$ be such 
subgraphs of $G$ stated in condition (c).
By Lemma~\ref{2-blocks}, 
%it can be shown that 
\beeq
F(G,\lambda)=
\frac{1}{\lambda-1}
F(G_1+u_1u_2,\lambda)F(G_2+u_1u_2,\lambda)
+F(G_1,\lambda)F(G_2,\lambda).
\eneq
As $G_i+u_1u_2$ is non-separable for $i=1,2$, we have 
\beeq
p(G_1+u_1u_2)+p(G_2+u_1u_2)
=|E(G_1)|+1-|V(G_1)|+|E(G_2)|+1-|V(G_2)|
=m-|V|=p(G).
\eneq
We also have
\beeqn
p(G_1)+p(G_2)
&=&|E(G_1)|-|V(G_1)|+b(G_1)-1+|E(G_2)|-|V(G_2)|+b(G_2)-1\\
&=&m-|V|-4+b(G_1)+b(G_2)
=p(G)-4+b(G_1)+b(G_2).
\eneqn
Since $b(G_1)+b(G_2)$ is even, 
\beeq
Q(G,\lambda)=
\frac{1}{\lambda-1}
Q(G_1+u_1u_2,\lambda)Q(G_2+u_1u_2,\lambda)
+Q(G_1,\lambda)Q(G_2,\lambda).
\eneq
As $|E(G_i)|\le m-2$, by the assumption $G$,
$Q(G_1+u_1u_2,\lambda)$, $Q(G_2+u_1u_2,\lambda)$,
$Q(G_1,\lambda)$ and $Q(G_2,\lambda)$ are all positive,
and so 
$Q(G,\lambda)>0$, a contradiction.

\inclaim Condition (d) of (iii) is not satisfied.

Suppose that condition (d) of (iii) is satisfied.
Assume that $G_1$ and $G_2$ are two subgraphs of $G$ 
as stated in condition (d),
as shown in Figure~\ref{f3}(a).
By Lemma~\ref{2-blocks}, 
\beeqn
& &F(G,\lambda)\\
&=&
\frac{1}{\lambda-1}
F(G_1+u_1u_2,\lambda)F(G_2+u_1u_2,\lambda)
+F(G_1,\lambda)F(G_2,\lambda)\\
&=&\frac{1}{\lambda-1}
F(G_1+u_1u_2,\lambda)F(G_2+u_1u_2,\lambda)
+[F(G_1/u_1u_2,\lambda)-F(G_1+u_1u_2,\lambda)]F(G_2,\lambda)\\
&=&F(G_1/u_1u_2,\lambda)F(G_2,\lambda)
+\frac{F(G_1+u_1u_2,\lambda)}{\lambda-1}
\left [ F(G_2+u_1u_2,\lambda)
-(\lambda-1)F(G_2,\lambda)\right ],
\eneqn 
and also by Lemma~\ref{2-blocks}, we have 
\beeq
F(G_2+2u_1u_2,\lambda)
=(\lambda-2)F(G_2+u_1u_2,\lambda)
+(\lambda-1)F(G_2,\lambda).
\eneq
Thus 
\beeq
F(G,\lambda)
=F(G_1/u_1u_2,\lambda)F(G_2,\lambda)
+F(G_1+u_1u_2,\lambda)
\left [F(G_2+u_1u_2,\lambda)
-\frac{F(G_2+2u_1u_2,\lambda)}{\lambda-1}\right ].
\eneq
Note that 
\beeqn
& &p(G_1/u_1u_2)+p(G_2)\\
&=&|E(G_1)|-(|V(G_1)|-1)+b(G_1/u_1u_2)-1
+|E(G_2)|-|V(G_2)|+b(G_2)-1\\
&=&|E(G)|-|V(G)|-3+b(G_1/u_1u_2)+b(G_2)\\
&=&p(G)-3+b(G_1/u_1u_2)+b(G_2).
\eneqn
As $b(G_1/u_1u_2)+b(G_2)$ is an odd number, 
$p(G_1/u_1u_2)+p(G_2)$ and $p(G)$ 
have the same parity (i.e., the sum of them is even).
%\footnote{Two integers are said to have the same parity if either both are odd numbers or both are even numbers}.
It can also be checked similarly that 
$p(G_1+u_1u_2)+p(G_2+u_1u_2)$ and $p(G)$ have the same parity,
but $p(G_1+u_1u_2)+p(G_2+2u_1u_2)$ and $p(G)$ have different parity.
Thus 
\beeq
Q(G,\lambda)=
Q(G_1/u_1u_2,\lambda)Q(G_2,\lambda)
+Q(G_1+u_1u_2,\lambda)
\left [ 
Q(G_2+u_1u_2,\lambda)
+\frac{Q(G_2+2u_1u_2,\lambda)}{\lambda -1}\right ].
\eneq
By the given conditions and the assumption on $G$, 
$Q(G_1/u_1u_2,\lambda)$, $Q(G_1+u_1u_2,\lambda)$, 
and $Q(G_2,\lambda)$, $Q(G_2+u_1u_2,\lambda)$
and $Q(G_2+2u_1u_2,\lambda)$ 
are all positive.
Hence
%if the inequality $Q(G_2+u_1u_2,\lambda)\ge (\lambda-1)Q(G_2,\lambda)$ holds,
$Q(G,\lambda)>0$,
a contradiction. 

Hence Claim~\thecountclaim\ holds. 
By the above claims, we know that 
$G$ is non-separable and does not satisfy condition (iii), 
contradicting the the given conditions. 
Thus the result holds.
\proofend

By Lemma~\ref{general-interval-le}, the following 
result is immediately obtained.

\begin{theo}\quad
\relabel{general-interval}
Let $\sets$ be a family of bridgeless connected graphs and 
$\beta$ a real number in $(1,2]$.
Assume that there exists $\sets'\subseteq\sets$ 
such that condition (i) in Lemma~\ref{general-interval-le} 
holds for all $\lambda\in (1,\beta)$
and both conditions (ii) and (iii) in Lemma~\ref{general-interval-le} hold, then 
$Q(G,\lambda)>0$ for all graphs $G\in \sets$
and all real $\lambda\in (1,\beta)$.
\proofend
\end{theo}

\resection{How to determine $\xi_k$}\relabel{sec4}

Recall that $\Psi_k$ is the set of bridgeless connected graphs $G$ 
with $|W(G)|\le k$ and 
$\xi_k$ is the supremum in $(1,2]$ 
such that every graph in $\Psi_k$ has no flow roots in $(1,\xi_k)$.
In this section, we will show that 
$\xi_k$ can be determined by considering the set of graphs in   
$\Theta$ with exactly $k$ vertices, 
where $\Theta$ is the set of graphs defined by the two steps below:
\begin{enumerate}
\item[(i)] $Z_3\in \Theta$, where  $Z_j$ is the graph with two vertices 
and $j$ parallel edges joining these two vertices; and 
\item[(ii)] $G(e)\in \Theta$ for every $G\in \Theta$ and 
every $e\in E(G)$, 
where $G(e)$ is the graph obtained from $G-e$ 
by adding a new vertex $w$ and adding two parallel edges 
joining $w$ and $u_i$ for both $i=1,2$, as shown in 
Figure~\ref{f4}.
\end{enumerate}

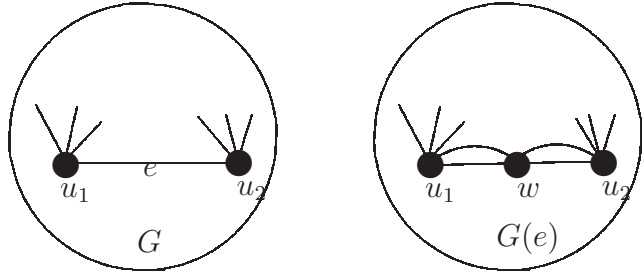
\begin{figure}[h!]
\centering 
%\scalebox{1.0}{\input{f4.pic}}
%\input f1.pic
%TeXCAD (http://texcad.sf.net/) Picture. File: [f4.pic]. Options on following lines.
%\grade{\on}
%\emlines{\off}
%\epic{\off}
%\beziermacro{\on}
%\reduce{\on}
%\snapping{\off}
%\pvinsert{% Your \input, \def, etc. here}
%\quality{8.000}
%\graddiff{0.005}
%\snapasp{1}
%\zoom{4.0000}
\unitlength 1mm % = 2.845pt
\linethickness{0.4pt}
\ifx\plotpoint\undefined\newsavebox{\plotpoint}\fi % GNUPLOT compatibility
\begin{picture}(89.005,38.755)(0,0)
\put(-130.5,-93.25){\circle*{3.354}}
\put(-82,-93.25){\circle*{3.354}}
\put(-81.75,-93.75){\circle*{3.5}}
\put(12.5,17.5){\line(1,0){22.75}}
\put(12.5,17.25){\circle*{3.536}}
\put(61,17.25){\circle*{3.536}}
\put(35.5,17.5){\circle*{3.536}}
\put(84,17.5){\circle*{3.536}}
\put(11.75,13){$u_1$}
\put(60.25,13){$u_1$}
\put(35.25,13.25){$u_2$}
\put(83.75,13.25){$u_2$}
\put(22,5.75){$G$}
%\circle(22.75,21){35.511}
\put(40.505,21){\line(0,1){.848}}
\put(40.485,21.848){\line(0,1){.846}}
\put(40.425,22.694){\line(0,1){.8422}}
\multiput(40.323,23.536)(-.028236,.167278){5}{\line(0,1){.167278}}
\multiput(40.182,24.373)(-.03016,.138116){6}{\line(0,1){.138116}}
\multiput(40.001,25.201)(-.031476,.117015){7}{\line(0,1){.117015}}
\multiput(39.781,26.02)(-.0324,.100956){8}{\line(0,1){.100956}}
\multiput(39.522,26.828)(-.033053,.088261){9}{\line(0,1){.088261}}
\multiput(39.224,27.622)(-.033508,.077923){10}{\line(0,1){.077923}}
\multiput(38.889,28.402)(-.030992,.063528){12}{\line(0,1){.063528}}
\multiput(38.517,29.164)(-.0313763,.0572084){13}{\line(0,1){.0572084}}
\multiput(38.109,29.908)(-.031639,.05167){14}{\line(0,1){.05167}}
\multiput(37.666,30.631)(-.0317992,.04676){15}{\line(0,1){.04676}}
\multiput(37.189,31.332)(-.0318713,.0423637){16}{\line(0,1){.0423637}}
\multiput(36.68,32.01)(-.0318665,.0383937){17}{\line(0,1){.0383937}}
\multiput(36.138,32.663)(-.0336638,.036828){17}{\line(0,1){.036828}}
\multiput(35.566,33.289)(-.0334185,.0332238){18}{\line(-1,0){.0334185}}
\multiput(34.964,33.887)(-.0370239,.0334481){17}{\line(-1,0){.0370239}}
\multiput(34.335,34.456)(-.0409903,.0336194){16}{\line(-1,0){.0409903}}
\multiput(33.679,34.994)(-.0453858,.0337316){15}{\line(-1,0){.0453858}}
\multiput(32.998,35.5)(-.046945,.0315255){15}{\line(-1,0){.046945}}
\multiput(32.294,35.972)(-.051854,.0313366){14}{\line(-1,0){.051854}}
\multiput(31.568,36.411)(-.062173,.033628){12}{\line(-1,0){.062173}}
\multiput(30.822,36.815)(-.0695,.033404){11}{\line(-1,0){.0695}}
\multiput(30.057,37.182)(-.078118,.033052){10}{\line(-1,0){.078118}}
\multiput(29.276,37.513)(-.088452,.032537){9}{\line(-1,0){.088452}}
\multiput(28.48,37.806)(-.101143,.03181){8}{\line(-1,0){.101143}}
\multiput(27.671,38.06)(-.117197,.030792){7}{\line(-1,0){.117197}}
\multiput(26.85,38.276)(-.138289,.029353){6}{\line(-1,0){.138289}}
\multiput(26.021,38.452)(-.16744,.027258){5}{\line(-1,0){.16744}}
\put(25.184,38.588){\line(-1,0){.8428}}
\put(24.341,38.684){\line(-1,0){.8464}}
\put(23.494,38.74){\line(-1,0){1.6959}}
\put(21.798,38.73){\line(-1,0){.8457}}
\put(20.953,38.664){\line(-1,0){.8416}}
\multiput(20.111,38.558)(-.16711,-.029212){5}{\line(-1,0){.16711}}
\multiput(19.276,38.412)(-.137937,-.030967){6}{\line(-1,0){.137937}}
\multiput(18.448,38.226)(-.116829,-.032159){7}{\line(-1,0){.116829}}
\multiput(17.63,38.001)(-.100765,-.032989){8}{\line(-1,0){.100765}}
\multiput(16.824,37.737)(-.088066,-.033568){9}{\line(-1,0){.088066}}
\multiput(16.032,37.435)(-.07066,-.030875){11}{\line(-1,0){.07066}}
\multiput(15.254,37.096)(-.063346,-.031363){12}{\line(-1,0){.063346}}
\multiput(14.494,36.719)(-.0570242,-.0317099){13}{\line(-1,0){.0570242}}
\multiput(13.753,36.307)(-.0514844,-.0319402){14}{\line(-1,0){.0514844}}
\multiput(13.032,35.86)(-.0465735,-.0320717){15}{\line(-1,0){.0465735}}
\multiput(12.333,35.379)(-.0421769,-.0321182){16}{\line(-1,0){.0421769}}
\multiput(11.659,34.865)(-.0382069,-.0320902){17}{\line(-1,0){.0382069}}
\multiput(11.009,34.319)(-.0345957,-.0319962){18}{\line(-1,0){.0345957}}
\multiput(10.386,33.744)(-.0330281,-.0336119){18}{\line(0,-1){.0336119}}
\multiput(9.792,33.139)(-.0332313,-.0372187){17}{\line(0,-1){.0372187}}
\multiput(9.227,32.506)(-.0333794,-.041186){16}{\line(0,-1){.041186}}
\multiput(8.693,31.847)(-.0334659,-.045582){15}{\line(0,-1){.045582}}
\multiput(8.191,31.163)(-.033483,-.0504946){14}{\line(0,-1){.0504946}}
\multiput(7.722,30.456)(-.0334204,-.0560389){13}{\line(0,-1){.0560389}}
\multiput(7.288,29.728)(-.033265,-.062369){12}{\line(0,-1){.062369}}
\multiput(6.888,28.979)(-.032998,-.069694){11}{\line(0,-1){.069694}}
\multiput(6.525,28.213)(-.032595,-.078309){10}{\line(0,-1){.078309}}
\multiput(6.199,27.429)(-.03202,-.088641){9}{\line(0,-1){.088641}}
\multiput(5.911,26.632)(-.031219,-.101328){8}{\line(0,-1){.101328}}
\multiput(5.662,25.821)(-.030107,-.117375){7}{\line(0,-1){.117375}}
\multiput(5.451,24.999)(-.028545,-.138458){6}{\line(0,-1){.138458}}
\multiput(5.28,24.169)(-.03285,-.2095){4}{\line(0,-1){.2095}}
\put(5.148,23.331){\line(0,-1){.8433}}
\put(5.057,22.487){\line(0,-1){.8467}}
\put(5.006,21.641){\line(0,-1){.8482}}
\put(4.996,20.793){\line(0,-1){.8477}}
\put(5.026,19.945){\line(0,-1){.8453}}
\put(5.096,19.1){\line(0,-1){.8409}}
\multiput(5.207,18.259)(.030188,-.166937){5}{\line(0,-1){.166937}}
\multiput(5.358,17.424)(.031772,-.137754){6}{\line(0,-1){.137754}}
\multiput(5.549,16.597)(.032841,-.116639){7}{\line(0,-1){.116639}}
\multiput(5.779,15.781)(.033577,-.100571){8}{\line(0,-1){.100571}}
\multiput(6.047,14.976)(.030674,-.079082){10}{\line(0,-1){.079082}}
\multiput(6.354,14.186)(.031287,-.070479){11}{\line(0,-1){.070479}}
\multiput(6.698,13.41)(.031732,-.063162){12}{\line(0,-1){.063162}}
\multiput(7.079,12.652)(.0320424,-.056838){13}{\line(0,-1){.056838}}
\multiput(7.496,11.914)(.0322404,-.0512969){14}{\line(0,-1){.0512969}}
\multiput(7.947,11.195)(.0323432,-.0463854){15}{\line(0,-1){.0463854}}
\multiput(8.432,10.5)(.032364,-.0419886){16}{\line(0,-1){.0419886}}
\multiput(8.95,9.828)(.0323129,-.0380188){17}{\line(0,-1){.0380188}}
\multiput(9.499,9.181)(.0321977,-.0344082){18}{\line(0,-1){.0344082}}
\multiput(10.079,8.562)(.0338043,-.0328312){18}{\line(1,0){.0338043}}
\multiput(10.687,7.971)(.0374121,-.0330134){17}{\line(1,0){.0374121}}
\multiput(11.323,7.41)(.0413802,-.0331383){16}{\line(1,0){.0413802}}
\multiput(11.985,6.88)(.0457767,-.0331991){15}{\line(1,0){.0457767}}
\multiput(12.672,6.382)(.0506893,-.0331875){14}{\line(1,0){.0506893}}
\multiput(13.382,5.917)(.0562331,-.0330925){13}{\line(1,0){.0562331}}
\multiput(14.113,5.487)(.062562,-.0329){12}{\line(1,0){.062562}}
\multiput(14.864,5.092)(.069886,-.03259){11}{\line(1,0){.069886}}
\multiput(15.632,4.734)(.078498,-.032137){10}{\line(1,0){.078498}}
\multiput(16.417,4.412)(.088826,-.031502){9}{\line(1,0){.088826}}
\multiput(17.217,4.129)(.101508,-.030626){8}{\line(1,0){.101508}}
\multiput(18.029,3.884)(.117549,-.029421){7}{\line(1,0){.117549}}
\multiput(18.852,3.678)(.166347,-.033283){5}{\line(1,0){.166347}}
\multiput(19.683,3.511)(.20968,-.03163){4}{\line(1,0){.20968}}
\put(20.522,3.385){\line(1,0){.8438}}
\put(21.366,3.299){\line(1,0){.847}}
\put(22.213,3.253){\line(1,0){.8482}}
\put(23.061,3.247){\line(1,0){.8475}}
\put(23.909,3.282){\line(1,0){.8449}}
\put(24.753,3.358){\line(1,0){.8403}}
\multiput(25.594,3.474)(.166758,.031162){5}{\line(1,0){.166758}}
\multiput(26.428,3.63)(.137566,.032576){6}{\line(1,0){.137566}}
\multiput(27.253,3.825)(.116446,.033522){7}{\line(1,0){.116446}}
\multiput(28.068,4.06)(.08922,.030368){9}{\line(1,0){.08922}}
\multiput(28.871,4.333)(.078901,.031135){10}{\line(1,0){.078901}}
\multiput(29.66,4.644)(.070295,.031698){11}{\line(1,0){.070295}}
\multiput(30.433,4.993)(.062976,.032101){12}{\line(1,0){.062976}}
\multiput(31.189,5.378)(.0566499,.0323739){13}{\line(1,0){.0566499}}
\multiput(31.925,5.799)(.0511078,.0325394){14}{\line(1,0){.0511078}}
\multiput(32.641,6.255)(.0461957,.0326136){15}{\line(1,0){.0461957}}
\multiput(33.334,6.744)(.0417988,.0326087){16}{\line(1,0){.0417988}}
\multiput(34.003,7.266)(.0378294,.0325344){17}{\line(1,0){.0378294}}
\multiput(34.646,7.819)(.0342196,.0323981){18}{\line(1,0){.0342196}}
\multiput(35.262,8.402)(.0326332,.0339955){18}{\line(0,1){.0339955}}
\multiput(35.849,9.014)(.0327943,.0376043){17}{\line(0,1){.0376043}}
\multiput(36.407,9.653)(.032896,.0415731){16}{\line(0,1){.0415731}}
\multiput(36.933,10.318)(.0329312,.0459698){15}{\line(0,1){.0459698}}
\multiput(37.427,11.008)(.0328909,.0508823){14}{\line(0,1){.0508823}}
\multiput(37.887,11.72)(.0327635,.0564254){13}{\line(0,1){.0564254}}
\multiput(38.313,12.454)(.032534,.062753){12}{\line(0,1){.062753}}
\multiput(38.704,13.207)(.032182,.070075){11}{\line(0,1){.070075}}
\multiput(39.058,13.977)(.031678,.078685){10}{\line(0,1){.078685}}
\multiput(39.374,14.764)(.030982,.089009){9}{\line(0,1){.089009}}
\multiput(39.653,15.565)(.030033,.101685){8}{\line(0,1){.101685}}
\multiput(39.894,16.379)(.033523,.137338){6}{\line(0,1){.137338}}
\multiput(40.095,17.203)(.032311,.166539){5}{\line(0,1){.166539}}
\multiput(40.256,18.036)(.0304,.20987){4}{\line(0,1){.20987}}
\put(40.378,18.875){\line(0,1){.8443}}
\put(40.459,19.719){\line(0,1){1.2807}}
%\end
%\circle(71.25,21){35.511}
\put(89.005,21){\line(0,1){.848}}
\put(88.985,21.848){\line(0,1){.846}}
\put(88.925,22.694){\line(0,1){.8422}}
\multiput(88.823,23.536)(-.028236,.167278){5}{\line(0,1){.167278}}
\multiput(88.682,24.373)(-.03016,.138116){6}{\line(0,1){.138116}}
\multiput(88.501,25.201)(-.031476,.117015){7}{\line(0,1){.117015}}
\multiput(88.281,26.02)(-.0324,.100956){8}{\line(0,1){.100956}}
\multiput(88.022,26.828)(-.033053,.088261){9}{\line(0,1){.088261}}
\multiput(87.724,27.622)(-.033508,.077923){10}{\line(0,1){.077923}}
\multiput(87.389,28.402)(-.030992,.063528){12}{\line(0,1){.063528}}
\multiput(87.017,29.164)(-.0313763,.0572084){13}{\line(0,1){.0572084}}
\multiput(86.609,29.908)(-.031639,.05167){14}{\line(0,1){.05167}}
\multiput(86.166,30.631)(-.0317992,.04676){15}{\line(0,1){.04676}}
\multiput(85.689,31.332)(-.0318713,.0423637){16}{\line(0,1){.0423637}}
\multiput(85.18,32.01)(-.0318665,.0383937){17}{\line(0,1){.0383937}}
\multiput(84.638,32.663)(-.0336638,.036828){17}{\line(0,1){.036828}}
\multiput(84.066,33.289)(-.0334185,.0332238){18}{\line(-1,0){.0334185}}
\multiput(83.464,33.887)(-.0370239,.0334481){17}{\line(-1,0){.0370239}}
\multiput(82.835,34.456)(-.0409903,.0336194){16}{\line(-1,0){.0409903}}
\multiput(82.179,34.994)(-.0453858,.0337316){15}{\line(-1,0){.0453858}}
\multiput(81.498,35.5)(-.046945,.0315255){15}{\line(-1,0){.046945}}
\multiput(80.794,35.972)(-.051854,.0313366){14}{\line(-1,0){.051854}}
\multiput(80.068,36.411)(-.062173,.033628){12}{\line(-1,0){.062173}}
\multiput(79.322,36.815)(-.0695,.033404){11}{\line(-1,0){.0695}}
\multiput(78.557,37.182)(-.078118,.033052){10}{\line(-1,0){.078118}}
\multiput(77.776,37.513)(-.088452,.032537){9}{\line(-1,0){.088452}}
\multiput(76.98,37.806)(-.101143,.03181){8}{\line(-1,0){.101143}}
\multiput(76.171,38.06)(-.117197,.030792){7}{\line(-1,0){.117197}}
\multiput(75.35,38.276)(-.138289,.029353){6}{\line(-1,0){.138289}}
\multiput(74.521,38.452)(-.16744,.027258){5}{\line(-1,0){.16744}}
\put(73.684,38.588){\line(-1,0){.8428}}
\put(72.841,38.684){\line(-1,0){.8464}}
\put(71.994,38.74){\line(-1,0){1.6959}}
\put(70.298,38.73){\line(-1,0){.8457}}
\put(69.453,38.664){\line(-1,0){.8416}}
\multiput(68.611,38.558)(-.16711,-.029212){5}{\line(-1,0){.16711}}
\multiput(67.776,38.412)(-.137937,-.030967){6}{\line(-1,0){.137937}}
\multiput(66.948,38.226)(-.116829,-.032159){7}{\line(-1,0){.116829}}
\multiput(66.13,38.001)(-.100765,-.032989){8}{\line(-1,0){.100765}}
\multiput(65.324,37.737)(-.088066,-.033568){9}{\line(-1,0){.088066}}
\multiput(64.532,37.435)(-.07066,-.030875){11}{\line(-1,0){.07066}}
\multiput(63.754,37.096)(-.063346,-.031363){12}{\line(-1,0){.063346}}
\multiput(62.994,36.719)(-.0570242,-.0317099){13}{\line(-1,0){.0570242}}
\multiput(62.253,36.307)(-.0514844,-.0319402){14}{\line(-1,0){.0514844}}
\multiput(61.532,35.86)(-.0465735,-.0320717){15}{\line(-1,0){.0465735}}
\multiput(60.833,35.379)(-.0421769,-.0321182){16}{\line(-1,0){.0421769}}
\multiput(60.159,34.865)(-.0382069,-.0320902){17}{\line(-1,0){.0382069}}
\multiput(59.509,34.319)(-.0345957,-.0319962){18}{\line(-1,0){.0345957}}
\multiput(58.886,33.744)(-.0330281,-.0336119){18}{\line(0,-1){.0336119}}
\multiput(58.292,33.139)(-.0332313,-.0372187){17}{\line(0,-1){.0372187}}
\multiput(57.727,32.506)(-.0333794,-.041186){16}{\line(0,-1){.041186}}
\multiput(57.193,31.847)(-.0334659,-.045582){15}{\line(0,-1){.045582}}
\multiput(56.691,31.163)(-.033483,-.0504946){14}{\line(0,-1){.0504946}}
\multiput(56.222,30.456)(-.0334204,-.0560389){13}{\line(0,-1){.0560389}}
\multiput(55.788,29.728)(-.033265,-.062369){12}{\line(0,-1){.062369}}
\multiput(55.388,28.979)(-.032998,-.069694){11}{\line(0,-1){.069694}}
\multiput(55.025,28.213)(-.032595,-.078309){10}{\line(0,-1){.078309}}
\multiput(54.699,27.429)(-.03202,-.088641){9}{\line(0,-1){.088641}}
\multiput(54.411,26.632)(-.031219,-.101328){8}{\line(0,-1){.101328}}
\multiput(54.162,25.821)(-.030107,-.117375){7}{\line(0,-1){.117375}}
\multiput(53.951,24.999)(-.028545,-.138458){6}{\line(0,-1){.138458}}
\multiput(53.78,24.169)(-.03285,-.2095){4}{\line(0,-1){.2095}}
\put(53.648,23.331){\line(0,-1){.8433}}
\put(53.557,22.487){\line(0,-1){.8467}}
\put(53.506,21.641){\line(0,-1){.8482}}
\put(53.496,20.793){\line(0,-1){.8477}}
\put(53.526,19.945){\line(0,-1){.8453}}
\put(53.596,19.1){\line(0,-1){.8409}}
\multiput(53.707,18.259)(.030188,-.166937){5}{\line(0,-1){.166937}}
\multiput(53.858,17.424)(.031772,-.137754){6}{\line(0,-1){.137754}}
\multiput(54.049,16.597)(.032841,-.116639){7}{\line(0,-1){.116639}}
\multiput(54.279,15.781)(.033577,-.100571){8}{\line(0,-1){.100571}}
\multiput(54.547,14.976)(.030674,-.079082){10}{\line(0,-1){.079082}}
\multiput(54.854,14.186)(.031287,-.070479){11}{\line(0,-1){.070479}}
\multiput(55.198,13.41)(.031732,-.063162){12}{\line(0,-1){.063162}}
\multiput(55.579,12.652)(.0320424,-.056838){13}{\line(0,-1){.056838}}
\multiput(55.996,11.914)(.0322404,-.0512969){14}{\line(0,-1){.0512969}}
\multiput(56.447,11.195)(.0323432,-.0463854){15}{\line(0,-1){.0463854}}
\multiput(56.932,10.5)(.032364,-.0419886){16}{\line(0,-1){.0419886}}
\multiput(57.45,9.828)(.0323129,-.0380188){17}{\line(0,-1){.0380188}}
\multiput(57.999,9.181)(.0321977,-.0344082){18}{\line(0,-1){.0344082}}
\multiput(58.579,8.562)(.0338043,-.0328312){18}{\line(1,0){.0338043}}
\multiput(59.187,7.971)(.0374121,-.0330134){17}{\line(1,0){.0374121}}
\multiput(59.823,7.41)(.0413802,-.0331383){16}{\line(1,0){.0413802}}
\multiput(60.485,6.88)(.0457767,-.0331991){15}{\line(1,0){.0457767}}
\multiput(61.172,6.382)(.0506893,-.0331875){14}{\line(1,0){.0506893}}
\multiput(61.882,5.917)(.0562331,-.0330925){13}{\line(1,0){.0562331}}
\multiput(62.613,5.487)(.062562,-.0329){12}{\line(1,0){.062562}}
\multiput(63.364,5.092)(.069886,-.03259){11}{\line(1,0){.069886}}
\multiput(64.132,4.734)(.078498,-.032137){10}{\line(1,0){.078498}}
\multiput(64.917,4.412)(.088826,-.031502){9}{\line(1,0){.088826}}
\multiput(65.717,4.129)(.101508,-.030626){8}{\line(1,0){.101508}}
\multiput(66.529,3.884)(.117549,-.029421){7}{\line(1,0){.117549}}
\multiput(67.352,3.678)(.166347,-.033283){5}{\line(1,0){.166347}}
\multiput(68.183,3.511)(.20968,-.03163){4}{\line(1,0){.20968}}
\put(69.022,3.385){\line(1,0){.8438}}
\put(69.866,3.299){\line(1,0){.847}}
\put(70.713,3.253){\line(1,0){.8482}}
\put(71.561,3.247){\line(1,0){.8475}}
\put(72.409,3.282){\line(1,0){.8449}}
\put(73.253,3.358){\line(1,0){.8403}}
\multiput(74.094,3.474)(.166758,.031162){5}{\line(1,0){.166758}}
\multiput(74.928,3.63)(.137566,.032576){6}{\line(1,0){.137566}}
\multiput(75.753,3.825)(.116446,.033522){7}{\line(1,0){.116446}}
\multiput(76.568,4.06)(.08922,.030368){9}{\line(1,0){.08922}}
\multiput(77.371,4.333)(.078901,.031135){10}{\line(1,0){.078901}}
\multiput(78.16,4.644)(.070295,.031698){11}{\line(1,0){.070295}}
\multiput(78.933,4.993)(.062976,.032101){12}{\line(1,0){.062976}}
\multiput(79.689,5.378)(.0566499,.0323739){13}{\line(1,0){.0566499}}
\multiput(80.425,5.799)(.0511078,.0325394){14}{\line(1,0){.0511078}}
\multiput(81.141,6.255)(.0461957,.0326136){15}{\line(1,0){.0461957}}
\multiput(81.834,6.744)(.0417988,.0326087){16}{\line(1,0){.0417988}}
\multiput(82.503,7.266)(.0378294,.0325344){17}{\line(1,0){.0378294}}
\multiput(83.146,7.819)(.0342196,.0323981){18}{\line(1,0){.0342196}}
\multiput(83.762,8.402)(.0326332,.0339955){18}{\line(0,1){.0339955}}
\multiput(84.349,9.014)(.0327943,.0376043){17}{\line(0,1){.0376043}}
\multiput(84.907,9.653)(.032896,.0415731){16}{\line(0,1){.0415731}}
\multiput(85.433,10.318)(.0329312,.0459698){15}{\line(0,1){.0459698}}
\multiput(85.927,11.008)(.0328909,.0508823){14}{\line(0,1){.0508823}}
\multiput(86.387,11.72)(.0327635,.0564254){13}{\line(0,1){.0564254}}
\multiput(86.813,12.454)(.032534,.062753){12}{\line(0,1){.062753}}
\multiput(87.204,13.207)(.032182,.070075){11}{\line(0,1){.070075}}
\multiput(87.558,13.977)(.031678,.078685){10}{\line(0,1){.078685}}
\multiput(87.874,14.764)(.030982,.089009){9}{\line(0,1){.089009}}
\multiput(88.153,15.565)(.030033,.101685){8}{\line(0,1){.101685}}
\multiput(88.394,16.379)(.033523,.137338){6}{\line(0,1){.137338}}
\multiput(88.595,17.203)(.032311,.166539){5}{\line(0,1){.166539}}
\multiput(88.756,18.036)(.0304,.20987){4}{\line(0,1){.20987}}
\put(88.878,18.875){\line(0,1){.8443}}
\put(88.959,19.719){\line(0,1){1.2807}}
%\end
\put(69.75,6.5){$G(e)$}
\put(22.75,15.75){$e$}
%\emline(12.5,17.75)(8.5,25.25)
\multiput(12.5,17.75)(-.033613445,.06302521){119}{\line(0,1){.06302521}}
%\end
%\emline(60.75,17.75)(56.75,25.25)
\multiput(60.75,17.75)(-.033613445,.06302521){119}{\line(0,1){.06302521}}
%\end
%\emline(12.25,17.25)(14,25)
\multiput(12.25,17.25)(.03365385,.14903846){52}{\line(0,1){.14903846}}
%\end
%\emline(60.5,17.25)(62.25,25)
\multiput(60.5,17.25)(.03365385,.14903846){52}{\line(0,1){.14903846}}
%\end
%\emline(12.25,17.75)(17.25,23)
\multiput(12.25,17.75)(.033557047,.035234899){149}{\line(0,1){.035234899}}
%\end
%\emline(60.5,17.75)(65.5,23)
\multiput(60.5,17.75)(.033557047,.035234899){149}{\line(0,1){.035234899}}
%\end
%\emline(35.5,17.25)(30,23.75)
\multiput(35.5,17.25)(-.033536585,.039634146){164}{\line(0,1){.039634146}}
%\end
%\emline(35.25,17.5)(37.25,24)
\multiput(35.25,17.5)(.03333333,.10833333){60}{\line(0,1){.10833333}}
%\end
%\emline(83.5,17.5)(85.5,24)
\multiput(83.5,17.5)(.03333333,.10833333){60}{\line(0,1){.10833333}}
%\end
%\emline(35.5,17.25)(33.75,24)
\multiput(35.5,17.25)(-.03365385,.12980769){52}{\line(0,1){.12980769}}
%\end
%\emline(83.75,17.25)(82,24)
\multiput(83.75,17.25)(-.03365385,.12980769){52}{\line(0,1){.12980769}}
%\end
%\emline(83.75,17.75)(80.25,23.5)
\multiput(83.75,17.75)(-.033653846,.055288462){104}{\line(0,1){.055288462}}
%\end
%\emline(60.75,17.25)(83.5,17.75)
\multiput(60.75,17.25)(1.5166667,.0333333){15}{\line(1,0){1.5166667}}
%\end
\put(72.5,17.25){\circle*{3.64}}
\qbezier(60.5,17.5)(66.625,21.75)(72.25,18)
\qbezier(72.25,18)(77.625,22)(83.5,18)
\put(72.5,13){$w$}
\end{picture}

\caption{\relabel{f4} Graphs $G$ and $G(e)$}
\end{figure}
As examples, 
we also determine the values of $\xi_k$ for  $k\le 5$:
$\xi_k=2$ for $k=0,1,2$, % and %find that 
$\xi_3=1.430159709\cdots$, 
$\xi_4=1.361103081\cdots$ and $\xi_5=1.317672196\cdots$,
where the last three numbers in $(1, 2)$ 
%$1.430159709\cdots$, $1.361103081\cdots$ and $1.317672196\cdots$
are the real zeros of $\lambda^3-5\lambda^2+10\lambda-7$,
$\lambda^3-4\lambda^2+8\lambda-6$ and 
$\lambda^3 -6\lambda^2+13\lambda-9$ 
in $(1,2)$ respectively.

For any bridgeless graph $G$, 
let $\eta(G)$ be the minimum flow root of $G$ in 
the interval $(1,2]$ if such root exists 
and $\eta(G)=2$ otherwise. 
By Theorem~\ref{Wakelin}, we have $32/27< \eta(G)\le 2$ for 
every bridgeless graph $G$. 
For any set $\sets$ of bridgeless graphs,  let
\begin{equation}
\eta(\sets)=
\left \{ 
\begin{array}{ll}
\inf \{\eta(G): G\in \sets\}, \quad &\mbox{if } 
\sets\ne \emptyset;\\
2, &\mbox{otherwise}.
\end{array}
\right.
\end{equation}
Thus $\xi_k=\eta(\Psi_k)$ and 
$\xi_0, \xi_1,\xi_2, \cdots$ is a non-increasing sequence. 

Let $\Phi$ be the set of non-separable graphs $G$
with $|V(G)|\ge 2$ such that 
the following conditions are all satisfied:
\begin{enumerate} \label{con-phi}
\item[(a')]  $G-e$ is non-separable for each edge $e$ in $G$;

\item[(b')] 
$b(G/e)$ is even for each edge $e$ in $G$; and

\item[(c')]  
if $G_1$ and $G_2$ are subgraphs of $G$  
such that $|E(G_i)|\ge 2$ for $i=1,2$, 
$V(G_1)\cap V(G_2)=\{u_1,u_2\}$,
$V(G_1)\cup V(G_2)=V(G)$,
$E(G_1)\cap E(G_2)=\emptyset$
and $E(G_1)\cup E(G_2)=E(G)$, as shown in Figure~\ref{f3}(a),
then the three integers $b(G_1/u_1u_2), b(G_1)-1$ and  $b(G_2)$ 
all have the same parity.   
\end{enumerate}

Instead we prove directly that $\xi_k$ can be determined by
considering the set of graphs in $\Theta$ with exactly $k$ vertices,
we will obtain this conclusion by proving that 
$\Theta$ is actually equal to the set $\Phi$
and $\xi_k=\eta(\Phi_k)$, 
where $\Phi_k$ is the set of graphs $G\in \Phi$ with $|V(G)|=k$.

We will first show that $\xi_k=\min \{\eta(\Phi_i): 2\le i\le k\}$
and the following result will be applied in proving it. 
For a graph $G=(V,E)$ and $x\in V$,  
let $N(x)=\{u: xu\in E(G)\}$. 
So $d(x)\ge |N(x)|$, where equality holds if and only if 
$G$ has no loops or parallel edges incident with $x$.

\begin{lem}\quad
\relabel{v-degree3}
Let $G=(V,E)$ be a non-separable graph
with $|V|\ge 3$ 
and $x\in V$ with $d(x)\le 3$.
If $G-e$ is non-separable for every edge $e$ incident with $x$,
then $G/e'$ is also non-separable for every edge $e'$ 
incident with $x$.
\end{lem}

\proof Suppose that $G/e'$ is separable for some edge $e'$ 
incident with $x$.

Suppose that $|N(x)|\le 2$.
Since $|V|\ge 3$ and $G$ is non-separable, $|N(x)|=2$.
As $d(x)\le 3$ and $|N(x)|=2$, there is a single edge 
incident with $x$, and observe that 
$G-e$ is separable for such an edge $e$, a contradiction. 
Thus $|N(x)|=3$, implying that $d(x)=3$ and no parallel edges are
incident with $x$. 

Since $G/e'$ is separable and $d(x)=|N(x)|=3$, 
$G-e$ must be separable 
for every edge $e$ which is different from $e'$ 
and is incident with $x$, a contradiction. 
\proofend

%Recall that $\Phi_k$ is the set of graphs $G$ in $\Phi$ with $|V(G)|=k$.

\begin{lem}
\relabel{determine xik}
For $k\ge 2$, 
$\xi_k=\min\{\eta(\Phi_i): i=2,3,\cdots,k\}$.
\end{lem}

\proof
We prove this result by applying Theorem~\ref{general-interval}.
Let $\sets=\Psi_k$ and 
\beeq
\sets'=\{L, Z_2\}\cup  \bigcup_{2\le i\le k}\Phi_i, 
\eneq
where $L$ is the graph with one vertex and one loop.
Let  $\beta=\min\{\eta(\Phi_i): i=2,3,\cdots,k\}$. 

By the definition on $\beta$, 
we have 
$Q(G,\lambda)>0$ for all $G\in \sets'$ and  all 
$\lambda\in (1,\beta)$.
Thus condition (i) of Lemma~\ref{general-interval-le} is 
satisfied for all $\lambda\in (1,\beta)$.

Observe that for any $G\in \sets$ ($=\Psi_k$), 
if $G$ is separable,
then $|W(B)|\le |W(G)|\le k$ for each block $B$ of $G$ 
and so each block of $G$ belongs to $\sets$. 
Hence  Condition (ii) of Lemma~\ref{general-interval-le} is also
satisfied.

If Condition (iii) of Lemma~\ref{general-interval-le} 
holds for every  non-separable graph 
$G\in \sets\setminus \sets'$, then this result holds by 
Theorem~\ref{general-interval}.
Now suppose that Condition (iii) of Lemma~\ref{general-interval-le}
does not hold for some non-separable graph 
$G\in \sets\setminus \sets'$. 
So none of conditions (a), (b), (c) and (d) of (iii) 
in Lemma~\ref{general-interval-le} is satisfied for $G$.
We shall show that $W(G)=V(G)$ and $G$ satisfies 
conditions (a'), (b') and (c') in page~\pageref{con-phi},
and thus $G\in \Phi$, implying that  $G\in \sets'$,
a contradiction.  

If $G$ does not satisfy condition (a'), 
then for some edge $e$ in $G$, 
$G-e$ has a cut-vertex $u$ for some edge $e$.
Then $W(G_i)\subseteq W(G)$ and so $G_i\in \Psi_k$
for $i=1,2$, where $G_1$ and $G_2$ 
are the two graphs stated in Lemma~\ref{v-edge}.
Thus condition (a) is satisfied, a contradiction. 
Hence $G$ satisfies condition (a').

Before we can show that $G$ satisfies conditions (b') and (c'),
we need to show that $W(G)=V(G)$. 
Suppose that $W(G)\ne V(G)$. 
Let $x\in V(G)\setminus W(G)$ and $u\in N(x)$.  
If $W(G)\ne \emptyset$, 
$x$ and $u$ are selected so that $u\in W(G)$.
It is clear that $|V(G)|\ge 3$; 
otherwise, $d(x)\le 3$ implies that $G=Z_2$ or $Z_3$ 
and so $G\in \sets'$, a contradiction.
As $G$ satisfies condition (a'), 
$G-e$ is non-separable for 
every edge $e$ incident with $x$,  
and so Lemma~\ref{v-degree3} implies that 
$G/xu$ is non-separable.
As $W(G-xu)\subseteq W(G)$, $G-xu\in \Psi_k$. 
If $W(G)=\emptyset$, 
then $G/xu\in \Psi_1\subseteq \Psi_k$;
if $W(G)\ne \emptyset$, 
then $u\in W(G)$ and so $|W(G/xu)|=|W(G)|$,
implying that $G/xu\in \Psi_k$.
Thus condition (b) of (iii) in 
Lemma~\ref{general-interval-le} is satisfied, 
a contradiction. 
Hence $W(G)=V(G)$. 

Since $W(G)=V(G)$, we have $|V(G)|\le k$
and thus any bridgeless connected minor of $G$ 
belongs to $\Psi_k$. 
Since $G$ does not satisfy condition (b) of (iii) 
in Lemma~\ref{general-interval-le},
it immediately follows that $G$ satisfies 
condition (b') in page~\pageref{con-phi}.

We now show that $G$ also satisfies 
condition (c') in page~\pageref{con-phi}.
Suppose that this is not true. 
Then $G$ has  subgraphs $G_1$ and $G_2$  
such that $|E(G_i)|\ge 2$ for $i=1,2$, 
$V(G_1)\cap V(G_2)=\{u_1,u_2\}$,
$V(G_1)\cup V(G_2)=V(G)$,
$E(G_1)\cap E(G_2)=\emptyset$
$E(G_1)\cup E(G_2)=E(G)$, as shown in Figure~\ref{f3}(a),
%then $b(G_1)+b(G_2)$ is odd, 
but the three integers $b(G_1/u_1u_2), b(G_1)-1$ and  $b(G_2)$ 
don't have the same parity. 
Note that both $G_i$ and $G_i+u_1u_2$ belong to $\Psi_k$ for $i=1,2$.
Since $G$ does not satisfies condition (c) of (iii) 
in Lemma~\ref{general-interval-le}, 
$b(G_1)+b(G_2)$ is an odd number, i.e., 
$b(G_1)-1$ and $b(G_2)$ have the same parity.
Thus $b(G_2)$ and $b(G_1/u_1u_2)$ don't have the same parity, i.e., $b(G_2)+b(G_1/u_1u_2)$ is odd.
Note that 
$G_1+u_1u_2$, $G_1/u_1u_2$, $G_2$, 
$G_2+u_1u_2$ and 
$G_2+2u_1u_2$ all belong to $\sets$.
Since 
$G$ does not satisfies condition (c) of (iii), 
we have $|E(G_1)|<3$. 
Thus $|E(G_1)|=2$.
Since deleting any edge from $G$ does not produce a separable 
graph, the only two edges in $G_1$ are parallel edges 
joining $u_1$ and $u_2$. 
Thus $b(G_1/u_1u_2)=2$ and $b(G_1)=1$,
implying that $b(G_1/u_1u_2)$ and $b(G_1)-1$ have the same parity
and hence $b(G_1/u_1u_2), b(G_1)-1$ and  $b(G_2)$ 
all have the same parity, a contradiction. 

Hence $G$ satisfies condition (c'). 
Then, by definition of $\Phi$, $G\in \Phi$, 
implying that $G\in \Phi_i$, where $i=|V(G)|\le k$.
Thus $G\in \sets'$, contradicting the assumption on $G$.
\proofend

Later we will show that $\eta(\Phi_0),\eta(\Phi_1),\eta(\Phi_2), \eta(\Phi_3),\cdots$ 
is a non-increasing sequence and so Lemma~\ref{determine xik} 
implies that $\xi_k=\eta(\Phi_k)$ for $k\ge 2$.
 
Now we are going to show that $\Theta$ and $\Phi$ are actually the same set.
To prove this result, we need to apply some properties on graphs 
in $\Theta$ and $\Phi$.

\begin{lem}\quad
\relabel{graphs in phi}
Let $G=(V,E)\in\Phi$. 
Then for any distinct vertices $u_1,u_2$ in $G$,
$b(G/u_1u_2)\in \{1,3\}$.
\end{lem}

\proof The result is true when $|V|=2$. 
Assume that $|V|\ge 3$ and $b(G/u_1u_2)\ge 2$.
So there are  subgraphs $G_1$ and $G_2$ such that  
$V(G_1)\cap V(G_2)=\{u_1,u_2\}$,
$V(G_1)\cup V(G_2)=V(G)$,
$E(G_1)\cap E(G_2)=\emptyset$
and $E(G_1)\cup E(G_2)=E(G)$, as shown in Figure~\ref{f3}(a).
If $b(G/u_1u_2)\ge 4$, then $G_1$ and $G_2$ can be non-separable,
then $|E(G_i)|\ge 2$ and $b(G_1)=b(G_2)=1$, 
contradicting condition (c') that 
$b(G_1)-1$ and  $b(G_2)$ have the same parity. 
Now assume that $b(G/u_1u_2)=2$. 
So $G_i/u_1u_2$ is non-separable for $i=1,2$.
By condition (b'), we have $|E(G_i)|\ge 2$ for $i=1,2$.
Thus, by condition (c'), 
$b(G_1)+b(G_2)$ is an odd number at least 3.
Then $b(G_1)$ or $b(G_2)$ is even. 
Assume that $b(G_2)$ is even. Then $b(G_1/u_1u_2)$ 
must be even by condition (c'),
contradicting the fact that $G_1/u_1u_2$ is non-separable.
\proofend

\begin{lem}\quad
\relabel{new graphs in phi}
Let $G=(V,E)\in\Phi$.
Assume that $G_1$ and $G_2$ are proper subgraphs of $G$ such that  
$V(G_1)\cap V(G_2)=\{u_1,u_2\}$,
$V(G_1)\cup V(G_2)=V(G)$,
$E(G_1)\cap E(G_2)=\emptyset$
and $E(G_1)\cup E(G_2)=E(G)$, as shown in Figure~\ref{f3}(a).
If $|E(G_2)|\ge 2$ and $G_2$ is non-separable, then 
$G_2+u_1u_2\in \Phi$.
\end{lem}

\proof Assume that $|E(G_2)|\ge 2$ and $G_2$ is non-separable.
If $|E(G_1)|=1$, then $G_1+u_1u_2$ is $G$ and so the result holds.

Now assume that $|E(G_1)|\ge 2$.
Since $G$ satisfies condition (c') and $G_2$ is non-separable,
$b(G_1)$ must be even. 
Because  $G$ satisfies condition (c') again, $b(G_2/u_1u_2)$ 
and $b(G_1)$ should have the same parity and so $b(G_2/u_1u_2)$ 
must be even.

Let $e'$ denote an edge joining $u_1$ and $u_2$. 
Thus $G_i+u_1u_2$ can be written as $G_i+e'$.
Note that $G_1+e'$ is non-separable, implying that
the following statement is true:
\begin{quote}
for any non-separable subgraph $H$ of $G_2+e'$ 
(or $(G_2+e')/v_1v_2$ for any vertices $v_1,v_2$ in $G_2+e'$) 
with $e'\in E(H)$, if $|E(H)|\ge 2$, then 
the subgraph obtained from $H$ by replacing 
$e'$ by $G_1$ with vertex $u_i$ of $G_1$ being identified 
with $u_i$ in $H$ for $i=1,2$ is also non-separable. 
\end{quote}
Because $G$ satisfies  conditions (a'), (b') and (c')
and the above statement  holds, 
to show that $G_2+e'$ satisfies conditions (a'), (b') and (c'),
it suffices to show that it satisfies conditions (a') and (b')
for the edge $e'$.

Observe that deleting $e'$ from $G_2\Delta +e'$ obtains $G_2$ which is non-separable
by the given condition.
Also $(G_2+e')/e'=G_2/u_1u_2$ has even blocks. 
Thus $G_2+e'$ satisfies conditions (a') and (b')
for the edge $e'$. 
\proofend

By the definition of $\Theta$, $\Theta$ has only one graph (i.e., $Z_3$)
with two vertices, one graph with three vertices 
and one graph with four vertices
respectively, as shown in Figure~\ref{f5}.

%%&&&%%%%%%%%begin of Figure %%%%%%%%%%
\begin{figure}[h!]
\begin{center}
%\input f5.pic
%TeXCAD (http://texcad.sf.net/) Picture. File: [f5.pic]. Options on following lines.
%\grade{\on}
%\emlines{\off}
%\epic{\off}
%\beziermacro{\on}
%\reduce{\on}
%\snapping{\off}
%\pvinsert{% Your \input, \def, etc. here}
%\quality{8.000}
%\graddiff{0.005}
%\snapasp{1}
%\zoom{4.0000}
\unitlength 1mm % = 2.845pt
\linethickness{0.4pt}
\ifx\plotpoint\undefined\newsavebox{\plotpoint}\fi % GNUPLOT compatibility
\begin{picture}(113.018,42.518)(0,0)
\put(-130.5,-93.25){\circle*{3.354}}
\put(-82,-93.25){\circle*{3.354}}
\put(-81.75,-93.75){\circle*{3.5}}
\put(12.5,17.5){\line(1,0){22.75}}
\put(88.25,26.25){\line(1,0){22.75}}
\put(12.5,17.25){\circle*{3.536}}
\put(88.25,26){\circle*{3.536}}
\put(35.5,17.5){\circle*{3.536}}
\put(111.25,26.25){\circle*{3.536}}
\put(24,32){\circle*{3.536}}
\put(99.75,40.75){\circle*{3.536}}
\put(99.75,12.25){\circle*{3.536}}
%\emline(12.25,17.5)(23.75,32.25)
\multiput(12.25,17.5)(.0337243402,.043255132){341}{\line(0,1){.043255132}}
%\end
%\emline(88,26.25)(99.5,41)
\multiput(88,26.25)(.0337243402,.043255132){341}{\line(0,1){.043255132}}
%\end
%\emline(23.75,32.25)(35.5,17.5)
\multiput(23.75,32.25)(.0336676218,-.0422636103){349}{\line(0,-1){.0422636103}}
%\end
%\emline(99.5,41)(111.25,26.25)
\multiput(99.5,41)(.0336676218,-.0422636103){349}{\line(0,-1){.0422636103}}
%\end
%\emline(111.25,26.25)(100,12.5)
\multiput(111.25,26.25)(-.0336826347,-.0411676647){334}{\line(0,-1){.0411676647}}
%\end
%\emline(100,12.5)(88,25.75)
\multiput(100,12.5)(-.0337078652,.0372191011){356}{\line(0,1){.0372191011}}
%\end
\qbezier(100,12.25)(109.375,15.375)(111.25,26)
\qbezier(35.5,17.25)(35.375,28.125)(23.75,32.5)
\qbezier(111.25,26)(111.125,36.875)(99.5,41.25)
\qbezier(23.75,31.75)(13.125,29)(12,17.25)
\qbezier(99.5,40.5)(88.875,37.75)(87.75,26)
\qbezier(88,26)(89.125,15.125)(99.75,12.75)
\qbezier(12.25,17.75)(21.5,9.375)(35.75,17.5)
\put(21.25,4.75){(a)}
\put(99.5,5.5){(b)}
\end{picture}

\end{center}
\caption{\relabel{f5} The only graphs in $\Theta$ 
with $3$ or $4$ vertices}
\end{figure}
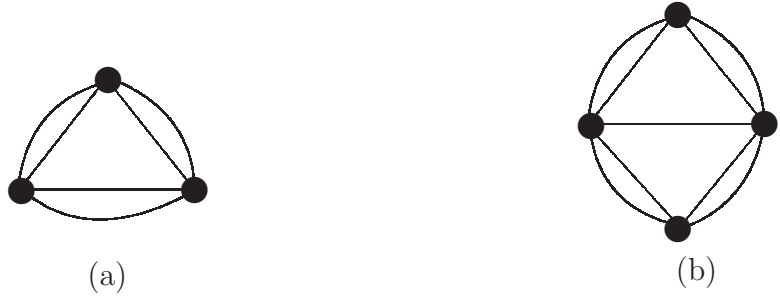

%The second target in this section is to show that $\xi_k$ can be determined by considering the set of graphs in $\Theta$ (i.e., $\Phi$) with exactly $k$ vertices.

It can be verified easily 
that every graph in $\Theta$ satisfies conditions 
(a'), (b') and (c') and thus $\Theta\subseteq \Phi$.
To show that $\Phi=\Theta$, we will prove by induction that 
every graph of $\Phi$ also belongs to $\Theta$.

Let $\Gamma(G)$ be the set of vertices $x$ in $G$ 
such that $d(x)=4$ and $|N(x)|=2$.
If $G\in \Phi$ and $|V(G)|\ge 3$, Lemma~\ref{graphs in phi}
implies that 
there are at most two parallel edges joining any two vertices in a graph 
of $\Phi$. Then for each $x\in \Gamma(G)$ with $N(x)=\{u_1,u_2\}$, 
there are exactly two parallel edges joining $x$ and $u_i$
for $i=1,2$.

\begin{lem}\quad
\relabel{on-phi'}
For any $G=(V,E)\in\Theta$, if $|V|\ge 3$, 
then $\delta(G)=4$; 
and if $|V|\ge 4$, then there are two 
non-adjacent vertices in $\Gamma(G)$.
\end{lem}

\proof We will prove this result by induction on $|V|$. 
By definition, the two graphs in  Figure~\ref{f5}
are the only graphs in $\Theta$ with three and four vertices 
respectively.  Thus the result holds when $|V|\le 4$.

Let $G=(V,E)\in \Phi$ with $|V|\ge 4$.
Assume that the result holds for $G$.
It is clear that $\delta(G(e))=4$ by the definition of $G(e)$
and the assumption that $\delta(G)=4$.

Assume that $u_1,u_2$ are the two ends of $e$.
As the result holds for $G$, 
there exists $w_1\in \Gamma(G)\setminus\{u_1,u_2\}$.
It is clear that $w_1\in \Gamma(G(e))$.
By the definition of $G(e)$, the new vertex $w$ of $G(e)$ 
is not adjacent to $w_1$ and also belongs to $\Gamma(G(e))$.
Thus the result holds for $G(e)$.
\proofend

Now we are going to prove that $\Phi$ and $\Theta$ are actually the same set.

\begin{theo}\relabel{equal-sets}
$\Phi=\Theta$.
\end{theo}

\proof It is easy to verify recursively 
that every graph in $\Theta$ satisfies conditions (a'), (b') and (c')
and so $\Theta\subseteq \Phi$. 

We will prove by induction on the number of vertices that 
every graph $G$ of $\Phi$ belongs to $\Theta$.
If $|V(G)|=2$, then $G=Z_3$ and so $G\in \Theta$.  
Assume that every graph of $\Phi$ with less than $m$ 
vertices belongs to $\Theta$, where $m\ge 3$.
Now let $G=(V,E)$ be a graph of $\Phi$ with $|V|=m$.
We first show that $\Gamma(G)\ne \emptyset$.

Assume that $u_1$ and $u_2$ are adjacent vertices in $G$.
As $G$ satisfies condition (b'), $b(G/u_1u_2)$ must be odd
and so $b(G/u_1u_2)=3$ by Lemma~\ref{graphs in phi}.
Then $G$ has the structure shown in Figure~\ref{f6}(a),
where $G_1$ and $G_2$ are two connected subgraphs of $G$ 
such that $V(G_1)\cap V(G_2)=\{u_1,u_2\}$,
$V(G_1)\cup V(G_2)=V(G)$,
$E(G_1)\cap E(G_2)=\emptyset$
and $E(G_1)\cup E(G_2)=E(G)\setminus \{e\}$,
where $e$ is an edge of $G$ joining $u_1$ and $u_2$.
\begin{figure}[h!]
\centering 
\scalebox{0.8}%{\input{f6.pic}}
{
%TeXCAD (http://texcad.sf.net/) Picture. File: [f8.pic]. Options on following lines.
%\grade{\on}
%\emlines{\off}
%\epic{\off}
%\beziermacro{\on}
%\reduce{\on}
%\snapping{\off}
%\pvinsert{% Your \input, \def, etc. here}
%\quality{8.000}
%\graddiff{0.005}
%\snapasp{1}
%\zoom{4.0000}
\unitlength 1mm % = 2.845pt
\linethickness{0.4pt}
\ifx\plotpoint\undefined\newsavebox{\plotpoint}\fi % GNUPLOT compatibility
\begin{picture}(151.625,66.75)(0,0)
\put(35,54.25){\circle*{4.61}}
\put(106.75,54.25){\circle*{4.61}}
\put(85,32.5){\circle*{4.61}}
\put(35.25,54.5){\circle*{4.61}}
\put(107,54.5){\circle*{4.61}}
\put(35,14.5){\circle*{4.61}}
\put(106.75,14.5){\circle*{4.61}}
\put(35.25,14.75){\circle*{4.61}}
\put(107,14.75){\circle*{4.61}}
\qbezier(35.25,54.75)(20.75,31.375)(35.25,14.5)
\qbezier(35,55)(49.5,31.625)(35,14.75)
\qbezier(106.75,55)(121.25,31.625)(106.75,14.75)
\qbezier(35,54)(-9.625,32.125)(35.25,14.75)
\qbezier(35.25,54.25)(79.875,32.375)(35,15)
\qbezier(107,54.25)(151.625,32.375)(106.75,15)
%\dashline{1}(29,48.75)(14.5,36.5)
\multiput(28.93,48.68)(-.0381579,-.0322368){19}{\line(-1,0){.0381579}}
\multiput(27.48,47.455)(-.0381579,-.0322368){19}{\line(-1,0){.0381579}}
\multiput(26.03,46.23)(-.0381579,-.0322368){19}{\line(-1,0){.0381579}}
\multiput(24.58,45.005)(-.0381579,-.0322368){19}{\line(-1,0){.0381579}}
\multiput(23.13,43.78)(-.0381579,-.0322368){19}{\line(-1,0){.0381579}}
\multiput(21.68,42.555)(-.0381579,-.0322368){19}{\line(-1,0){.0381579}}
\multiput(20.23,41.33)(-.0381579,-.0322368){19}{\line(-1,0){.0381579}}
\multiput(18.78,40.105)(-.0381579,-.0322368){19}{\line(-1,0){.0381579}}
\multiput(17.33,38.88)(-.0381579,-.0322368){19}{\line(-1,0){.0381579}}
\multiput(15.88,37.655)(-.0381579,-.0322368){19}{\line(-1,0){.0381579}}
%\end
%\dashline{1}(41.25,49)(55.75,36.75)
\multiput(41.18,48.93)(.0381579,-.0322368){19}{\line(1,0){.0381579}}
\multiput(42.63,47.705)(.0381579,-.0322368){19}{\line(1,0){.0381579}}
\multiput(44.08,46.48)(.0381579,-.0322368){19}{\line(1,0){.0381579}}
\multiput(45.53,45.255)(.0381579,-.0322368){19}{\line(1,0){.0381579}}
\multiput(46.98,44.03)(.0381579,-.0322368){19}{\line(1,0){.0381579}}
\multiput(48.43,42.805)(.0381579,-.0322368){19}{\line(1,0){.0381579}}
\multiput(49.88,41.58)(.0381579,-.0322368){19}{\line(1,0){.0381579}}
\multiput(51.33,40.355)(.0381579,-.0322368){19}{\line(1,0){.0381579}}
\multiput(52.78,39.13)(.0381579,-.0322368){19}{\line(1,0){.0381579}}
\multiput(54.23,37.905)(.0381579,-.0322368){19}{\line(1,0){.0381579}}
%\end
%\dashline{1}(113,49)(127.5,36.75)
\multiput(112.93,48.93)(.0381579,-.0322368){19}{\line(1,0){.0381579}}
\multiput(114.38,47.705)(.0381579,-.0322368){19}{\line(1,0){.0381579}}
\multiput(115.83,46.48)(.0381579,-.0322368){19}{\line(1,0){.0381579}}
\multiput(117.28,45.255)(.0381579,-.0322368){19}{\line(1,0){.0381579}}
\multiput(118.73,44.03)(.0381579,-.0322368){19}{\line(1,0){.0381579}}
\multiput(120.18,42.805)(.0381579,-.0322368){19}{\line(1,0){.0381579}}
\multiput(121.63,41.58)(.0381579,-.0322368){19}{\line(1,0){.0381579}}
\multiput(123.08,40.355)(.0381579,-.0322368){19}{\line(1,0){.0381579}}
\multiput(124.53,39.13)(.0381579,-.0322368){19}{\line(1,0){.0381579}}
\multiput(125.98,37.905)(.0381579,-.0322368){19}{\line(1,0){.0381579}}
%\end
%\dashline{1}(14.5,36.5)(14.5,36.5)
\put(14.43,36.43){\line(0,1){0}}
%\end
%\dashline{1}(86.25,36.5)(86.25,36.5)
\put(86.18,36.43){\line(0,1){0}}
%\end
%\dashline{1}(55.75,36.75)(55.75,36.75)
\put(55.68,36.68){\line(0,1){0}}
%\end
%\dashline{1}(127.5,36.75)(127.5,36.75)
\put(127.43,36.68){\line(0,1){0}}
%\end
%\dashline{1}(28.75,44.5)(14.25,32)
\multiput(28.68,44.43)(-.0381579,-.0328947){19}{\line(-1,0){.0381579}}
\multiput(27.23,43.18)(-.0381579,-.0328947){19}{\line(-1,0){.0381579}}
\multiput(25.78,41.93)(-.0381579,-.0328947){19}{\line(-1,0){.0381579}}
\multiput(24.33,40.68)(-.0381579,-.0328947){19}{\line(-1,0){.0381579}}
\multiput(22.88,39.43)(-.0381579,-.0328947){19}{\line(-1,0){.0381579}}
\multiput(21.43,38.18)(-.0381579,-.0328947){19}{\line(-1,0){.0381579}}
\multiput(19.98,36.93)(-.0381579,-.0328947){19}{\line(-1,0){.0381579}}
\multiput(18.53,35.68)(-.0381579,-.0328947){19}{\line(-1,0){.0381579}}
\multiput(17.08,34.43)(-.0381579,-.0328947){19}{\line(-1,0){.0381579}}
\multiput(15.63,33.18)(-.0381579,-.0328947){19}{\line(-1,0){.0381579}}
%\end
%\dashline{1}(41.5,44.75)(56,32.25)
\multiput(41.43,44.68)(.0381579,-.0328947){19}{\line(1,0){.0381579}}
\multiput(42.88,43.43)(.0381579,-.0328947){19}{\line(1,0){.0381579}}
\multiput(44.33,42.18)(.0381579,-.0328947){19}{\line(1,0){.0381579}}
\multiput(45.78,40.93)(.0381579,-.0328947){19}{\line(1,0){.0381579}}
\multiput(47.23,39.68)(.0381579,-.0328947){19}{\line(1,0){.0381579}}
\multiput(48.68,38.43)(.0381579,-.0328947){19}{\line(1,0){.0381579}}
\multiput(50.13,37.18)(.0381579,-.0328947){19}{\line(1,0){.0381579}}
\multiput(51.58,35.93)(.0381579,-.0328947){19}{\line(1,0){.0381579}}
\multiput(53.03,34.68)(.0381579,-.0328947){19}{\line(1,0){.0381579}}
\multiput(54.48,33.43)(.0381579,-.0328947){19}{\line(1,0){.0381579}}
%\end
%\dashline{1}(113.25,44.75)(127.75,32.25)
\multiput(113.18,44.68)(.0381579,-.0328947){19}{\line(1,0){.0381579}}
\multiput(114.63,43.43)(.0381579,-.0328947){19}{\line(1,0){.0381579}}
\multiput(116.08,42.18)(.0381579,-.0328947){19}{\line(1,0){.0381579}}
\multiput(117.53,40.93)(.0381579,-.0328947){19}{\line(1,0){.0381579}}
\multiput(118.98,39.68)(.0381579,-.0328947){19}{\line(1,0){.0381579}}
\multiput(120.43,38.43)(.0381579,-.0328947){19}{\line(1,0){.0381579}}
\multiput(121.88,37.18)(.0381579,-.0328947){19}{\line(1,0){.0381579}}
\multiput(123.33,35.93)(.0381579,-.0328947){19}{\line(1,0){.0381579}}
\multiput(124.78,34.68)(.0381579,-.0328947){19}{\line(1,0){.0381579}}
\multiput(126.23,33.43)(.0381579,-.0328947){19}{\line(1,0){.0381579}}
%\end
%\dashline{1}(14.25,32)(14.25,32)
\put(14.18,31.93){\line(0,1){0}}
%\end
%\dashline{1}(86,32)(86,32)
\put(85.93,31.93){\line(0,1){0}}
%\end
%\dashline{1}(56,32.25)(56,32.25)
\put(55.93,32.18){\line(0,1){0}}
%\end
%\dashline{1}(127.75,32.25)(127.75,32.25)
\put(127.68,32.18){\line(0,1){0}}
%\end
%\dashline{1}(27.5,38.5)(15.75,28.25)
\multiput(27.43,38.43)(-.0383987,-.0334967){18}{\line(-1,0){.0383987}}
\multiput(26.047,37.224)(-.0383987,-.0334967){18}{\line(-1,0){.0383987}}
\multiput(24.665,36.018)(-.0383987,-.0334967){18}{\line(-1,0){.0383987}}
\multiput(23.283,34.812)(-.0383987,-.0334967){18}{\line(-1,0){.0383987}}
\multiput(21.9,33.606)(-.0383987,-.0334967){18}{\line(-1,0){.0383987}}
\multiput(20.518,32.4)(-.0383987,-.0334967){18}{\line(-1,0){.0383987}}
\multiput(19.136,31.194)(-.0383987,-.0334967){18}{\line(-1,0){.0383987}}
\multiput(17.753,29.989)(-.0383987,-.0334967){18}{\line(-1,0){.0383987}}
\multiput(16.371,28.783)(-.0383987,-.0334967){18}{\line(-1,0){.0383987}}
%\end
%\dashline{1}(42.75,38.75)(54.5,28.5)
\multiput(42.68,38.68)(.0383987,-.0334967){18}{\line(1,0){.0383987}}
\multiput(44.062,37.474)(.0383987,-.0334967){18}{\line(1,0){.0383987}}
\multiput(45.444,36.268)(.0383987,-.0334967){18}{\line(1,0){.0383987}}
\multiput(46.827,35.062)(.0383987,-.0334967){18}{\line(1,0){.0383987}}
\multiput(48.209,33.856)(.0383987,-.0334967){18}{\line(1,0){.0383987}}
\multiput(49.591,32.65)(.0383987,-.0334967){18}{\line(1,0){.0383987}}
\multiput(50.974,31.444)(.0383987,-.0334967){18}{\line(1,0){.0383987}}
\multiput(52.356,30.239)(.0383987,-.0334967){18}{\line(1,0){.0383987}}
\multiput(53.739,29.033)(.0383987,-.0334967){18}{\line(1,0){.0383987}}
%\end
%\dashline{1}(114.5,38.75)(126.25,28.5)
\multiput(114.43,38.68)(.0383987,-.0334967){18}{\line(1,0){.0383987}}
\multiput(115.812,37.474)(.0383987,-.0334967){18}{\line(1,0){.0383987}}
\multiput(117.194,36.268)(.0383987,-.0334967){18}{\line(1,0){.0383987}}
\multiput(118.577,35.062)(.0383987,-.0334967){18}{\line(1,0){.0383987}}
\multiput(119.959,33.856)(.0383987,-.0334967){18}{\line(1,0){.0383987}}
\multiput(121.341,32.65)(.0383987,-.0334967){18}{\line(1,0){.0383987}}
\multiput(122.724,31.444)(.0383987,-.0334967){18}{\line(1,0){.0383987}}
\multiput(124.106,30.239)(.0383987,-.0334967){18}{\line(1,0){.0383987}}
\multiput(125.489,29.033)(.0383987,-.0334967){18}{\line(1,0){.0383987}}
%\end
%\dashline{1}(27,32.75)(18.25,24.75)
\multiput(26.93,32.68)(-.0354251,-.0323887){19}{\line(-1,0){.0354251}}
\multiput(25.584,31.449)(-.0354251,-.0323887){19}{\line(-1,0){.0354251}}
\multiput(24.237,30.218)(-.0354251,-.0323887){19}{\line(-1,0){.0354251}}
\multiput(22.891,28.987)(-.0354251,-.0323887){19}{\line(-1,0){.0354251}}
\multiput(21.545,27.757)(-.0354251,-.0323887){19}{\line(-1,0){.0354251}}
\multiput(20.199,26.526)(-.0354251,-.0323887){19}{\line(-1,0){.0354251}}
\multiput(18.853,25.295)(-.0354251,-.0323887){19}{\line(-1,0){.0354251}}
%\end
%\dashline{1}(43.25,33)(52,25)
\multiput(43.18,32.93)(.0354251,-.0323887){19}{\line(1,0){.0354251}}
\multiput(44.526,31.699)(.0354251,-.0323887){19}{\line(1,0){.0354251}}
\multiput(45.872,30.468)(.0354251,-.0323887){19}{\line(1,0){.0354251}}
\multiput(47.218,29.237)(.0354251,-.0323887){19}{\line(1,0){.0354251}}
\multiput(48.564,28.007)(.0354251,-.0323887){19}{\line(1,0){.0354251}}
\multiput(49.91,26.776)(.0354251,-.0323887){19}{\line(1,0){.0354251}}
\multiput(51.257,25.545)(.0354251,-.0323887){19}{\line(1,0){.0354251}}
%\end
%\dashline{1}(115,33)(123.75,25)
\multiput(114.93,32.93)(.0354251,-.0323887){19}{\line(1,0){.0354251}}
\multiput(116.276,31.699)(.0354251,-.0323887){19}{\line(1,0){.0354251}}
\multiput(117.622,30.468)(.0354251,-.0323887){19}{\line(1,0){.0354251}}
\multiput(118.968,29.237)(.0354251,-.0323887){19}{\line(1,0){.0354251}}
\multiput(120.314,28.007)(.0354251,-.0323887){19}{\line(1,0){.0354251}}
\multiput(121.66,26.776)(.0354251,-.0323887){19}{\line(1,0){.0354251}}
\multiput(123.007,25.545)(.0354251,-.0323887){19}{\line(1,0){.0354251}}
%\end
%\dashline{1}(27.5,27.25)(22,21.5)
\multiput(27.43,27.18)(-.0321637,-.0336257){19}{\line(0,-1){.0336257}}
\multiput(26.207,25.902)(-.0321637,-.0336257){19}{\line(0,-1){.0336257}}
\multiput(24.985,24.624)(-.0321637,-.0336257){19}{\line(0,-1){.0336257}}
\multiput(23.763,23.346)(-.0321637,-.0336257){19}{\line(0,-1){.0336257}}
\multiput(22.541,22.069)(-.0321637,-.0336257){19}{\line(0,-1){.0336257}}
%\end
%\dashline{1}(42.75,27.5)(48.25,21.75)
\multiput(42.68,27.43)(.0321637,-.0336257){19}{\line(0,-1){.0336257}}
\multiput(43.902,26.152)(.0321637,-.0336257){19}{\line(0,-1){.0336257}}
\multiput(45.124,24.874)(.0321637,-.0336257){19}{\line(0,-1){.0336257}}
\multiput(46.346,23.596)(.0321637,-.0336257){19}{\line(0,-1){.0336257}}
\multiput(47.569,22.319)(.0321637,-.0336257){19}{\line(0,-1){.0336257}}
%\end
%\dashline{1}(114.5,27.5)(120,21.75)
\multiput(114.43,27.43)(.0321637,-.0336257){19}{\line(0,-1){.0336257}}
\multiput(115.652,26.152)(.0321637,-.0336257){19}{\line(0,-1){.0336257}}
\multiput(116.874,24.874)(.0321637,-.0336257){19}{\line(0,-1){.0336257}}
\multiput(118.096,23.596)(.0321637,-.0336257){19}{\line(0,-1){.0336257}}
\multiput(119.319,22.319)(.0321637,-.0336257){19}{\line(0,-1){.0336257}}
%\end
%\dashline{1}(22,21.5)(22,21.5)
\put(21.93,21.43){\line(0,1){0}}
%\end
%\dashline{1}(93.75,21.5)(93.75,21.5)
\put(93.68,21.43){\line(0,1){0}}
%\end
%\dashline{1}(48.25,21.75)(48.25,21.75)
\put(48.18,21.68){\line(0,1){0}}
%\end
%\dashline{1}(120,21.75)(120,21.75)
\put(119.93,21.68){\line(0,1){0}}
%\end
%\dashline{1}(28,22.75)(25.5,19.5)
\multiput(27.93,22.68)(-.0333333,-.0433333){15}{\line(0,-1){.0433333}}
\multiput(26.93,21.38)(-.0333333,-.0433333){15}{\line(0,-1){.0433333}}
\multiput(25.93,20.08)(-.0333333,-.0433333){15}{\line(0,-1){.0433333}}
%\end
%\dashline{1}(42.25,23)(44.75,19.75)
\multiput(42.18,22.93)(.0333333,-.0433333){15}{\line(0,-1){.0433333}}
\multiput(43.18,21.63)(.0333333,-.0433333){15}{\line(0,-1){.0433333}}
\multiput(44.18,20.33)(.0333333,-.0433333){15}{\line(0,-1){.0433333}}
%\end
%\dashline{1}(114,23)(116.5,19.75)
\multiput(113.93,22.93)(.0333333,-.0433333){15}{\line(0,-1){.0433333}}
\multiput(114.93,21.63)(.0333333,-.0433333){15}{\line(0,-1){.0433333}}
\multiput(115.93,20.33)(.0333333,-.0433333){15}{\line(0,-1){.0433333}}
%\end
%\dashline{1}(25.5,19.5)(25.5,19.5)
\put(25.43,19.43){\line(0,1){0}}
%\end
%\dashline{1}(97.25,19.5)(97.25,19.5)
\put(97.18,19.43){\line(0,1){0}}
%\end
%\dashline{1}(44.75,19.75)(44.75,19.75)
\put(44.68,19.68){\line(0,1){0}}
%\end
%\dashline{1}(116.5,19.75)(116.5,19.75)
\put(116.43,19.68){\line(0,1){0}}
%\end
\put(34.25,62.25){\makebox(0,0)[cc]{$u_1$}}
\put(106,62.25){\makebox(0,0)[cc]{$u_1$}}
\put(33.75,7.75){\makebox(0,0)[cc]{$u_2$}}
\put(105.5,7.75){\makebox(0,0)[cc]{$u_2$}}
\put(15.5,66.5){\makebox(0,0)[]{}}
\put(87.25,66.5){\makebox(0,0)[]{}}
\put(54.75,66.75){\makebox(0,0)[cc]{}}
\put(126.5,66.75){\makebox(0,0)[cc]{}}
\put(34,.75){\makebox(0,0)[cc]{(a) $G$}}
\put(105.75,.75){\makebox(0,0)[cc]{(b)  $G$}}
\put(6.25,33){$G_1$}
\put(78,33){$w$}
\put(63.75,34.5){$G_2$}
\put(135.5,34.5){$G_2$}
\put(35,54.75){\line(0,-1){40.25}}
\put(106.75,54.75){\line(0,-1){40.25}}
\put(36.5,34.75){$e$}
\put(108.25,34.75){$e$}
\qbezier(106.75,54.75)(84,50.125)(85.25,33)
\qbezier(85.25,33)(100.375,38.75)(107,54.5)
\qbezier(85.25,32.75)(87,12.75)(106.75,14.75)
\qbezier(106.75,14.75)(102.5,28.25)(85.25,32.75)
%\dashline{1}(87.25,42.25)(97.5,50.5)
\multiput(87.18,42.18)(.0401961,.0323529){17}{\line(1,0){.0401961}}
\multiput(88.546,43.28)(.0401961,.0323529){17}{\line(1,0){.0401961}}
\multiput(89.913,44.38)(.0401961,.0323529){17}{\line(1,0){.0401961}}
\multiput(91.28,45.48)(.0401961,.0323529){17}{\line(1,0){.0401961}}
\multiput(92.646,46.58)(.0401961,.0323529){17}{\line(1,0){.0401961}}
\multiput(94.013,47.68)(.0401961,.0323529){17}{\line(1,0){.0401961}}
\multiput(95.38,48.78)(.0401961,.0323529){17}{\line(1,0){.0401961}}
\multiput(96.746,49.88)(.0401961,.0323529){17}{\line(1,0){.0401961}}
%\end
%\dashline{1}(86,38.5)(103.5,52.5)
\multiput(85.93,38.43)(.0405093,.0324074){18}{\line(1,0){.0405093}}
\multiput(87.388,39.596)(.0405093,.0324074){18}{\line(1,0){.0405093}}
\multiput(88.846,40.763)(.0405093,.0324074){18}{\line(1,0){.0405093}}
\multiput(90.305,41.93)(.0405093,.0324074){18}{\line(1,0){.0405093}}
\multiput(91.763,43.096)(.0405093,.0324074){18}{\line(1,0){.0405093}}
\multiput(93.221,44.263)(.0405093,.0324074){18}{\line(1,0){.0405093}}
\multiput(94.68,45.43)(.0405093,.0324074){18}{\line(1,0){.0405093}}
\multiput(96.138,46.596)(.0405093,.0324074){18}{\line(1,0){.0405093}}
\multiput(97.596,47.763)(.0405093,.0324074){18}{\line(1,0){.0405093}}
\multiput(99.055,48.93)(.0405093,.0324074){18}{\line(1,0){.0405093}}
\multiput(100.513,50.096)(.0405093,.0324074){18}{\line(1,0){.0405093}}
\multiput(101.971,51.263)(.0405093,.0324074){18}{\line(1,0){.0405093}}
%\end
%\dashline{1}(87.25,36.5)(103.75,49.5)
\multiput(87.18,36.43)(.0421995,.0332481){17}{\line(1,0){.0421995}}
\multiput(88.614,37.56)(.0421995,.0332481){17}{\line(1,0){.0421995}}
\multiput(90.049,38.691)(.0421995,.0332481){17}{\line(1,0){.0421995}}
\multiput(91.484,39.821)(.0421995,.0332481){17}{\line(1,0){.0421995}}
\multiput(92.919,40.951)(.0421995,.0332481){17}{\line(1,0){.0421995}}
\multiput(94.354,42.082)(.0421995,.0332481){17}{\line(1,0){.0421995}}
\multiput(95.788,43.212)(.0421995,.0332481){17}{\line(1,0){.0421995}}
\multiput(97.223,44.343)(.0421995,.0332481){17}{\line(1,0){.0421995}}
\multiput(98.658,45.473)(.0421995,.0332481){17}{\line(1,0){.0421995}}
\multiput(100.093,46.604)(.0421995,.0332481){17}{\line(1,0){.0421995}}
\multiput(101.528,47.734)(.0421995,.0332481){17}{\line(1,0){.0421995}}
\multiput(102.962,48.864)(.0421995,.0332481){17}{\line(1,0){.0421995}}
%\end
%\dashline{1}(90,29.75)(102.75,19.75)
\multiput(89.93,29.68)(.0416667,-.0326797){17}{\line(1,0){.0416667}}
\multiput(91.346,28.569)(.0416667,-.0326797){17}{\line(1,0){.0416667}}
\multiput(92.763,27.457)(.0416667,-.0326797){17}{\line(1,0){.0416667}}
\multiput(94.18,26.346)(.0416667,-.0326797){17}{\line(1,0){.0416667}}
\multiput(95.596,25.235)(.0416667,-.0326797){17}{\line(1,0){.0416667}}
\multiput(97.013,24.124)(.0416667,-.0326797){17}{\line(1,0){.0416667}}
\multiput(98.43,23.013)(.0416667,-.0326797){17}{\line(1,0){.0416667}}
\multiput(99.846,21.902)(.0416667,-.0326797){17}{\line(1,0){.0416667}}
\multiput(101.263,20.791)(.0416667,-.0326797){17}{\line(1,0){.0416667}}
%\end
%\dashline{1}(88,28.25)(102.75,16.75)
\multiput(87.93,28.18)(.0409722,-.0319444){18}{\line(1,0){.0409722}}
\multiput(89.405,27.03)(.0409722,-.0319444){18}{\line(1,0){.0409722}}
\multiput(90.88,25.88)(.0409722,-.0319444){18}{\line(1,0){.0409722}}
\multiput(92.355,24.73)(.0409722,-.0319444){18}{\line(1,0){.0409722}}
\multiput(93.83,23.58)(.0409722,-.0319444){18}{\line(1,0){.0409722}}
\multiput(95.305,22.43)(.0409722,-.0319444){18}{\line(1,0){.0409722}}
\multiput(96.78,21.28)(.0409722,-.0319444){18}{\line(1,0){.0409722}}
\multiput(98.255,20.13)(.0409722,-.0319444){18}{\line(1,0){.0409722}}
\multiput(99.73,18.98)(.0409722,-.0319444){18}{\line(1,0){.0409722}}
\multiput(101.205,17.83)(.0409722,-.0319444){18}{\line(1,0){.0409722}}
%\end
%\dashline{1}(88.5,24.25)(99,16.25)
\multiput(88.43,24.18)(.04375,-.0333333){16}{\line(1,0){.04375}}
\multiput(89.83,23.113)(.04375,-.0333333){16}{\line(1,0){.04375}}
\multiput(91.23,22.046)(.04375,-.0333333){16}{\line(1,0){.04375}}
\multiput(92.63,20.98)(.04375,-.0333333){16}{\line(1,0){.04375}}
\multiput(94.03,19.913)(.04375,-.0333333){16}{\line(1,0){.04375}}
\multiput(95.43,18.846)(.04375,-.0333333){16}{\line(1,0){.04375}}
\multiput(96.83,17.78)(.04375,-.0333333){16}{\line(1,0){.04375}}
\multiput(98.23,16.713)(.04375,-.0333333){16}{\line(1,0){.04375}}
%\end
\put(93.5,44){\makebox(0,0)[cc]{$H_1$}}
\put(94,22.25){\makebox(0,0)[cc]{$H_2$}}
\end{picture}
}

\caption{\relabel{f6} Graph $G$}
\end{figure}

Since $|V|\ge 3$, we may assume that $|V(G_1)|\ge 3$.
As $b(G_2+e)=1$ and $G$ satisfies condition (c'), $b(G_1)$ is even.
Thus $G_1$ can be divided into two edge-disjoint subgraphs 
$H_1$ and $H_2$ such that 
$V(H_1)\cap V(H_2)=\{v\}$,
$V(H_1)\cup V(H_2)=V(G_1)$,
and $E(G_1)\cup E(G_2)=E(G_1)$,
as shown in Figure~\ref{f6}(b).
By Lemma~\ref{graphs in phi}, it can be deduced  
that $b(G_1+e)+b(G_2)=3$, implying that both 
$H_1$ and $H_2$ are non-separable. 
As $G$ satisfies condition (a'), we have $|E(H_i)|\ge 2$.

If $|V|=3$, then each $H_i$ has exactly two edges 
and $G_2$ has just one edge, and so $G$ is the graph $Z_3(e')$ 
for some edge $e'$ in $Z_3$, and hence $\Gamma(G)=V(G)$.
Now assume that $|V|\ge 4$.
At least one of the three subgraphs $H_1, H_2$ and $G_2+e$ 
contains at least three vertices. 

Consider the case that $G_2+e$ has at least three vertices.
Lemma~\ref{new graphs in phi} implies that 
$G_2+e+u_1u_2\in \Phi$.
Since this graph has less vertices than $G$, 
by inductive assumption, $G_2+e+u_1u_2\in\Theta$.
Then, by Lemma~\ref{on-phi'}, 
there exists $x\in \Gamma(G_2+e+u_1u_2)\setminus \{u_1,u_2\}$.
It is clear that $x\in \Gamma(G)$.

Now assume that $V(G_2)=\{u_1,u_2\}$.
As $G/u_1u_2$ has exactly 3 blocks by Lemma~\ref{graphs in phi},
$G_2$ is the graph $Z_1$.
Then we may assume that $|V(H_1)|\ge 3$.
Lemma~\ref{new graphs in phi} implies that 
$H_1+u_1w\in \Phi$.
Since $H_1+u_1w$ has less vertices than $G$, 
by inductive assumption, $H_1+u_1w\in\Theta$.
Thus the graph $(H_1+u_1w)(f)\in \Theta$ by the definition of $\Theta$,
where $f$ is an edge of $H_1+u_1w$ joining $u_1$ and $w$.
By Lemma~\ref{on-phi'}, either $u_1\in \Gamma((H_1+u_1w)(f))$
or there is $x\in \Gamma((H_1+u_1w)(f))\setminus \{w,u_1\}$.
Thus either $u_1\in \Gamma(G)$ or $x\in \Gamma(G)$.

Hence  $\Gamma(G)\ne \emptyset$.
Let $w$ be a vertex in $\Gamma(G)$.
Let $N(w)=\{v_1,v_2\}$.
Then there are exactly two parallel edges joining $w$ and $v_i$
for $i=1,2$.
Since $G$ satisfies condition (c'), $G-w$ is non-separable and 
has at least two edges.
By Lemma~\ref{new graphs in phi}, $(G-w)+v_1v_2\in \Phi$.
By inductive assumption, $(G-w)+v_1v_2\in\Theta$.
Hence $G\in \Theta$ by the definition of $\Theta$.
\proofend

By Theorem~\ref{equal-sets} and the definition of $\Theta$, 
we have $\Phi_0=\Phi_1=\emptyset$, 
$\Phi_2=\{Z_3\}$
and $\Phi_{k+1}=\{G(e): G\in \Phi_k, e\in E(G)\}$
for $k\ge 2$. 

Now it remains to show that 
$\eta(\Phi_0), \eta(\Phi_1),\eta(\Phi_2),\eta(\Phi_3), \cdots$ is non-increasing and so Lemma~\ref{determine xik} implies that 
$\xi_k=\eta(\Phi_k)$.

\begin{theo}\quad
\relabel{free-interval}
$\eta(\Phi_0),\eta(\Phi_1),
\eta(\Phi_2), \eta(\Phi_3),  \cdots$
is a non-increasing sequence 
and $\xi_k=\eta(\Phi_k)$ for $k=0,1,2,\cdots$.
\end{theo}

\proof 
Since $\Phi_2$ has only one graph, i.e., $Z_3$, 
we have $\xi_2=\eta(\Phi_2)=2$ by Lemma~\ref{determine xik}.
Thus $\xi_i=2=\eta(\Phi_i)$ for $i=0,1,2$. 
We need to apply the following claim.

\noindent {\bf Claim A}:
For $k\ge 2$, if $\xi_k=\eta(\Phi_k)$, then 
$\eta(\Phi_{k+1})\le \eta(\Phi_k)$.

Suppose that $\eta(\Phi_k)< \eta(\Phi_{k+1})$.
Then there exists $G\in \Phi_k$ 
such that $\eta(G)=\eta(\Phi_k)$.  
As $\eta(\Phi_{k+1})\le 2$, $\eta(\Phi_k)<2$, and so 
$F(G,\eta(\Phi_k))=0$.

Let $e$ be an edge of $G$ joining two vertices 
$u_1$ and $u_2$. 
By Lemma~\ref{2-blocks}, we have  
\beeq
F(G(e),\lambda)=(\lambda-1)^2F(G-e,\lambda)
+(\lambda-2)^2F(G,\lambda)
\eneq
and 
\beeq
F(G+u_1u_2,\lambda)=(\lambda-1)F(G-e,\lambda)
+(\lambda-2)F(G,\lambda).
\eneq
Thus
\beeq
F(G(e),\lambda)=(\lambda-1)F(G+u_1u_2,\lambda)
+(2-\lambda)F(G,\lambda).
\eneq
Since $G, G+u_1u_2$ and $G(e)$ are all non-separable, 
$p(G), p(G(e))$ and $p(G+u_1u_2)-1$ all have the same parity.
Thus 
\beeq
(2-\lambda)Q(G,\lambda)=Q(G(e),\lambda)+(\lambda-1)Q(G+u_1u_2,\lambda).
\eneq
As $G\in \Phi_k$, we have $G(e)\in \Phi_{k+1}$.
Since $\eta(\Phi_{k+1})>\eta(\Phi_k)$ by assumption, 
we have $Q(G(e),\eta(\Phi_k))>0$.
As $G+u_1u_2\in \Psi_k$ and $\xi_k=\eta(\Phi_k)$, we have  
$Q(G+u_1u_2,\eta(\Phi_k))\ge 0$. 
Hence $Q(G,\eta(\Phi_k))>0$, contradicting the assumption that 
$F(G,\eta(\Phi_k))=0$.

So Claim A holds.
Now assume that for integer $k$ with $k\ge 2$, %the sequence 
$\eta(\Phi_0),\eta(\Phi_1),
\eta(\Phi_2),  \cdots,\eta(\Phi_k)$
is non-increasing 
and $\xi_i=\eta(\Phi_i)$ for $i=0,1,2,\cdots,k$.
By Claim A, $\eta(\Phi_{k+1})\le \eta(\Phi_k)$.
Then, Lemma~\ref{determine xik} implies that 
$\xi_{k+1}=\eta(\Phi_{k+1})$.
Hence this theorem holds.
\proofend

Before the end of this section, we try to find the values 
of $\xi_k$ (i.e., $\eta(\Phi_k)$) for some $k$.
By Theorem~\ref{free-interval} and the fact that 
$\Phi_{k+1}=\{G(e): G\in \Phi_k, e\in E(G)\}$ for $k\ge 2$, 
it is not hard to find the value of $\xi_k$ for small $k$.
As an example, we will determine $\xi_k$ for $0\le k\le 5$.

\begin{theo}
\relabel{some value of eta}
$\xi_k=2$ for $k=0,1,2$, 
$\xi_3=1.430159709\cdots$,
$\xi_4=1.361103081\cdots$,
and $\xi_5=1.317672196\cdots$,
where the last three numbers 
are the real roots of $\lambda^3-5\lambda^2+10\lambda-7$,
$\lambda^3-4\lambda^2+8\lambda-6$ 
and $\lambda^3-6\lambda^2+13\lambda-9$
in $(1,2)$ respectively.
\end{theo}

\proof By Theorem~\ref{free-interval}, $\xi_k=\eta(\Phi_k)$.
As $Z_3$ is the only graph of $\Phi_2$, we have 
$\xi_2=\eta(\Phi_2)=2$. Thus $\xi_0=\xi_1=2$.
Note that the two graphs in Figure~\ref{f5} are the only graphs of 
$\Phi_3$ and $\Phi_4$.
Their flow polynomials are  
\beeq
(\lambda-1)(\lambda^3-5\lambda^2+10\lambda-7),
\eneq
and 
\beeq
(\lambda-1)(\lambda-2)^2
(\lambda^3-4\lambda^2+8\lambda-6).
\eneq
Each of the above polynomials has only one real root in $(1,2)$:
\beeq
1.430159709\cdots\quad \mbox{and}\quad 1.361103081\cdots.
\eneq
Thus the result holds for $\xi_3$ and $\xi_4$.
Because $\Phi=\Theta$, 
$\Phi_5$ has only two different graphs, as shown in Figure~\ref{f7}.
\begin{figure}[h!]
\centering 
\scalebox{0.9}%{\input{f7.pic}}
{
%TeXCAD (http://texcad.sf.net/) Picture. File: [f7.pic]. Options on following lines.
%\grade{\on}
%\emlines{\off}
%\epic{\off}
%\beziermacro{\on}
%\reduce{\on}
%\snapping{\off}
%\pvinsert{% Your \input, \def, etc. here}
%\quality{8.000}
%\graddiff{0.005}
%\snapasp{1}
%\zoom{4.0000}
\unitlength 1mm % = 2.845pt
\linethickness{0.4pt}
\ifx\plotpoint\undefined\newsavebox{\plotpoint}\fi % GNUPLOT compatibility
\begin{picture}(92.375,47)(0,0)
\put(17.5,44){\circle*{3}}
\put(75.25,44.75){\circle*{3}}
\put(17.5,12.75){\circle*{3}}
\put(75.25,13.5){\circle*{3}}
\put(3,29.25){\circle*{3}}
\put(60.5,30){\circle*{3}}
\put(17.5,29.5){\circle*{3}}
\put(87.25,41.25){\circle*{3}}
\put(32.5,29.5){\circle*{3}}
\put(90.25,30.25){\circle*{3}}
\qbezier(32.75,29.5)(21.375,34.25)(17.5,44)
\qbezier(90.5,30.25)(79.125,35)(75.25,44.75)
\qbezier(32.75,29.75)(21.5,24.75)(17.25,29.75)
\qbezier(17.25,29.75)(24.25,31.375)(32.25,29.5)
\qbezier(32.25,29.5)(19.25,20.625)(17.25,13.25)
\qbezier(90,30.25)(77,21.375)(75,14)
\qbezier(17.25,13.25)(27,18)(32.75,29.75)
\qbezier(75,14)(84.75,18.75)(90.5,30.5)
\qbezier(86.5,41.25)(80.25,40)(75,44.75)
\qbezier(75,44.75)(81.5,47)(87,41.25)
\qbezier(87,41.25)(85.875,32.75)(90.25,30.25)
\qbezier(90.25,30.25)(92.375,36.625)(87,41.5)
\put(17,3.5){\makebox(0,0)[cc]{(a)}}
\put(75.25,3.5){\makebox(0,0)[cc]{(b)}}
\qbezier(60.25,30)(63,42.5)(74.75,45)
\qbezier(3,29.5)(5.75,42)(17.5,44.5)
\qbezier(17.25,44.25)(12.5,35.75)(2.75,29.25)
\qbezier(74.75,45)(70,36.5)(60.25,30)
\qbezier(2.75,29.25)(10.875,33)(17.5,29.75)
\qbezier(17.5,29.75)(8.875,27)(2.75,29.25)
\qbezier(2.75,29.25)(6.375,15.5)(17.5,12.75)
\qbezier(60.25,30)(63.875,16.25)(75,13.5)
\qbezier(17.5,12.75)(12.5,20)(2.5,29.25)
\qbezier(75,13.5)(70,20.75)(60,30)
%\emline(88.25,28.5)(89.5,29.75)
\multiput(88.25,28.5)(.03289474,.03289474){38}{\line(0,1){.03289474}}
%\end
\put(90.25,30.25){\line(-1,0){29.75}}
\put(13.25,5){\line(0,1){.25}}
\qbezier(17.5,44)(29.5,40.125)(32.5,29.75)
\end{picture}
}

\caption{\relabel{f7} The only two graphs in $\Phi_5$}
\end{figure}
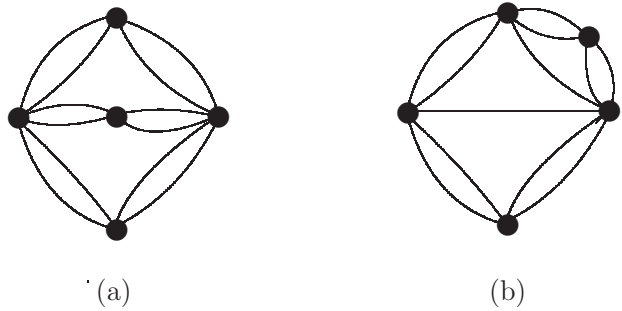

Their flow polynomials are 
\beeq
(\lambda-1)(\lambda^3-6\lambda^2+13\lambda-9)
(\lambda^4-5\lambda^3+12\lambda^2-16\lambda+9),
\eneq
\beeq
(\lambda-1)(\lambda-2)
(\lambda^6-9\lambda^5+37\lambda^4-89\lambda^3+132\lambda^2-112\lambda+41).
\eneq
Their smallest roots in $(1,2)$ are $1.317672196\cdots $ and 
$1.335087886\cdots$ respectively. 
Thus the result holds.
\proofend

%\newpage

\resection{Integral Flow Roots}\relabel{sec6}

It is known that there exist graphs whose chromatic roots are all 
integers, for example, chordal graphs. 
There are also graphs which have all real chromatic roots 
but also include non-integral chromatic roots.
For any integer $n$ with $n\ge 2$, let $H_n$ be the graph 
obtained from the complete graph $K_n$ by 
subdividing some edge in $K_n$ once.
Observe that 
\beeq
P(H_n,\lambda)=\lambda (\lambda-1)\cdots (\lambda -n+2)
(\lambda^2-n\lambda+2n-3).
\eneq 
When $n\ge 7$, all roots of $P(H_n,\lambda)$ are real,
but some roots are not integral. 

In this section we consider the problem of 
whether there is a graph whose flow roots also 
have similar properties, 
i.e., all flow roots are real but some of them are not integral.
We shall show that if there is such a graph $G=(V,E)$,
then this graph must satisfy various conditions 
(see Theorem~\ref{main-lem1})
%, one of which is $\max\{|V|+17, 2|V|+4\}\le |E|\le (32|V|-49)/5$. 

%For a connected graph $G$ and $E'\subseteq E(G)$, if $G-E'$ is disconnected, then $E'$ is called an edge-cut  or a $k$-edge-cut, where $k=|E'|$. 
Let $G=(V,E)$ be a bridgeless connected graph. 
If $G$ has no $2$-edge-cut, 
it can be  proved by induction and by applying (\ref{eq1-2})
that $F(G,\lambda)$ is a polynomial of order $r$, where 
$r=|E|-|V|+1$,
and if $F(G,\lambda)=\sum\limits_{0\le i\le r}b_i\lambda^i$, 
then 
\beeq\relabel{leading-coe-eq}
b_r=1, b_{r-1}=-|E|\mbox{ and }b_{r-2}={|E|\choose 2}-\gamma,
\eneq
where $\gamma$ is the number of $3$-edge-cuts of $G$. 
Applying the technique used in the proof of 
Lemma 4.2 in \cite{kung}, 
a lower bound on $\gamma$ in terms of $|E|$ 
and $r$ can be obtained. 
We need to apply the following result whose proof 
can be found in \cite{kung}.

\begin{lem}[\cite{kung}]\quad 
\relabel{le0}
Assume that the polynomial 
\beeq
P(\lambda)=\sum_{i=0}^n (-1)^ia_i\lambda ^{n-i},
\eneq
where $a_0=1$, has only positive real roots.
Then for each $i: 2\le i\le n$, 
\beeq
0<a_i\le {n\choose i}\left (\frac {a_1} n\right )^i,
\eneq
where equality holds if and only if 
$P(\lambda)=(\lambda -a_1/n)^n$.
\end{lem}

\begin{lem}
\quad \relabel{le1}
Let $G=(V,E)$ be a bridgeless connected graph
which has no $2$-edge-cut. 
Assume that all roots of $F(G,\lambda)$ are real numbers.
Let $\gamma$ be the number of $3$-edge-cuts of $G$.
Then
\beeq\relabel{le-eq1}
\gamma\ge \frac{(|E|-r)(|E|-1)}{2(r-1)},
\eneq
where the inequality is strict if $r-1$ does not divide $|E|-1$.
\end{lem}

\proof 
By Theorem~\ref{Wakelin}, 
$F(G,\lambda)$ has a root $1$. 
Write 
\beeq\relabel{le-eq3}
F(G,\lambda)
=(\lambda -1)
F_0(G,\lambda).
\eneq
Let $1, -a_1$ and $a_2$ be the three leading coefficients 
of $F_0(G,\lambda)$.
By (\ref{leading-coe-eq}), 
%Since $G$ is a connected graph with no  2-edge-cut, $F(G,\lambda)$ is a polynomial of order $r$, where $r=|E|-|V|+1$, with leading coefficients $1, -|E|$ and ${|E|\choose 2}-\gamma$.Thus 
\beeq\relabel{le-eq4}
a_1+1=|E|\quad \mbox{and}\quad
a_2+a_1={|E|\choose 2}-\gamma,
\eneq
and so 
\beeq\relabel{le-eq5}
\gamma={|E|\choose 2}-a_2-(|E|-1).
\eneq
Since all roots of $F(G,\lambda)$ are real,
by Theorem~\ref{Wakelin}, all roots of 
$F(G,\lambda)$ are positive real numbers.
Then, by Lemma~\ref{le0}, we have 
\beeq\relabel{le-eq6}
a_2
\le {r-1\choose 2} 
\left (\frac{|E|-1}{r-1}\right )^2,
\eneq
where the equality holds if and only if 
$F_0(G,\lambda)=(\lambda -(|E|-1)/(r-1))^{r-1}$, 
which is impossible if $(|E|-1)/(r-1)$ is not an integer
as every rational root of $F(G,\lambda)$ is integral.
Hence (\ref{le-eq1}) follows from (\ref{le-eq5}) and (\ref{le-eq6}).
\proofend

In Lemma~\ref{main-lem0} and Theorem~\ref{main-lem1}
below, let $G=(V,E)$ be 
a non-separable graph in $\Psi_k\setminus \Psi_{k-1}$, 
where $k\ge 0$ and $\Psi_{-1}=\emptyset$, 
such that
$G$ has no $2$-edge-cut nor proper $3$-edge-cut.

Let $v_i$ be the number of vertices of degree $i$ in $G$,
$r=|E|-|V|+1$ and 
\beeq
\alpha=\sum_{i\ge 3}(i-3)v_i.
\eneq
If $G\ne K_2$, then $\delta(G)\ge 3$
and so $\alpha=2|E|-3|V|$.

\begin{lem}\quad \relabel{main-lem0}
%Let $G=(V,E)$ be any non-separable graph in $\Psi_k\setminus \Psi_{k-1}$, where $k\ge 0$ and $\Psi_{-1}=\emptyset$.Assume that $|V|\ge 3$, $G$ has no $2$-edge-cut nor proper $3$-edge-cutand all flow roots of $G$ are real.
If $|V|\ge 3$,  then the following results hold:
\begin{enumerate}
\item  $r\ge \max\{3,8k-6\}$ and $|V|\ge 2k$;
\item if $k=1$, then $\alpha\ge r-2$;
otherwise, $\alpha\ge r+2k-3$.
%\item  Furthermore, if $\alpha=r+2k-3+c$ for some non-negative integer $c$,  then $r\ge 8k-7+(c+3-2k)^2/(c+1)$. 
\end{enumerate}
\end{lem}

\proof As  $G$ is non-separable and has no $2$-edge-cut, 
we have $v_i=0$ for $i<3$. 
Since $\alpha=2|E|-3|V|$ and $r=|E|-|V|+1$, 
we can then deduce that 
$|V|=2r-2-\alpha$ and $|E|=3r-3-\alpha$.

As $|V|=2r-\alpha-2$, $\alpha\ge 0$ and $|V|\ge 3$, 
we have $r\ge 3$.
Since $G$ is non-separable, has no $2$-edge-cut nor proper $3$-edge-cut,
we have $v_3=\gamma$, where $\gamma$ is the number of 
$3$-edge-cuts of $G$. 
Thus Lemma~\ref{le1} implies that 
\beeq\relabel{main-lem1-eq3}
v_3\ge \frac{(2r-3-\alpha)(3r-4-\alpha)}{2(r-1)},
\eneq
where the inequality is strict if $r-1$ does not divide $3r-4-\alpha$.
Since $G\in \Psi_k\setminus \Psi_{k-1}$, 
we have $v_3=|V|-k$ and so inequality (\ref{main-lem1-eq3}) is equivalent 
to the following one:
\beeq\relabel{main-lem1-eq4}
2r-2-\alpha-k\ge  \frac{(2r-3-\alpha)(3r-4-\alpha)}{2(r-1)}.
\eneq
Inequality (\ref{main-lem1-eq4}) is again equivalent to 
\beeq\relabel{main-lem1-eq5}
(r-1)(r-8k+7)\ge (2\alpha -3r+5)^2.
\eneq
We can show that $(r-1)(r-8k+7)>0$.
By (\ref{main-lem1-eq5}), we need only consider  
the case that $2\alpha -3r+5=0$.
So $3r-4-\alpha=1.5(r-1)$, implying that 
inequalities (\ref{main-lem1-eq3}),
(\ref{main-lem1-eq4}) and (\ref{main-lem1-eq5})
are all strict and thus $(r-1)(r-8k+7)>0$.
As $r\ge 3$,  $r\ge 8k-6$ and so $r\ge \max\{3, 8k-6\}$. 
Inequality (\ref{main-lem1-eq5}) is equivalent to 
\beeq\relabel{main-lem1-eq7}
(\alpha-r-2k+4)(\alpha-2r+2k+1)+4(k-1)^2\le 0.
\eneq
Since $r+2k-4\le 2r-2k-1$, inequality (\ref{main-lem1-eq7})
yields that 
$\alpha\ge r-2$ if $k=1$ and 
$\alpha\ge r+2k-3$ otherwise. 
As $|V|\ge 3$, 
it is clear that $|V|\ge 2k$ when $k\le 1$.
If $k\ge 2$, then inequality (\ref{main-lem1-eq7})
implies that $\alpha\le 2r-2k-2$ and so $|V|\ge 2k$ follows.
\proofend

For any bridgeless graph $H$, 
let $\setr(H)$ be the multiset of real roots of $F(H,\lambda)$ in
$(1,2)$. Let 
\beeq
\omega (H)=
\sum_{u\in \setr(H)} (2-u).
\eneq
So $\omega(H)\ge 0$ where equality holds if and only if 
$\setr(H)=\emptyset$.
For any multiset $A$, 
let $|A|=\sum_{a\in A} n_a$, where $n_a$ is the number of times 
that $a$ appears in $A$. 
So $\omega(H)=0$ if and only if $|A|=0$, i.e., 
$\setr(H)=\emptyset$.
Now we are going to prove the main result of this section.

\begin{theo}\quad \relabel{main-lem1}
If $G\ne K_2$ and all flow roots of $G$ are real, 
then one of the following statements holds:
\begin{enumerate}
\item $G$ is $K_4$;  

\item every flow root of $G$ is in the set $\{1,2\}$; and 

\item $k\ge 3$, $|V|\ge 2k$, 
$\omega(G)\ge |E|-2|V|+1\ge 2k-1$, 
$$
|\setr(G)| 
\ge (2k-1)/(2-\xi_k)
\ge 
\left \{
\begin{array}{ll}
27k/11+12/11, &\mbox{if }3\le k\le 5,\\
27k/11-27/22,
&\mbox{if }k\ge 6;
\end{array}
\right.
$$
and 
$\max\{|V|+8k-7, 2|V|+2k-2\}\le |E|
\le ((|V|+1)\xi_k-3)/(\xi_k-1)
< (32|V|-49)/5$.
\end{enumerate}
\end{theo}

\proof 
If $|V|=1$, then $G$ is the graph with one vertex and one loop, and so  $G$ has one root  only (i.e., $1$).
Then consider the case that $|V|=2$. 
$G$ is the graph with two vertices 
and $|E|$ parallel edges joining these two vertices.
Only when $|E|\le 3$, all flow roots of $G$ are real.
As $G\ne K_2$, we have $2\le |E|\le 3$ and thus (ii) holds.

Now assume that $|V|\ge 3$ and 
both (i) and (ii) are not true.
We first prove two claims below before show that (iii) holds.

\noindent {\bf Claim 1}: 
if $k=0$, then $\alpha\ge r-2$. 

Suppose that $k=0$ and $\alpha\le r-3$. 
Then Lemma~\ref{main-lem0} (ii) implies that $\alpha=r-3$. 
However, as $G\in \Psi_0$,  
$G$ is cubic by the given conditions and so $\alpha=0$ and $|E|=\frac 32 |V|$. 
Thus 
$$
3=r=|E|-|V|+1=\frac 32 |V|-|V|+1, 
$$
implying that $|V|=4$ and $|E|=6$. 
Since $G$ is non-separable and has no 2-edge-cut, 
$G$ has no multiedges and so $G\cong K_4$,
contradicting the assumption that (i) does not hold.

\noindent {\bf Claim 2}: $k\ge 3$
and $\omega(G)\ge |E|-2|V|+1\ge 2k-1$, 
where the inequality is strict if 
$F(G,\lambda)$ has some real roots in $(2,\infty)$.

Let $t=|\setr(G)|$, i.e., $t$ is the number of real roots 
of $F(G,\lambda)$ in the interval $(1,2)$,
where the repeated roots are also counted.  
Since $G$ is non-separable,
$F(G,\lambda)$ has one root equal to $1$ 
by Theorem~\ref{Wakelin}.
As all flow roots of $G$ are real, 
Theorem~\ref{Wakelin} also implies that 
all roots of $F(G,\lambda)$ are in $[1,\infty)$. 
As $|E|$ is the sum of all flow roots of $G$
and $F(G,\lambda)$ has exactly $r$ roots, 
one of which is $1$, exactly $t$ of which are in $(1,2)$
and $(r-t-1)$ are at least $2$,  
we have 
\beeq\relabel{main-lem1-eq1}
|E|=3r-\alpha-3\ge 
1+2t-\omega(G)+2(r-1-t)
=2r-1-\omega(G),
\eneq
implying that $\omega(G)\ge \alpha+2-r=|E|-2|V|+1$, 
where the inequality is strict if 
$F(G,\lambda)$ has some real roots in $(2,\infty)$.

Assume that $k\le 2$. 
Then $\setr(G)=\emptyset$ by Theorem~\ref{free-interval},
and so $t=0=\omega(G)$.
Since (ii) does not hold for $G$, 
$F(G,\lambda)$ has some roots in $(2,\infty)$, 
implying that (\ref{main-lem1-eq1}) is strict and so $\alpha<r-2$,
contradicting Claim 1 and the result of Lemma~\ref{main-lem0} (ii). 

As $k\ge 3$, Lemma~\ref{main-lem0} (ii) implies that $\alpha+2-r\ge 2k-1$.
So Claim 2 holds.

Now we are going to show that (iii) holds.

By Claim 2, it remains to show that the bounds for $|\setr(G)|$ 
and $|E|$ in (iii) hold.
Note that $\omega(G)\le (2-\xi_k)|\setr(G)|$.
By Claim 2, $\omega(G)\ge 2k-1$ and so 
$|\setr(G)| \ge (2k-1)/(2-\xi_k)$.
It is known that  $\xi_3\ge 1.430$, $\xi_4\ge 1.361$,
$\xi_5\ge 1.317$ 
and $\xi_k\ge 32/27$ for all $k\ge 5$.
Thus we have 
$$
(2k-1)/(2-\xi_k)
\ge 
\left \{
\begin{array}{ll}
27k/11+12/11, &\mbox{if }3\le k\le 5,\\
27k/11-27/22,
&\mbox{if }k\ge 6.
\end{array}
\right.
$$
By Lemma~\ref{main-lem0} (i), we have $r\ge 8k-6$ and so 
$|E|\ge |V|+8k-7$.
By Lemma~\ref{main-lem0} (ii), we have $\alpha\ge r+2k-3$ and so 
$|E|\ge 2|V|+2k-2$.
Thus the lower bound for $|E|$ in (iii) holds.

Since $|V|\ge 2k$ by Lemma~\ref{main-lem0}, 
$G$ has at least $k\ge 3$  vertices of degree $3$ and thus 
$2$ is its flow root\footnote{It is known that $2$ is a flow root of $G$
if and only if $G$ has at least one vertex of odd degree,
because $G$ has a nowhere-zero ${\mathbb Z}_2$-flow if and only if every
vertex of $G$ has even degree.}.
As $F(G,\lambda)$ is a polynomial of order $r$, 
we have $|\setr(G)|\le r-2$. 
Thus
\beeq\label{main-lem1-cor00-eq40}
|E|-|V|-1=r-2\ge |\setr(G)|.
\eneq 
On the other hand, 
\beeq %\label{main-lem1-cor00-eq40}
|\setr(G)|(2-\xi_k)\ge \omega(G)\ge |E|-2|V|+1,
\eneq 
where the last inequality is from Claim 2.
So we have 
\beeq   %\label{main-lem1-cor00-eq40}
(|E|-|V|-1)(2-\xi_k)\ge |E|-2|V|+1.
\eneq 
Then it follows that 
\beeq
|E|\le \frac{(|V|+1)\xi_k-3}{\xi_k-1}< \frac{32|V|-49}{5},
\eneq
where the last inequality follows from the fact that 
$\xi_k> 32/27$.
\proofend

By Theorem~\ref{main-lem1},  
if $|\setr(G)|\le 8$ or $|E|\le |V|+16$, then 
each flow root of $G$ is contained in $\{1,2,3\}$.
In fact, for such a result, 
the condition that $G$ has no $2$-edge-cut nor proper $3$-edge-cut
is not necessary.

\begin{theo}\quad \relabel{main-lem1-cor}
Let $G=(V,E)$ be any  bridgeless graph.
Assume that all roots of $F(G,\lambda)$ are real.
If either $|E|\le |V|+16$ or 
$|\setr(G)|\le 8$,  
then every root of $F(G,\lambda)$ is in $\{1,2,3\}$.
\end{theo}

\proof Note that if $|E|\le |V|+16$, then 
$|E(B)|\le |V(B)|+16$ for every block $B$ of $G$;
and if $|\setr(G)|\le 8$, then $|\setr(B)|\le 8$
for each block $B$ of $G$.
Thus we need only to prove the result for 
all non-separable graphs.

Let $\setz$ be the set of non-separable bridgeless graphs 
$G$ such that all roots of $F(G,\lambda)$ are real,
and either $|E|\le |V|+16$ or 
$|\setr(G)|\le 8$.
Suppose that there is a graph $G\in \setz$ such that 
some flow root of $G$ is not in $\{1,2,3\}$.
We may assume that $|V|$ has the minimum value among all 
such graphs and that $G\in \Psi_k\setminus \Psi_{k-1}$.
We shall complete the proof by showing the following claims.

\inclaim $G$ contains a $2$-edge-cut or 
a proper $3$-edge-cut.

Suppose that the claim is wrong. 
By the assumption on the minimality of $|V|$,
$\delta(G)\ge 3$. 
Then Theorem~\ref{main-lem1} implies that 
$k\ge 3$, $|\setr(G)|\ge 9$ 
and $|E|\ge |V|+8k-8$,
contradicting the given condition. 
Hence the claim holds.

\inclaim all flow roots of $G$ are in $\{1,2,3\}$,
contradicting the assumption on $G$.

By Claim 1, $G$ contains a $2$-edge-cut or 
a proper $3$-edge-cut.
Let $S$ be such an edge-cut. 
By Lemma~\ref{2-edge} or~\ref{3-edge}, 
\beeq\relabel{main-lem1-cor-eq1}
F(G,\lambda)=\frac{F(G_1,\lambda)F(G_2,\lambda)}{\lambda -1}
\quad \mbox{ or }\quad  
F(G,\lambda)=\frac{F(G_1,\lambda)F(G_2,\lambda)}
{(\lambda -1)(\lambda -2)},
\eneq
where $G_1$ and $G_2$ are the 
graphs stated in  Lemma~\ref{2-edge} or~\ref{3-edge}.
By (\ref{main-lem1-cor-eq1}),
as $G$ has real flow roots only, 
$G_i$ has real flow roots only for $i=1,2$;
if $|\setr(G)|\le 8$, 
then $|\setr(G_i)|\le 8$ for $i=1,2$. 

It is obvious that $G_1$ is non-separable,
and $|E(G_1)|-|V(G_1)|\le |E(G)|-|V(G)|$.
Thus, if $|E|\le |V|+16$, then $|E(G_1)|\le |V(G_1)|+16$.
Hence $G_1\in \setz$ and similarly $G_2\in \setz$.

It is clear that $|V(G_i)|<|V|$ for $i=1,2$.
By the assumption on $G$, 
every flow root of $G_i$ is contained in $\{1,2,3\}$ for $i=1,2$.
Hence 
(\ref{main-lem1-cor-eq1}) implies that 
every flow roots of $G$ is contained in $\{1,2,3\}$.

Therefore claim 2 is true and the result holds.
\proofend

Recently, Kung and Royle~\cite{kung} proved 
a very interesting result.  

\begin{theo}
[Kung and Royle~\cite{kung}]
\quad
\relabel{kung-theo}
If $G$ is a bridgeless graph, then 
its flow roots are integral
if and only if $G$ is the dual of a planar chordal graph. 
\end{theo} 

By Theorems~\ref{main-lem1-cor} and~\ref{kung-theo},
we immediately obtain the following result.

\begin{cor}\quad
\relabel{cor2-th1}
Let $G=(V,E)$ be any bridgeless graph
which has only real flow roots.
If $G$ is not the dual of a planar chordal graph, 
then $|E|\ge |V|+17$ and $|\setr(G)|\ge 9$,
i.e., 
$G$ has at least $9$ flow roots in $(1,2)$.
\proofend
\end{cor}

\begin{cor}
\quad
\relabel{cor3-th1}
For a connected planar graph $G=(V,E)$, 
if $G$ has real chromatic roots only and 
$G$ is not chordal, then
$|V|\ge 19$ and $G$ has at least $9$ chromatic roots in $(1,2)$
(counting multiplicity for each root). 
\end{cor}

\proof We have $P(G,\lambda)=\lambda F(G^{* },\lambda)$.
As $G$ is not chordal, by Theorem~\ref{kung-theo}, 
$G^{* }$ has some non-integral real flow roots.
By Corollary~\ref{cor2-th1},
$|E(G^{* })|\ge |V(G^{* })|+17$ 
and $G^{* }$ has at least $9$ flow roots in $(1,2)$,
where the latter implies that 
$G$ has at least $9$ chromatic roots in $(1,2)$.
Notice that 
$$
|E(G^{* })|=|E(G)|
\quad \mbox{and}\quad
|V(G^{* })|=|E(G)|-|V(G)|+2,
$$
thus $|E(G^{* })|\ge |V(G^{* })|+17$
implies that $|V(G)|\ge 19$. 
\proofend

\def \setl {{\cal L}}

We would like to propose the following conjecture
to end this article.

\begin{con}\quad \label{con2}
For any bridgeless graph $G$, 
if all flow roots of $G$ are real, 
then all flow roots of $G$ are contained in 
$\{1,2,3\}$.
\end{con}

Let $\setl$ be the family of non-separable graphs 
which have no $2$-edge-cut nor $3$-edge-cut.
Lemmas~\ref{block-factor}, \ref{2-edge} 
and~\ref{3-edge} imply that 
 Conjecture~\ref{con2} holds 
if and only if it holds for all graphs in $\setl$.
If Conjecture~\ref{con2} does not hold  for 
some graph $G\in \setl$, 
then Theorem~\ref{main-lem1} implies that 
$G$ has at least $27k/11-27/22$ flow roots in $(1,2)$,
where $k=|W(G)|\ge 3$.
This is a reason why this conjecture is proposed. 

\vspace{0.5 cm}

\noindent {\bf Acknowledgement.} 
The author wishes to thank the referees 
for their very helpful comments and suggestions.

\end{document}